\documentclass[preprint,3p,times]{elsarticle}

\usepackage{amssymb}
\usepackage{amsthm}

\usepackage{lineno}

\usepackage{amsmath, amsfonts}%, amsthm, amssymb}
\usepackage{mathabx}   % upup-..-dowdow arrows
\usepackage[makeroom]{cancel}    % cancellation symbol
\usepackage{tipa}      % fonetica
\usepackage{textcomp}  % caratteri extra
\usepackage{latexsym}  % simboli vari
\usepackage{mathrsfs}  % caratteri corsivi (agg. io); per utilizzarlo: \mathscr{•}

\usepackage{graphicx}
\usepackage{tikz}
\usepackage{pgfplots}
\pgfplotsset{compat=1.6}%{width=10cm, compat=1.9} % retrocompatibilità e proprietà varie plot
\usepgfplotslibrary{groupplots}
\usepackage{pgf}
\usetikzlibrary{positioning,shapes,arrows,fit,calc,automata,perspective,3d}
\usepackage{neuralnetwork}
\usepackage{tikz-network}

\usepackage{hyperref}
\hypersetup{
	colorlinks=true,       % false: boxed links; true: colored links
	linkcolor=blue,%red,          % color of internal links (change box color with linkbordercolor)
	citecolor=red,%green,        % color of links to bibliography
	filecolor=magenta,      % color of file links
	urlcolor=cyan           % color of external links
}
\usepackage{cleveref}

\usepackage{hhline}
\usepackage{multirow}
\usepackage{colortbl} 
\usepackage{algpseudocode}
\usepackage{algorithm}
\usepackage{subcaption} 
\usepackage{enumitem}% http://ctan.org/pkg/enumitem

\usepackage{xspace}

\theoremstyle{definition}
\newtheorem{remark}{Remark}[section]

\newcommand{\cyan}[1]{{\color{cyan}{#1}}}
\newcommand{\mgnt}[1]{{\color{magenta}{#1}}}
\newcommand{\green}[1]{{\color{green}{#1}}}
\newcommand{\orng}[1]{{\color{orange}{#1}}}

\newcommand{\Norm}[1]{\left\Vert{#1}\right\Vert}   % Norma altezza adattabile

\newcommand{\pe}[1]{^{(#1)}}	% exponent within parenthesis

\newcommand{\N}{\ensuremath{\mathbb{N}}}    %   N:  insieme dei numeri naturali
\newcommand{\R}{\ensuremath{\mathbb{R}}}	%   R:  insieme dei numeri Reali
\newcommand{\wh}{\widehat}	%  td:	widetilde

\renewcommand{\v}{\ensuremath{\boldsymbol}}	%   v:	vettore (grassetto)

\newcommand{\bmu}{\v{\mu}}

\newcommand{\eqdot}{\mathrel{\vcentcolon=}}
\newcommand{\dd}{\ensuremath{\mathrm{d}}}
\newcommand{\nn}{\ensuremath{\mathrm{n}}}
\newcommand{\pNN}{$\v{\mu}$BC-NN\xspace}
\newcommand{\pGINN}{$\v{\mu}$BC-GINN\xspace}
\newcommand{\pNNs}{$\v{\mu}$BC-NNs\xspace}
\newcommand{\pFCNNs}{$\v{\mu}$BC-FCNNs\xspace}
\newcommand{\pFCNN}{$\v{\mu}$BC-FCNN\xspace}
\newcommand{\pGINNs}{$\v{\mu}$BC-GINNs\xspace}

\newcommand{\mubcenc}{\ensuremath{\v{\mu}^{\rm (enc)}}}

\newcommand{\MLE}{\ensuremath{\text{me}_{L^2}}}
\newcommand{\MlE}{\ensuremath{\text{me}_{\ell_2}}}
\newcommand{\MHE}{\ensuremath{\text{me}_{H^1}}}
\newcommand{\MLRE}{\ensuremath{\text{mre}_{L^2}}}
\newcommand{\MlRE}{\ensuremath{\text{mre}_{\ell_2}}}
\newcommand{\MHRE}{\ensuremath{\text{mre}_{H^1}}}

\newcommand{\solh}{\ensuremath{\mathsf{u}_h}}
\newcommand{\nnsolh}{\ensuremath{\widehat{\mathsf{u}}_h}}

\journal{arXiv}

\begin{document}

\begin{frontmatter}



\title{Graph-Instructed Neural Networks for parametric problems with varying boundary conditions}


\author[inst1,inst2]{Francesco {Della Santa}\corref{cor1}}
\cortext[cor1]{Corresponding author}

\affiliation[inst2]{organization={Department of Mathematical Sciences, Politecnico di Torino},%Department and Organization
            addressline={Corso Duca degli Abruzzi 24}, 
            % city={Turin},
            postcode={10129}, 
            state={Turin},
            country={Italy}}

\affiliation[inst2]{organization={Gruppo Nazionale per il Calcolo Scientifico INdAM},%Department and Organization
            addressline={Piazzale Aldo Moro 5}, 
            % city={Rome},
            postcode={00185}, 
            state={Rome},
            country={Italy}}

\author[inst1,inst2]{Sandra Pieraccini}

\author[inst1,inst2]{Maria Strazzullo}

\begin{abstract}
This work addresses the accurate and efficient simulation of physical phenomena governed by parametric Partial Differential Equations (PDEs) characterized by varying boundary conditions, where parametric instances modify not only the physics of the problem but also the imposition of boundary constraints on the computational domain.

In such scenarios, classical Galerkin projection-based reduced order techniques encounter a fundamental bottleneck. Parametric boundaries typically necessitate a re-formulation of the discrete problem for each new configuration, and often, these approaches are unsuitable for real-time applications. To overcome these limitations, we propose a novel methodology based on Graph-Instructed Neural Networks (GINNs). The GINN framework effectively learns the mapping between the parametric description of the computational domain and the corresponding PDE solution. Our results demonstrate that the proposed GINN-based models, can efficiently represent highly complex parametric PDEs, serving as a robust and scalable asset for several applied-oriented settings when compared with fully connected architectures. 
\end{abstract}


% \begin{highlights}
% \end{highlights}


\end{frontmatter}

% \linenumbers

\section{Introduction}\label{sec:intro}

Parametric partial differential equations (PDEs) are essential for reliably describing complex physical phenomena in many scientific and industrial contexts, playing an ubiquitous role in applied sciences. The computational costs related to numerical simulations of parametric systems increase for complex problems featuring high variability with respect to the parametric instance. Clearly, the computational burden may be unacceptable when numerical simulations are related to time-consuming activities, such as parameter estimation, uncertainty quantification, or control. 
For this reason, many efforts have been made to conceive numerical methods to efficiently solve parametric PDEs.

An intriguing enhancement to the simulation of parametric PDEs is the efficient treatment of \emph{parametric boundary conditions} (BCs) with the final goal of a fast adaptation to several behaviors of the model without any re-meshing or re-assembling of the system. This contribution focuses on these models, where the parameter defines the portion of the boundary of the computational domain where the type and the value of the BCs might vary.

The formulation directly adapts from \cite{StrazzulloVicini}, where the authors address a parametric changing of the position of Neumann boundary control. We here propose to extend the model to various combinations of BCs of Neumann, Dirichlet, and mixed types. Varying parametric settings can be used to describe many interesting phenomena.

For example, the variation of baffle geometries can increase the performance of heat exchangers. Indeed, modifying the position and the length of the baffles drastically affects the heat conduction \cite{BICER2020106417, YU2019351}.
Another example is the simulation of flows in porous fractured media, which is influenced by the position and the lengths of the intersections among the fractures \cite{BCPS,UQFR}. In addition, groundwater flow models often involve boundary conditions with uncertain locations, such as spring zones or extraction wells \cite{Rubin2003,Tartakovsky2007}. Varying boundary models are crucial in this setting.
Moreover, in aerodynamics, variable actions of deflectors and flaps integrated into turbine blade profiles are commonly used to control loads and improve the energetic performance \cite{Hansen2015,Barlas2010}. 
These models can also find useful applications in urban and building managing, where the position of the openings directly affects the dispersion of pollutants and/or heat \cite{Blocken2015,Chen2009}.

When the parameter represents different physical or geometrical affine deformations, reduced order methods (ROMs)  \cite{hesthaven2015certified,QuarteroniReducedBasisMethods2016} have been successfully used in many settings, such as medical, industrial, and environmental applications (for an overview on ROMs and their applications, see, e.g., \cite{Pichi2025,bookmathlab,Saluzzi2025,StrazzulloZainib} and the references therein). ROMs build a low-dimensional framework based on properly chosen parametric instances of a reference model to solve each new parametric instance in a fast and reliable manner.
However, ROMs might fail and are not straightforwardly applicable in these varying boundary PDEs, being inefficient and inaccurate with respect to the reference solutions.

In this work, we propose a novel NN-based strategy that addresses, in real-time, parametric PDE problems with varying boundary PDEs on a domain with fixed geometry. The method consists of building a NN model that predicts the values of the PDE solution at the nodes of a fixed mesh, given any configuration of BCs and physical parameters belonging to an admissible set on which the NN has been trained; we refer to such models as \emph{Parametric Boundary Conditions Neural Networks} (\pNNs). There are no specific restrictions on the NN architecture, except for the input and output layers; indeed, the input layer must be able to ``read'' all the parameters and BCs information, while the output layer must return the solution values at each mesh node. Due to the strict connection of the method with the definition of a mesh on the domain, we adopt the Graph-Instructed Neural Network (GINN) architecture for defining our models (see \cite{GINN,DELLASANTA2025ewginn,sparseGIlayers}), denoted as $\boldsymbol{\mu}$BC-GINN.
Deep Learning (DL) and Neural Network (NN) models prove to be effective instruments for building surrogate models of highly complex, parametrized problems. These NN-based approaches can be of many different types, ranging from Physics-Informed models that exploit the differential equations of the problem \cite{PINN1,PINN2,PINN3,PINN5,PINN6}, to Neural Operators \cite{NO1,NO2,NO3}, to data-driven approaches \cite{PDENN16,Berrone2021,Millevoi2025,AGLIETTI2026108681,AGLIETTI2026137435}.

However, for ML-guided techniques in the context of varying boundary conditions, the investigation is still limited. For example, we refer to \cite{mousavi2026imposingboundaryconditionsneural} for the application of Neural Operators for data collecting different BCs. We strongly differentiate from the previous literature in several important aspects:
\begin{itemize}
    \item the proposed model of varying boundary PDEs in a parametric setting, which describes not only the action of different BCs acting on the spatial domain and the changing of the physics of the problem;
    \item the integration of parameter variation in the GINN architecture, building a scalable and robust surrogate real-time model for different BCs actions based on sparse local information, linked to the mesh discretization. 
\end{itemize}

The effectiveness of the proposed approach is validated through several numerical tests and compared with non-GINN, more traditional neural NN models based on 1D-Convolutional layers and Fully Connected (FC) layers. In particular, these comparisons highlight the advantages of using GINN architectures for this problem, also on a reduced number of training simulations. Indeed, embedding the mesh structure into the model through Graph-Instructed layers, together with the sparsity of the corresponding weight matrices, allows the construction of deeper models with improved performance compared to other models with a comparable number of parameters but based on 1D-Convolutional and FC layers.

The paper is outlined as follows. 
In Section~\ref{sec:prob}, we present the varying boundary PDEs at the continuous and discrete levels. Section~\ref{sec:ROMsota} is devoted to a literature review on ROMs and surrogate models for PDEs, highlighting the limitations of current approaches in the proposed setting. 
In Section~\ref{sec:GINNgen} we introduce the usage of NNs and their extension to the treated parametric setting, with a focus on GINNs. In Section \ref{sec:results}, we test the robustness and the effectiveness of the $\boldsymbol \mu$BC-GINN though several tests: a linear PDE with parametric varying Neumann boundaries, a linear PDE with parametric varying mixed boundaries, and a nonlinear PDE for parametric varying mixed boundaries. Conclusions follow in Section \ref{sec:conc}.
%%%%%%%%%%%%%%%%%%%%%%%%%%%%%%%%%%%%%%%%%%%%%%%%%%%%%%%%%%%%%%%%%%%%%%%%%%%%%%%%%%%%%%%%%%%

\section{Problem Formulation}\label{sec:formulations}

This section provides the continuous and discrete formulation of the problem with varying boundary conditions and their main features.

\subsection{Continuous problem}\label{sec:prob}

Let us consider a parametric PDE $G: \mathbb U \times \mathcal P \rightarrow \mathbb U^*$, for a suitable Hilbert space $\mathbb U$ and a proper parametric space $\mathcal P$. We denote the solution of the PDE as $u$.
The general problem formulation reads: given a parameter $\bmu \in \mathcal P$, and a forcing term $f \in \mathbb U^*$, find $u = u(\bmu) \in \mathbb U$ such that
\begin{equation}
    \label{eq:strong_pde}
        G(u,\bmu) = f(\bmu) \quad \text{in } \mathbb U^*,
\end{equation}
for $G(u,\bmu) \in \mathbb U^*$, possibly nonlinear in $u$, equipped with boundary conditions. Problem \eqref{eq:strong_pde} can be interpreted in a variational formulation as: for $\bmu \in \mathcal P$, find $u \in \mathbb U$ such that
\begin{equation}\label{eq:weak_pde}
        g(u,v\,; \bmu) - f(v\,; \bmu) \eqdot 
            \left\langle G(u,\bmu), v \right\rangle - \left\langle F(\bmu), v \right\rangle = 0 \quad \text{for all } v \in \mathbb{U},
\end{equation}
with $\left \langle \cdot, \cdot \right \rangle$ denoting the duality pairing between $\mathbb U^*$ and $\mathbb U$, and $F(\bmu) \in \mathbb U^*$ comprises the contributions of the forcing term and the boundary conditions. In particular, defining $\Gamma^{\bmu_b} \subseteq \partial \Omega$ the parameter varying boundary, we consider a parameter of the form $\bmu = (\bmu_{\phi}, \mu_b, \mu_v) \in \mathcal{P} = \mathbb R^{p_{\phi}} \times \mathcal{F}_{pwc}(\Gamma^{\bmu_b}) \times \mathcal{F}_{pwc}(\Gamma^{\bmu_b})$, with $p_{\phi} \geq 1$, and $\mathcal{F}_{pwc}(\Gamma^{\bmu_b}) $ being the space of piece-wise constant functions defined on $\Gamma^{\bmu_b}$. Specifically, $\bmu_{\phi}$ acts on the physics of the problem, $\mu_b$ is a piece-wise constant function that determines \emph{where and what kind of} boundary conditions are applied (e.g., value 0 for Dirichlet, value 1 for Neumann), and $\mu_v$ assigns the value of the boundary condition, say $z_N(\bmu)$ and $z_D(\bmu)$ for Neumann and Dirichlet conditions, respectively.
The whole boundary is the union of the portion affected by the action of $\mu_b$ and $\mu_v$, and a possibly empty portion with fixed boundary conditions $\Gamma_{\text{fix}}$, i.e., $\partial \Omega = \Gamma^{\bmu_b} \cup \Gamma_{\text{fix}}$. 

In other words, $\Gamma_{\text{fix}}$ is a boundary region, of Neumann or Dirichlet type, that is fixed and does not change the nature and properties with respect to $\mu_b$ and $\mu_v$, with boundary value $z_{\text{fix}}(\bmu)$, which can only depend on the physical parameter $\bmu_{\phi}$. The parametric varying boundary $\Gamma^{\bmu_b}$ is defined by two boundary regions, $\Gamma_N^{\bmu_b}$ and $\Gamma_D^{\bmu_b}$, which may not be connected, featuring Neumann and Dirichlet boundary conditions, respectively. The presented formulation admits that or $\Gamma_N^{\bmu_b}$ or $\Gamma_D^{\bmu_b}$ can be possibly empty. The description of $\Gamma^{\bmu_b}$ is encoded in $\mu_b$. Clearly, $\Gamma_{\text{fix}} \cap \Gamma_N^{\bmu_b} \cap \Gamma_D^{\bmu_b} = \emptyset$.
From this definition, we specify the forcing term as
\begin{equation}
\left \langle F(\bmu), v \right \rangle
= 
\left \langle f(\bmu), v \right \rangle + \left \langle f_N(\bmu), v \right \rangle + \left \langle f_D(\bmu), v \right \rangle  + \left \langle f_{\text{fix}}(\bmu), v \right \rangle,
\end{equation}
where 
\begin{equation}\label{eq:fmu}
    \begin{aligned}
        \left\langle f_N(\bmu), v \right\rangle
        & = \int_{\Gamma_N^{\bmu_b}} z_N(\bmu) v \ \text{ds}, \\ 
        \left\langle f_D(\bmu), v \right\rangle & = \int_{\Gamma_D^{\mu_b}}z_D(\bmu) v \ \text{ds} \\
        \left\langle f_{\text{fix}}(\bmu), v \right\rangle & = \int_{\Gamma_{\text{fix}}}
        z_{\text{fix}}(\bmu) v \ \text{ds}.
    \end{aligned}
\end{equation}

\begin{figure}[htb!]
\centering
\resizebox{0.85\textwidth}{!}{
% \begin{tikzpicture}
% \node at (0,0) {$\Omega$};%
% \node at (-1.7,0.9) {$\Gamma_D^{\bmu_b^*}$};%
% \node at (-1.8,-0.9) {$\Gamma_N^{\bmu_b^*}$};%
% \begin{scope}%
% \clip (-2cm,0) rectangle (0,2cm);
%     \draw  (0,0) ellipse(2cm and 1cm);
% \end{scope}%
% \begin{scope}%
% % \fill[green] (2cm,0) rectangle (0,-2cm);
% \clip (2cm,0) rectangle (0,2cm);
%     \draw[densely dotted] (0,0) ellipse(2cm and 1cm);
% \end{scope}%
% \begin{scope}%
% \clip (2cm,0) rectangle (0,-2cm);
%     \draw[densely dotted] (0,0) ellipse(2cm and 1cm);
% \end{scope}%
% \begin{scope}%
% \clip (-2cm,0) rectangle (0,-2cm);
%     \draw[densely dotted] (0,0) ellipse(2cm and 1cm);
% \end{scope}%
% \begin{scope}%
% \clip (2cm,0) rectangle (6,-2cm);
%     \draw[densely dotted] (6,0) ellipse(2cm and 1cm);
% \end{scope}%
% \begin{scope}%
% \clip (6,0) rectangle ++(2,1);
%     \draw[densely dotted] (6,0) ellipse(2cm and 1cm);
% \end{scope}%
% \begin{scope}%
% \clip (6,0) rectangle ++(2,-1);
%     \draw (6,0) ellipse(2cm and 1cm);
% \end{scope}%
% \begin{scope}%
% \clip  (6,0) rectangle ++(-2,1);
%     \draw (6,0) ellipse(2cm and 1cm);
% \end{scope}%
% \node at (6,0) {$\Omega$};%
% \node at (4.3,0.9) {$\Gamma_D^{\overline{\bmu_b}}$};%
% \node at (4.2,-0.9) {$\Gamma_N^{\overline{\bmu_b}}$};%
% \end{tikzpicture}

\begin{tikzpicture}
% Ellisse di sinistra
\node at (0,0) {$\Omega$};%
\node at (-1.7,0.9) {$\Gamma_D^{\bmu_b^*}$};%
\node at (-1.8,-0.9) {$\Gamma_N^{\bmu_b^*}$};%
\begin{scope}
\clip (-2cm,0) rectangle (0,2cm);
    \draw (0,0) ellipse(2cm and 1cm);
\end{scope}
\begin{scope}
\clip (2cm,0) rectangle (0,2cm);
    \draw[densely dotted] (0,0) ellipse(2cm and 1cm);
\end{scope}
\begin{scope}
\clip (2cm,0) rectangle (0,-2cm);
    \draw[densely dotted] (0,0) ellipse(2cm and 1cm);
\end{scope}
\begin{scope}
\clip (-2cm,0) rectangle (0,-2cm);
    \draw[densely dotted] (0,0) ellipse(2cm and 1cm);
\end{scope}

% Ellisse di destra
\begin{scope}
\clip (2cm,0) rectangle (6,-2cm);
    \draw[densely dotted] (6,0) ellipse(2cm and 1cm);
\end{scope}
\begin{scope}
\clip (6,0) rectangle ++(2,1);
    \draw[densely dotted] (6,0) ellipse(2cm and 1cm);
\end{scope}
\begin{scope}
\clip (6,0) rectangle ++(2,-1);
    \draw (6,0) ellipse(2cm and 1cm);
\end{scope}
\begin{scope}
\clip (6,0) rectangle ++(-2,1);
    \draw (6,0) ellipse(2cm and 1cm);
\end{scope}
\node at (6,0) {$\Omega$};%
\node at (4.3,0.9) {$\Gamma_D^{\overline{\bmu_b}}$};%
\node at (4.2,-0.9) {$\Gamma_N^{\overline{\bmu_b}}$};%

% Aggiunta porzione dashed su entrambi i bordi esistenti

% Su ellisse sinistra
\begin{scope}
\clip (0,0) ++(0.1,0.3) rectangle ++(2,1);
\draw[dashed] (0,0) ellipse(2cm and 1cm);
\end{scope}
\node at (1.9,0.9) {$\Gamma_{\text{fix}}$};

% Su ellisse destra
\begin{scope}
\clip (6,0) ++(0.1,0.3) rectangle ++(2,1);
\draw[dashed] (6,0) ellipse(2cm and 1cm);
\end{scope}
\node at (7.9,0.9) {$\Gamma_{\text{fix}}$};

\end{tikzpicture}      
}
\caption{Representation of a general domain $\Omega$ for two values of $\mu_b$, i.e., $\mu_b = \mu_b^*$ (left) and $\mu_b = \overline{\mu_b}$ (right). The solid, the dotted, and the dashed-dotted lines represent $\Gamma_D^{\bmu_b}, \Gamma_N^{\bmu_b}$, and $\Gamma_{\text{fix}}$, respectively.}
\label{fig:genericdomain}
\end{figure}
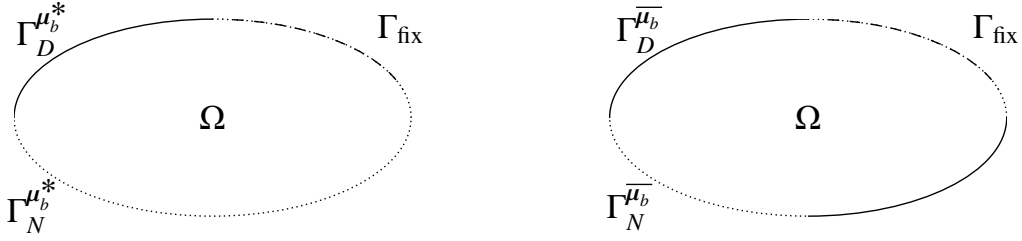

Figure \ref{fig:genericdomain} schematically represents a generic spatial domain $\Omega$ with the boundary $\partial \Omega$ under the action of two different parameters $\bmu^* = (\bmu_{\phi}^*, \mu_b^*, \mu_v^*)$ and $\overline{\bmu} = (\overline{\bmu_{\phi}}, \overline{\mu_b}, \overline{\mu_v})$.   

\begin{remark}
    We remark that the continuous formulation reflects the problems treated in the numerical experiments. As a consequence, the discrete formulation, described in the following Section \ref{sec:discrete_form}, also reflects this characteristic. However, the proposed strategies can be straightforwardly extended to other complex scenarios, where the physical parameter might be non-constant in space, as well as $z_N(\bmu)$ and $z_D(\bmu)$ be more involved depending on the quality of the data or the discetization used.
\end{remark}

\subsection{Discrete formulation}\label{sec:discrete_form} 

We now focus on the algebraic discrete formulation of the weak problem \eqref{eq:weak_pde} and of the discrete representation of parametric space $\mathcal P$, which we denote with $\mathcal P_h$, which will be detailed later for the sake of presentation. Given $\bmu \in \mathcal P_h$ and a tesselation, i.e., a mesh of the physical domain $\Omega$ with $N_h$ nodes, and considering $N_{\text{dof}}$ degrees of freedoms (dofs), one wants to find the approximated \emph{solution vector} $\mathsf u = \mathsf u(\bmu) \in \mathbb{R}^{N_{\text{dof}}}$
 such that
\begin{equation}
    \label{eq:res_vec}
        \mathcal G(\mathsf u; \bmu) = \mathcal F (\bmu),
\end{equation}
where $\mathcal G: \mathbb R^{N_{\text{dof}}} \times \mathcal P \rightarrow \mathbb R^{N_{\text{dof}}}$ and $\mathcal F( \boldsymbol \bmu) \in \mathbb R^{N_{\text{dof}}}$ are the discrete \emph{parametric system} and the \emph{parametric forcing term} of the system considered. In this setting, $N_{\text{dof}}$ and $N_h$ might not coincide.  We remark that $\mathsf u(\bmu) = [\mathsf u^1(\bmu), \cdots, \mathsf u^{N_{\text{dof}}}(\bmu)]^\top$ represents the value of the discrete solution for each dof. Moreover, we use $\mathsf u_h(\bmu) = [\mathsf u_h^1(\bmu), \cdots, \mathsf u_h^{N_h}(\bmu)]^\top$  to denote the value of the solution in the nodes of the mesh, collected in the set $\mathcal M$. 
Since $\bmu$BC-NNs do not depend on the numerical discretization strategy employed to solve \eqref{eq:res_vec}, we here skip the details concerning the numerical approximation, and we postpone them to the numerical investigation section.

\begin{figure}[htb!]
\centering
\resizebox{0.55\textwidth}{!}{
% (nel preambolo)
% \usetikzlibrary{calc}

\begin{tikzpicture}[scale=3.0]

%-------------------
% COMANDI MARCATORI
%-------------------
\newcommand{\sqgray}[2]{%
  \filldraw[black, fill=gray!50]
  ($(#1,#2)+(-0.015,-0.015)$) rectangle ++(0.04,0.04);
}
\newcommand{\circblack}[2]{%
  \filldraw[black, fill=black] (#1,#2) circle (0.02);
}
\newcommand{\sqwhite}[2]{%
  \filldraw[black, fill=white]
  ($(#1,#2)+(-0.015,-0.015)$) rectangle ++(0.04,0.04);
}

%-------------------
% COMANDO PER X
%-------------------
\newcommand{\xmark}[2]{\node[text=black] at (#1,#2) {\Large $\times$};}

%-------------------
% GRIGLIA INTERNA
%-------------------
\draw[black] (0,0) -- (0,1);
\draw[black, dashed] (0.2,0) -- (0.2,1);
\draw[black, dashed] (0.4,0) -- (0.4,1);
\draw[black, dashed] (0.6,0) -- (0.6,1);
\draw[black, dashed] (0.8,0) -- (0.8,1);

\draw[black, dashed] (0,0.2) -- (1,0.2);
\draw[black, dashed] (0,0.4) -- (1,0.4);
\draw[black, dashed] (0,0.6) -- (1,0.6);
\draw[black, dashed] (0,0.8) -- (1,0.8);

%-------------------
% BORDI GRIGI (zone con label μ)
%-------------------
\draw[gray!40, line width=3pt] (0.2,0) -- (1,0);   % μ_b^(1) - μ_b^(5)
\draw[gray!40, line width=3pt] (1,0) -- (1,0.4);   % μ_b^(5) - μ_b^(7)

%-------------------
% BORDI NERI SOTTILI (zone con X)
%-------------------
% lato sinistro (tutto X)
\draw[black] (0,0) -- (0,1);

% lato alto (tutto X)
\draw[black] (0,1) -- (1,1);

% lato basso (solo tratto con X)
\draw[black] (0,0) -- (0.2,0);

% lato destro (solo parte con X)
\draw[black] (1,0.4) -- (1,1);

%-------------------
% X (ex cerchi neri sul bordo sinistro/alto + alcuni a destra)
%-------------------
\xmark{1}{1}

\xmark{0}{0}
\xmark{0}{0.2}
\xmark{0}{0.4}
\xmark{0}{0.6}
\xmark{0}{0.8}
\xmark{0}{1}

\xmark{1}{0.6}
\xmark{1}{0.8}

\xmark{0.2}{1}
\xmark{0.4}{1}
\xmark{0.6}{1}
\xmark{0.8}{1}

%-------------------
% MARCATORI ASSOCIATI ALLE μ_b^(i)
%-------------------
% μ_b^(1), μ_b^(2): quadrati grigi
\sqgray{0.2}{0}
\sqgray{0.4}{0}

% μ_b^(3), μ_b^(4), μ_b^(5): cerchi neri
\circblack{0.6}{0}
\circblack{0.8}{0}
\circblack{1}{0}

% μ_b^(6), μ_b^(7): quadrati bianchi (sul bordo destro)
\sqwhite{1}{0.2}
\sqwhite{1}{0.4}

%-------------------
% LABEL
%-------------------
\node at (0.2,-0.09) {\tiny $\mu_b^{(1)}$};
\node at (0.4,-0.09) {\tiny $\mu_b^{(2)}$};
\node at (0.6,-0.09) {\tiny $\mu_b^{(3)}$};
\node at (0.8,-0.09) {\tiny $\mu_b^{(4)}$};
\node at (1.0,-0.09) {\tiny $\mu_b^{(5)}$};

\node at (1.15,0.2) {\tiny $\mu_b^{(6)}$};
\node at (1.15,0.4) {\tiny $\mu_b^{(7)}$};

\node at (1.4,0.2) {$\Gamma^{\bmu_b}$};
\node at (-0.3,0.6) {$\Gamma_{\text{fix}}$};

\end{tikzpicture}
}
\caption{Schematic representation of a mesh discretization on a square domain $\Omega$ with a portion of the boundary under the action of the parameters $\bmu_b$ and $\bmu_v$. This portion is $\Gamma^{\bmu_b}$, the continuous, grey, and thick line. The nodes that may change features are denoted with square and circle dots, while the nodes of the boundary with fixed features are denoted with black crosses. the color of the square and circle dotes denote the value of the BCs described by $\bmu_v$.}
\label{fig:mub}
\end{figure}
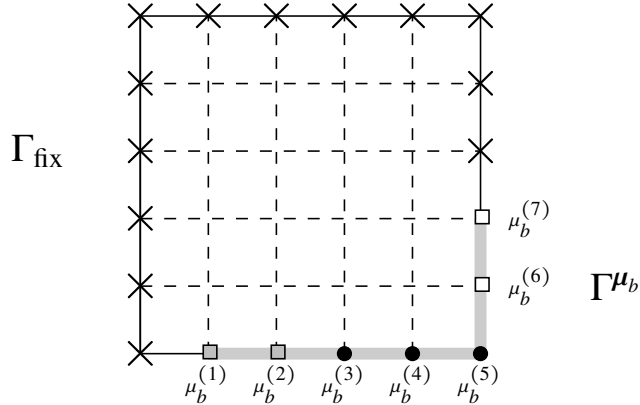

The major aspect of the discretization is to provide a suitable discrete interpretation of $\bmu \in \mathcal P$ for the surrogate model. Indeed, at the discrete level, we can define $\bmu = (\bmu_{\phi}, \bmu_b, \bmu_v)$ in $\mathcal P_h \subset \R^P$, with $P = p_{\phi} + p_b + p_v$, consisting of:
\begin{itemize}
    \item $\bmu_{\phi} \in \mathbb R^{p_\phi}$, where $p_{\phi}$ is the number of physical parameters of the problem;
    \item $\bmu_b \in \mathbb R^{p_b}$, where $p_b$ is the number of mesh nodes in $\Gamma^{\bmu_b}$. We recall that this parameter indicates where the boundary conditions change and the type of boundary condition, Dirichlet or Neumann. We will give more details about the encoding strategy in Section \ref{sec:results}.
    \item $\bmu_v \in \mathbb R^{p_v}$, which indicates the value of the boundary condition related to the nodes lying on $\Gamma^{\bmu_b}$. {In particular, $p_v$ depends on the number $n_v$ of values in $\R^k$ that the function $\mu_v$ assumes on $\Gamma^{\bmu_b}$, where $k>1$ in case of vectorial PDEs with BCs characterized by vectors in $\R^k$; for example, with scalar BCs, if $\Gamma^{\bmu_b}$ is partitioned in one homogeneous Dirichlet BC, one homogeneous Neumann BC, and one non-homogeneous Neumann BC, we have $p_v=n_v \cdot k = 3\cdot 1 = 3$.}
\end{itemize} 
For the sake of clarity, in Figure \ref{fig:mub}, we illustrate an example for a mesh on a square domain. The black crosses denote the nodes of $\mathcal M$ located in $\Gamma_{\text{fix}}$ (the continuous black line), where the boundary conditions are prescribed, while the circle and square dots denote the nodes of $\mathcal M$ that may vary the boundary conditions according to the parameters $\bmu_b = (\mu_b\pe{1}, \ldots,\mu_b\pe{7}) \in \mathbb R^{7}$ and $\bmu_v$; i.e., they are the nodes on $\Gamma^{\bmu_b}$ (continuous, grey, thick line). For the example illustrated in Figure \ref{fig:mub}, squares denote Neumann BCs, while circles denote Dirichlet BCs; the color of the dots represents the BC (scalar) value contained in $\bmu_v\in\R^3$.

In the next section, we propose a literature overview of intrusive and non-intrusive strategies to solve parametric PDEs, providing the core motivation of the employment of $\bmu$BC-NNs in this varying boundary setting.

%%%%%%%%%%%%%%%%%%%%%%%%%%%%%%%%%%%%%%%%%%%%%%%%%%%%%%%%%%%%%%%%%%%%%%%%%%%%%%%%%%%%%%%%%%%

\section{State of the art for efficiently solving parametric PDEs}\label{sec:ROMsota}

In several fields based on PDEs, standard discretization techniques can result in high computational costs due to the large number of degrees of freedom $N_{\text{dof}}$, and may be unsuitable for real-time and many-query applications. 
This section proposes an overview on the numerical surrogate models proposed to address the specific issue of efficient multiple parametric evaluation.

\subsection{Surrogate Models based on ROMs} 

Classically, ROMs have being conceived to build a linear low-dimensional framework of dimension $N \ll N_{\text{dofs}}$ based on proper chosen parametric solutions to provide a surrogate reduced representation of $\mathsf u_h$, say $\mathsf u_N = \mathsf u_N(\bmu)$, with $\mathsf u \approx \mathsf V \mathsf u_N$, where $\mathsf V \in \mathbb R^{N_{\text{dofs}} \times N}$ collects, column-wise, the basis functions of the reduced space. We stress that $\bmu$ has always been interpreted as a physical or geometrical parameter for fixed boundaries in a compact subset $\mathcal D$ of $\mathbb R^P$, with $P \geq 1$. The ROM strategy is based on an efficient separation into:
\begin{itemize}
\item an \emph{offline phase}, that builds the matrix $\mathsf V$ and assembles quantities that are not related to $\bmu$.
\item an \emph{online stage}, where the reduced system is assembled and solved for each new parametric instance $\bmu^*$, finding $\mathsf u_N(\bmu^*)$ efficiently, exploiting the pre-computed offline quantities.
\end{itemize}
This \emph{offline-online} decomposition is based on the affinity assumption. Namely, for the parameter $\bmu \in \mathcal D$, the system \eqref{eq:weak_pde} must be represented as 
\begin{equation}
\label{eq:affine}
    g(u,v; \bmu) - f(v; \bmu) = \sum_{i=1}^{Q_g} \theta_i(\bmu) g_i(u, v) -\sum_{i=1}^{Q_f} \theta_i(\bmu) f_i(v) = 0,
\end{equation}
for $Q_g$ and $Q_f$ in $\mathbb N$, with $\theta_i: \mathcal D \rightarrow \mathbb R$ smooth real functions and $g_i(\cdot, \cdot)$ and $f_i(\cdot)$ forms not depending on the parameter, comprising forcing terms and boundary conditions action. Assumption \eqref{eq:affine} is crucial, but it is not verified in our case. Indeed, the forcing terms in \eqref{eq:fmu} depend on $\mu_b$ and $\mu_v$, thus, need to be re-assembled for each new parametric instance.

Many strategies can be employed to tackle this issue. For example, the Empirical Interpolation Method (EIM) \cite{barrault2004eim}, can be used to recover the affine decomposition of the system, allowing the application of standard ROMs techniques such as Proper Orthogonal Decomposition (POD) and Greedy \cite{hesthaven2015certified}. Another approach builds the reduced vector coefficient $\mathsf u_N$ not relying on a linear reduced system, but non-intrusively, as done in POD with interpolation \cite{Bui-Thanh2003,Demo2019873} or by POD-NN \cite{HESTHAVEN201855}. 
Despite the non-intrusiveness, these approaches still rely on linear compressions and approximations. This might represent a bottleneck for problems with low reducibility \cite{GREIF2019216}. For this reason, nonlinear ROMs have been extensively and successfully employed. They are based on a paradigm change in the reduced representation given by $\mathsf u \approx \psi (\mathsf u_N)$, with $\psi: \mathbb R^{N} \rightarrow \mathbb R^{N_{\text{dofs}}}$ being a nonlinear map. Among nonlinear ROMs, we mention local POD \cite{Amsallem2012891}, kernel POD \cite{Díez20217306}, registration methods \cite{Taddei2020A997}, shifted POD \cite{ReissShift}, and nonlinear ROMs based on autoencoders have emerged and have been applied in several fields, see, e.g., \cite{Franco2023,Fresca2021,Hernandez2021,Lee2020,Romor2023,Pichi2024}.

However, in the proposed context, several issues arise in terms of the suitability of ROMs. For example, classical ROMs can be applied only when $\Gamma_D^{\bmu_b}=\emptyset$ and even in this simpler settings the EIM strategies might lead to a large number of interpolation points around $\Gamma_N^{\bmu_b}$ to properly describe the system at hand, yielding inefficient online solutions, as already shown in \cite{StrazzulloVicini} for varying Neumann control problems.  
In the same work, the complex behavior related to different Neumann boundary conditions results in a large reduced space, even once the affine assumption is recovered, and tailored strategies based on local POD must be employed to increase the accuracy of the solution. Another challenge adds to the problem under investigation: even when $\Gamma_D^{\bmu_b}\neq \emptyset$, we deal with a great variability in the parametric boundary conditions. The primary issue is the difficult representation of a globally coherent latent space for the solution varying with respect $\bmu$. We believe that this challenge affects the direct applicability of ROMs based on autoencoder architectures. Similarly, other nonlinear reduction techniques, such as registration-based approaches or shifted POD, rely on the existence of smooth alignment maps between solution fields. However, the nature of such mappings is not evident in the case we address. Finally, we remark that the case of changing the location of Dirichlet boundary conditions changes the number $N_{\text{dofs}}$ from one parameter to another. Thus, it cannot be addressed by any Galerkin projection-based ROM approach and, indeed, lacks an investigation in the ROMs community.

\subsection{Surrogate Models based on ML} Lately, ML-based surrogate models have been extensively exploited in numerical simulations, and NNs have been applied as direct solvers for approximated systems, offering highly efficient solutions in many applications \cite{PDENN16,Berrone2021,Millevoi2025,AGLIETTI2026137435}. 
One of the most interesting innovation in using NN in the PDE solution was to integrate the information of the residual of the PDE in the loss function, as introduced in \cite{PDENN20}, i.e., working with Physics informed Neural Networks (PINNs). While PINNs proved successful in many contexts \cite{PINN1,PINN2,PINN3,PINN5,PINN6}, they suffer in representing solutions with high frequencies \cite{WANG2022110768}, convergence issues due to different loss components \cite{WANG2022110768,DEMO2023383}, and complex dynamics \cite{Stiff}. All these issues are amplified in the extension to parametric problems and do not make PINNs a suitable option for the problem we are dealing with.
Recently, there has been growing interest in Neural Operators (NOs), conceived as an extension of PINN to learn maps between function spaces \cite{JMLR:v24:21-1524}. NOs proved really effective in many contexts, see, e.g., \cite{NO1,NO2,NO3}. Often these architectures are based on convolutional layers and, thus, they strongly rely on structured data. This hypothesis is restrictive in a PDE context, based on simulations on usually unstructured meshes, or for a real-data case scenario.  
For this reason, the application of Graph Neural Networks (GNN) has become more and more popular in the context of PDEs. GNN naturally suits the mesh design of nodes and edges. Recently, many contributions have been made to the field of GNN-based surrogates, see \cite{Franco2023MINNs,GNN1,belbute2020combining}. In the parametric setting, we list these notable contributions \cite{MORRISON2024117458,Franco2023, Pichi2024}. However, the problem of varying boundary poses several challenges, as already highlighted in Section \ref{sec:ROMsota}.

% \medskip
For these reasons, we decided to address the peculiar parametric nature of varying boundary conditions with novel architectures based on Graph-Instructed Neural Networks, which will be detailed in the following section.

%%%%%%%%%%%%%%%%%%%%%%%%%%%%%%%%%%%%%%%%%%%%%%%%%%%%%%%%%%%%%%%%%%%%%%%%%%%%%%%%%%%%%%%%%%%

\section{\pNNs: NNs for Parametric PDEs with Parametric Boundary Conditions}\label{sec:muBCNNs}
\label{sec:GINNgen}

In this section, we describe the NN architectures that we use as surrogate models for predicting a solution $\solh(\v{\mu})$, given the parameters $\v{\mu}\in\mathcal{P}_h$ and the mesh with $N_h\in\N$ nodes for the domain $\Omega$ (see Section \ref{sec:discrete_form}). {A Parametric Boundary Conditions NN} (\pNN) is a NN characterized by an architecture that is able to ``read'' the discrete formulation described in Section \ref{sec:discrete_form} of a parametric PDE problem, and that returns an array $\nnsolh(\v{\mu})\in\R^{N_h\times m}$ such that $\nnsolh(\v{\mu})\approx \solh(\v{\mu})$, where $m \geq 1$ allows for vector valued solutions. While the output encoding is trivial (the output layer just returns a vector if $m=1$, a matrix otherwise), the input encoding must consider all the information that the parameter $\v{\mu}\in\mathcal{P}_h$ implicitly contains about the problem, and how it can vary in $\mathcal{P}_h$. After preliminary investigations, {the input encoding adopted} for a generic problem (see Section \ref{sec:discrete_form}) is based on a characterization of all the mesh nodes with a vector of values $(\v{\beta}\pe{i},{\dd\pe{i},\nn\pe{i}}, \v{\phi}\pe{i})$, for {$i=1, \ldots, N_h$}, such that:
\begin{itemize}
    \item {$\dd\pe{i}\in\{0,1\}$, and $\nn\pe{i}\in\{0,1\}$} describe if the node is a Dirichlet node ($\dd\pe{i}=1$) or a Neumann node ($\nn\pe{i}=1$), respectively. If both are zeros, the node is not a boundary node.
    \item $\v{\beta}\pe{i}\in\R^k$ is the vector of the BC values (Dirichlet or Neumann, depending on {$\dd\pe{i}$ and $\nn\pe{i}$}, values). We have $k>1$ if BC values are denoted by vectors (i.e., when solving vectorial PDEs). By convention, $\v{\beta}\pe{i}=\v{0}$ if {$\dd\pe{i}=\nn\pe{i}=0$} (i.e., for a non-boundary mesh node).
    \item $\v{\phi}\pe{i}\in\R^q$ is the vector of physical parameter values at the $i$-th mesh node. If these parameters are not depending on the mesh node's coordinates in the domain (see Section \ref{sec:formulations}), then $q=p_{\phi}$ and we have $\v{\phi}\pe{i}=\v{\phi}=\bmu_{\phi}$ constant for each $i=1,\ldots ,N_h$. From preliminary analyses and experiments, the redundancy of information helps the NN training, therefore this input encoding is also used for the case of constant $\v{\phi}\pe{i}$ (i.e., repeating $\bmu_{\phi}$ for each node). By convention, if there are no physical parameter values, each vector $\v{\phi}\pe{i}$ becomes a scalar value $\phi=1$.
\end{itemize}

In summary, the input of a \pNN is a matrix
\begin{equation}\label{eq:input_pNN}
    \mubcenc \coloneqq
    \begin{bmatrix}
        \v{\beta}^{(1)\,T} & \dd\pe{1} & \nn\pe{1} & \v{\phi}^{(1)\,T}\\
        \vdots & \vdots & \vdots & \vdots \\
        \v{\beta}^{(N_h)\,T} & \dd\pe{N_h} & \nn\pe{N_h} & \v{\phi}^{(N_h)\,T}\\
    \end{bmatrix}\in\R^{N_h\times (k+2+q)}\,,
\end{equation}
where the order of parameters $\v{\beta}\pe{i}, {\dd\pe{i}, \nn\pe{i}}, \v{\phi}\pe{i}$ illustrated in \eqref{eq:input_pNN} is the one adopted in the implementations for the numerical experiments.

Concerning the hidden layers of a \pNN, there are no specific restrictions. Actually, any sequence of NN layers that is able to transform inputs like \eqref{eq:input_pNN} into arrays $\nnsolh(\v{\mu})\in\R^{N_h\times m}$ can be used. Nonetheless, not all NN layers are able to process inputs like \eqref{eq:input_pNN}; in this work, for doing so, we consider the following types of layers: Graph-Instructed (GI) layers and 1-dimensional Convolutional (Conv1D) layers. These layers share the characteristic of seeing the $N_h$-by-$(k+2+q)$ input matrix as a sequence of $N_h$ vectors in $\R^{k+2+q}$, and they process this sequence according to these principles:
\begin{itemize}
    \item GI layers: process and transform the input sequence of vectors following the adjacency of nodes in the mesh (for more details, see \cite{GINN,DELLASANTA2025ewginn,sparseGIlayers} and Section \ref{sec:GINNs} in the following);
    \item Conv1D layers: process and transform the input sequence of vectors by index proximity, depending on the size of the kernel (for more details, see \cite{Goodfellow-et-al-2016}). Since mesh nodes' indexing $i=1,\ldots ,N_h$ may not reflect proximity in the domain, we use Conv1D layers with kernel size equal to $N_h$ in order to gather all the mesh nodes together.
\end{itemize}

From these two types of layers, we develop two distinct NN architecture archetypes, respectively{{: one based on GI layers and one based on Conv1D and Fully Connected (FC) layers.} Between the two NN architectures, for approximately the same number of parameters, the one based on GI layers can achieve significantly greater depth; i.e., a substantially larger number of layers. This enables a more effective exploitation of deep (and residual) architectures. This is made possible by the sparsity of the mesh adjacency matrix, which allows the construction of GI layers with sparse weight tensors, thereby increasing model depth without increasing the parameter count}. Moreover, numerical experiments show that the GI-based NN archetype is also the one with the best performances (see Section \ref{sec:results}), thanks to the embedding of the mesh's geometrical properties into the NN architecture.

We continue this section by briefly recalling the mechanism of GI layers; then, {we describe the NN architecture archetypes based on GI layers and Conv1D-FC layers}, respectively (Section \ref{sec:architectures}). In the end, we conclude by describing the loss function used for training these NN models (Section \ref{sec:trainingloss}).

% -------------------------------------------------------

\subsection{Graph-Instructed Layers and Graph-Instructed NNs}\label{sec:GINNs}

\begin{figure}[htb]
    \centering
    \resizebox{0.65\textwidth}{!}{
    \begin{tikzpicture}[multilayer=3d,rotate=90]
    \Vertices{Figures/ginnfilter_verts_v2b.csv}
    \Edges{Figures/ginnfilter_edges_v2b.csv}
    \end{tikzpicture}
	}
    \caption{Visual representation of \eqref{eq:ginn_node_action}. Example with $i=1$ and a non-directed graph of four nodes.}
    \label{fig:ginnfilter}
\end{figure}
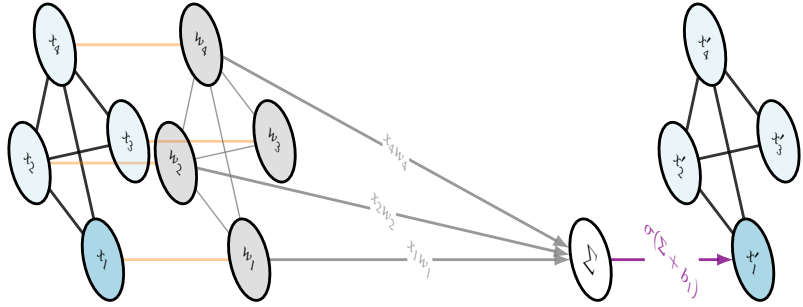

A GINN is a NN characterized by Graph-Instructed (GI) Layers. These layers are defined by an alternative graph-convolution operation introduced for the first time in \cite{GINN}. Briefly, given a graph $G=(V, E)${, with nodes $V$, edges $E$, without self-loops and an associated} input feature $x_i\in\R$ for each graph's node $v_i\in V$, the graph-convolution operation of a GI layer consists in a message-passing procedure where each node $v_i\in V$ sends to its neighbors a message equal to the input feature $x_i$ rescaled by a weight $w_i\in\R$ assigned to the node by the GI layer; then, the output feature corresponding to each node is obtained by summing up all the messages received by the neighbors, adding the bias, and applying the activation function.
Therefore, the message-passing interpretation can be summarized by Figure \ref{fig:ginnfilter} and the following node-wise equation
\begin{equation}\label{eq:ginn_node_action}
x_{i}' = \sigma \Big( \sum_{j \in \mathrm{N}_{\text{in}}(i)\cup \{i\}} x_j \, w_j  + b_i \Big)\,,
\end{equation}
where $x_i'$ is the output feature corresponding to node $v_i$ and $\mathrm{N}_{\rm in}(i)$ is the set of indices such that $j\in\mathrm{N}_{\rm in}(i)$ if and only if $e_{ij}=\{v_i,v_j\}$ is an edge of the graph (i.e., $e_{ij}\in E$). The symbol $\sigma$ denotes the activation function $\sigma:\R\rightarrow\R$.

Assuming to apply \eqref{eq:ginn_node_action} to each node of the graph, we can formally define the GI layers through a compact formulation. Given {the adjacency matrix $A\in\R^{\mathrm{V}\times \mathrm{V}}$ of the graph $G$, $\mathrm{V}=|V|$}, the basic and simplest version of GI layer with respect to $G$ is a NN layer that process one input feature per node and returns one output feature per node, and that is described by a function $\mathcal{L}^{\rm GI}:{\R^{\mathrm{V}}\rightarrow\R^{\mathrm{V}}}$, such that
\begin{equation}\label{eq:GI_action_simple}
    \begin{aligned}
    \mathcal{L}^{\rm GI}(\v{x}) 
    &
    = \v{\sigma}\left( \wh{W}^T \, \v{x} + \v{b}\right) = \\
    &
    = 
    \v{\sigma}\left( \,\left(\,\mathrm{diag}(\v{w}) (A + \mathbb{I}_{{\mathrm{V}}})\,\right)^T\, \v{x} + \v{b}\right),
    \end{aligned}
    \,
\end{equation}
where $\v{x}\in\R^{{\mathrm{V}}}$ denotes the vector of input features and:
\begin{itemize}
    \item {$\v{w}\in\R^{\mathrm{V}}$ is the weight vector, with the component $w_i$ associated to node $v_i\in V$, for each $i=1,\ldots ,\mathrm{V}$. The vector $\v{w}$ is used to build the weight matrix $\wh{W}:=\mathrm{diag}(\v{w}) (A + \mathbb{I}_{\mathrm{V}})$. We remark that the diagonal matrix $\mathrm{diag}(\v{w})$ in the definition of $\wh{W}$ is used only for describing in matrix form the multiplication of the $i$-th row of $(A + \mathbb{I}_{\mathrm{V}})$ by the weight $w_i$, for each $i=1,\ldots ,\mathrm{V}$.}

    \item $\v{\sigma}:\R^{{\mathrm{V}}}\rightarrow\R^{{\mathrm{V}}}$ is the element-wise application of the activation function $\sigma$;
    
    \item $\v{b}\in\R^{{\mathrm{V}}}$ is the bias vector{, such that $b_i$ is the bias associated to node $v_i$, for each $i=1,\ldots ,\mathrm{V}$}.
\end{itemize}

\begin{figure}[htb]
    \centering
    \subcaptionbox{FC representation of GI layer\label{fig:GIa}}[0.45\textwidth]{
        \centering
        \begin{tikzpicture}[x=1cm, y=1cm, >=stealth]

% HIDDEN 3

\node [circle,fill=white!50,minimum size=0.75cm] (hidden3-1) at (0,1.5) {$x_1$};
\node [] (w1) at (0.75,1.75) {{\color{cyan!50}$w_1$}};

\node [circle,fill=white!50,minimum size=0.75cm] (hidden3-2) at (0,0.5) {$x_2$};
\node [] (w2) at (0.7,0.8) {{\color{magenta!50}$w_2$}};

\node [circle,fill=white!50,minimum size=0.75cm] (hidden3-3) at (0,-0.5) {$x_3$};
\node [] (w3) at (0.7,-0.8) {{\color{orange!50}$w_3$}};

\node [circle,fill=white!50,minimum size=0.75cm] (hidden3-4) at (0,-1.5) {$x_4$};
\node [] (w4) at (0.75,-1.75) {{\color{green!50}$w_4$}};

% CHAR.FUNC.

%\node [] (L) at (1,1.75) {$\mathcal{L}$};

% OUTPUT
\node [] (L) at (5,2.25) {$L^{GI}$};

\node [circle,fill=blue!30,minimum size=0.75cm] (output-1) at (5,1.5) {};
%\draw [->] (output-1) -- ++(1,0);
%\node [] (yhat-1) at (2+3,1) {$(\mathcal{L}(\v{x}))_1 = f(W_{\cdot \, 1}^\top \v{x} + b_1)$};

\node [circle,fill=blue!30,minimum size=0.75cm] (output-2) at (5,0.5) {};
%\draw [->] (output-2) -- ++(1,0);
%\node [] (yhat-1) at (2+3,0) {$\vdots$};

\node [circle,fill=blue!30,minimum size=0.75cm] (output-3) at (5,-0.5) {};
%\draw [->] (output-3) -- ++(1,0);
%\node [] (yhat-1) at (2+3,-1) {$(\mathcal{L}(\v{x}))_d = f(W_{\cdot \, d}^\top \v{x} + b_d)$};

\node [circle,fill=blue!30,minimum size=0.75cm] (output-4) at (5,-1.5) {};

% EDGES

\draw [->,color=cyan!100] (hidden3-1) -- (output-1);
\draw [->,color=cyan!100] (hidden3-1) -- (output-2);
\draw [->,color=cyan!100] (hidden3-1) -- (output-4);

\draw [->,color=magenta!100] (hidden3-2) -- (output-1);
\draw [->,color=magenta!100] (hidden3-2) -- (output-2);
\draw [->,color=magenta!100] (hidden3-2) -- (output-3);
\draw [->,color=magenta!100] (hidden3-2) -- (output-4);

\draw [->,color=orange!100] (hidden3-3) -- (output-2);
\draw [->,color=orange!100] (hidden3-3) -- (output-3);
\draw [->,color=orange!100] (hidden3-3) -- (output-4);

\draw [->,color=green!100] (hidden3-4) -- (output-1);
\draw [->,color=green!100] (hidden3-4) -- (output-2);
\draw [->,color=green!100] (hidden3-4) -- (output-3);
\draw [->,color=green!100] (hidden3-4) -- (output-4);

\end{tikzpicture}
    }
    \subcaptionbox{Weight matrix of GI layer\label{fig:GIb}}[0.45\textwidth]{
        \centering
        \input{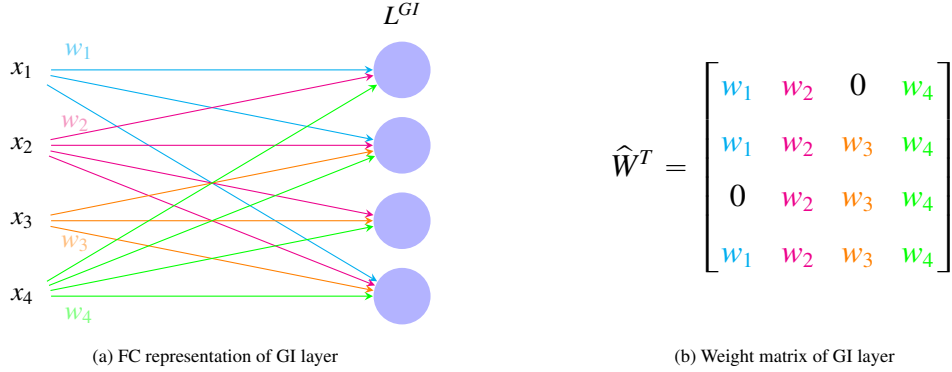}
        \vspace{0.35cm}
    }
    \caption{
    Visual representation of a GI layer as a ``constrained'' FC layer (subfigure ($A$)), with weight matrix defined by \eqref{eq:GI_weights_simple} (subfigure ($B$)). This figure is based on the same graph illustrated in Figure \ref{fig:ginnfilter}.
    }
    \label{fig:GIasFC}
\end{figure}

Looking better at equation \eqref{eq:GI_action_simple}, we observe that it is equivalent to the action of a ``constrained'' FC layer where the weights are the same if the connection is outgoing from the same unit, whereas it is zero if two units correspond to graph nodes that are not connected (see Figure \ref{fig:GIasFC}); more precisely:
\begin{equation}\label{eq:GI_weights_simple}
    \widehat{w}_{ij}=
    \begin{cases}
    w_i\,,\quad & \text{if }a_{ij}\neq 0 \text{ or }i=j\\
    0\,,\quad & \text{otherwise}
    \end{cases}
    \,,
\end{equation}
where $a_{ij},\widehat{w}_{ij}$ denote the $(i,j)$-th element of $A, \widehat{W}$, respectively.

GI layers characterized by \eqref{eq:GI_action_simple} can be generalized to read any arbitrary number $K\geq 1$ of input features per node and to return any arbitrary number $F\geq 1$ of output features per node. This property is crucial for the possibility of using GI layers as building blocks of our \pNN models, since they can process and transform easily matricial inputs like \eqref{eq:input_pNN}.
Then, the action of a general GI layer is a function $\mathcal{L}^{\rm GI}: \R^{{\mathrm{V}}\times K}\rightarrow\R^{{\mathrm{V}}\times F}$, with $K\geq1$ and $F\geq 1$. Additionally, pooling and other operations can be added; in general, for more details about GI layers, see \cite{GINN,DELLASANTA2025ewginn,sparseGIlayers}.
In particular, we point out that the number of weights of a GI layer is equal to ${\mathrm{V}}KF+{\mathrm{V}}F$, while the number of weights of a FC layer of $n_{\rm out}$ units, that is reading the outputs of a layer of $n_{\rm in}$ units, is equal to $n_{\rm out}n_{\rm in} + n_{\rm out}$; therefore, if we consider the case of $n_{\rm in}=n_{\rm out}={\mathrm{V}}$ and $KF + F < {\mathrm{V}}+1$ (satisfied for sufficiently large graphs, in general), the GI layer has fewer weights to be trained if compared with the FC layer. Moreover, the adjacency matrices are typically sparse and, therefore, the weight tensor is typically sparse too. Then, it is possible to exploit the sparsity of this tensor to reduce the memory cost of the GINN implementation{, storing the trainable parameters only} (see \cite{sparseGIlayers}). 

Now, we point the attention of the reader to the fact that the GI layer definition is very general and does not ask for a non-weighted graph. Indeed, one of the main advantages of using a GI layer is that it exploits the zeros of the adjacency matrix $A$, ``pruning'' the connections between layer units according to the graph edges; therefore, the zero elements of $A$ are in the same position independently on the usage of weights for the graph edges.
This observation is important to build GINNs with respect to the mesh of a given domain. In particular, we prefer a weighted graph because we want to emphasize the difference between connections according to the distance of the mesh nodes in the domain $\Omega$. Then, using a weighted graph with weights inversely proportional to the mesh edge lengths, we have that the action of the GI layers' weights is rescaled with respect to the non-trainable edge weights $\omega_{ij}\in (0, 1]${, corresponding to the non-zero elements of the graph's adjacency matrix.}

{We conclude this subsection recalling that, in a GINN, the relationship between the input features of a node and the output features of another node depends on the GINN’s depth and on the distance between these nodes in the graph}. Specifically, from Proposition 1 in \cite{GINN}, it holds that the input features corresponding to node $v_i$ contribute to the computation of the output features corresponding to node $v_j$ if the number $H$ of consecutive GI layers of the GINN is greater than or equal to the distance between $v_i$ and $v_j$ in the graph $G$ (in number of edges); i.e., if $H\geq\mathrm{dist}_G(v_i,v_j)$. Additionally, it is trivial to prove that this input-output relationship holds only if $H\geq\mathrm{dist}_G(v_i,v_j)$. The importance of the dependency of the input-output relationship on the GINN’s depth and the distances between nodes in the graph is important for selecting the depth hyper-parameter $H$ of the GINN.

% -------------------------------------------------------

\subsection{Architecture Archetypes for \pNNs}\label{sec:architectures}

Here, we describe two NN architecture archetypes based on GI layers and Conv1D-FC layers, respectively, for building \pNNs; i.e., surrogate models for approximating the solution $\solh(\v{\mu})$ of a parametric PDE with parametrized BCs, given the vector of parameters $\v{\mu}\in\mathcal{P}_h${, and a fixed mesh of $N_h$ nodes (i.e., a graph with $|V|=\mathrm{V}=N_h$)}. All the hyper-parameters of these architectures (e.g., activation functions, number of filters, etc.) have been chosen after preliminary analyses and studies. For each parametric PDE problem, we recall that the {two} NN architecture archetypes are characterized by the same input encoding $\mubcenc\in\R^{N_h\times(k+2+q)}$, see \eqref{eq:input_pNN}, and the same output encoding (a matrix in $\R^{N_h\times m}$); therefore, we will focus on the hidden layers that characterize and distinguish the architectures. Moreover, the architectures' archetypes have in common also the usage of Batch-Normalization layers in the residual blocks, and a set of layers for preprocessing and postprocessing operations (see Section \ref{sec:prepost_ops} below).

\subsubsection{{\pGINN:} \pNN Based on GI Layers} 

This NN archetype is made of residual blocks (see \cite{He2016_ResidualNN}) of GI layers based on the weighted adjacency matrix of the problem's mesh (weights $\omega_{ij}\in (0, 1]$, inversely proportional to edge lengths); then, we denote these NNs as \pGINNs. An important characteristic of our \pGINNs\ is the (preprocessed) inputs re-feeding every $\rho\in\N$ hidden GI layers; indeed, during preliminary analyses and studies, this strategy demonstrated to improve considerably the training performances of the \pGINN model. The main idea of the input re-feeding is to let the model remember the BC nature of the mesh nodes and to combine it with the values computed by the intermediate layers.

Let $H\in\N$ be the hyper-parameter characterizing the number of hidden GI layers of the \pGINN model; then, we set $H=\left[ \frac{3}{2} \mathrm{diam}(G) \right ]$, where $G$ is the mesh graph and $\mathrm{diam}(G)$ is the maximum distance between two nodes in $G$. This choice of $H$ is motivated by the outputs' dependency from the input features, based on the GINN's depth previously mentioned (see Section \ref{sec:GINNs}); indeed, with the first $H'=\mathrm{diam}(G)$ layers, we have the guarantee that any input feature contributes to the computation of the predicted output features of each graph node, while with the remaining $H''=\left[\frac{1}{2} \mathrm{diam}(G)\right]$ layers, we give to the NN extra layers for transforming the features and improving the approximation of the target solutions.

In summary, the architecture of a \pGINN is characterized by the hyper-parameters $\rho$, $H$, and the hyper-parameters shared by all the GI layers, i.e.: the activation function $\sigma$ and the number $F\in\N$ of output features per node.

\subsubsection{{\pNN Based on Conv1D and FC Layers}}

This NN archetype is made of residual blocks (see \cite{He2016_ResidualNN}) of FC layers with $N_h$ units each, and we denote them as \pFCNNs. However, before the FC layers, Conv1D layers are adopted to process the input $\mubcenc\in\R^{N_h\times (k+2+q)}$. In particular, the number of Conv1D layers used to transform an input $\mubcenc\in\R^{N_h\times (k+2+q)}$ into a vector in $\R^{N_h}$ is equal to the number of input features (i.e., $k+2+q$ Conv1D layers).
Concerning the \pFCNNs models, $H\in\N$ denotes the hyper-parameter characterizing the number of FC layers in the model's main structure. The number $H$ of FC layers used is set $H\geq 3$, and such that the total number of trainable parameters of the \pFCNN is approximately the same, or slightly greater than, the number of trainable parameters of a \pGINN designed for the same PDE problem. 
We set the lower bound of three FC layers for having at least one residual block of FC layers; indeed, due to the large number of weights that characterizes FC layers, for large meshes the number of weights of a \pFCNN increases rapidly, and deep models are almost never obtained. After this sequence of FC-based residual blocks, an extra FC layer is added.

In the end, if $m=1$ (i.e., $\solh(\v{\mu})\in\R^{N_h}$), after the FC layer of index $H+1$, there is a last FC layer. Otherwise, if $m>1$, we continue with other $m$ parallel sequences of $m$ FC layers with a residual connection to the $(H+1)$-th FC layer; then, the $m$ branches' outputs are concatenated into an $N_h$-by-$m$ matrix.

Therefore, the architecture of a \pFCNN is characterized by the hyper-parameter $H$, the number $F$ of filters for the Conv1D layers, and the activation function $\sigma$ of the Conv1D/FC layers.

\begin{remark}
    The input re-feeding, characteristic of \pGINNs, is not used in \pFCNNs. Indeed, the usage of FC layers makes it prohibitive to re-insert the (preprocessed) input $\mubcenc$, because FC layers are not able to ``read'' such kind of inputs. The only option would be to use Conv1D layers again, but this approach would increase too much the complexity of the NN architecture and, as a consequence, the number of trainable parameters.
\end{remark}

\subsubsection{Preprocessing and Postprocessing Operations}\label{sec:prepost_ops}

Both the NN architecture archetypes share a set of layers that can be used for performing preprocessing operations to the inputs $\mubcenc$, preprocessing operations to the targets, and/or postprocessing operations to the outputs.

The embedded preprocessing operations can be of many kinds, from classical ones (e.g., standardization or rescaling of the input/target values) to specific operations identified during preliminary analyses (e.g., change of the constant homogeneous Dirichlet/Neumann conditions into constant unitary values, switch the 0/1 convention for BC type to $\pm 1$ convention, etc.). The decision to embed the preprocessing operation is made in order to speed up the online phase of trained models and to ease the imposition of the Dirichlet BCs in the outputs. For more details about the advantages of embedded preprocessing, see Remark \ref{rem:adv_embedded_preproc} below.

On the other hand, the embedded postprocessing operations are divided into two types: the first type is the inverse of the preprocessing operations applied, during the training, to the targets (if there are any); the second type constrains the NN's predictions to satisfy the Dirichlet BCs indicated by the inputs.

\begin{remark}\label{rem:adv_embedded_preproc}
    Despite preprocessing operations are typically performed ``outside'' of the model and directly on the data, we preferred to include them ``inside'' the model in order to maintain the same input encoding $\mubcenc$ (see \eqref{eq:input_pNN}) for all the \pNN models, and avoiding the need of exporting alternate versions of the data or exporting preprocessing operators in case of need of reproducibility. This choice reduces the risk of performing wrong preprocessing for a third user and, in practical contexts, permits exporting and easily sharing a trained model.
    Moreover, embedded preprocessing operations are particularly helpful for the postprocessing operation dedicated to imposing the Dirichlet BCs. Indeed, having these operations embedded in the model, the input is always defined by the convention described in \eqref{eq:input_pNN}; therefore, the postprocessing operation used for imposing the Dirichlet BCs do not need to invert any preprocessing operation and only have to ``read-and-apply'' the BC values for Dirichlet nodes.
\end{remark}

% -------------------------------------------------------

\subsection{Training Loss for \pNNs}\label{sec:trainingloss}

Thanks to the embedded preprocessing and postprocessing operations, \pNN models can return an approximation of the target solution $\solh(\v{\mu})$ without the need of applying any inverse transformation on their outputs. Moreover, they can authomatically impose the Dirichlet BCs described in the input $\mubcenc$ on their predicted output $\nnsolh(\v{\mu})$. For this reason, we train \pNNs on a modified version of the Mean Square Error (MSE) function, for reducing the importance of prediction errors on Dirichlet nodes.

For better understanding the criteria behind the new loss function, we have to distinguish four types of mesh nodes, given a PDE problem: 
\begin{itemize}
    \item Internal mesh nodes;
    \item Boundary nodes always associated to Neumann BC (denoted as only-Neumann nodes);
    \item Boundary nodes always associated to Dirichlet BC (denoted as only-Dirichlet nodes);
    \item Boundary nodes that can be associated either to Dirichlet BC or to Neumann BC (denoted as BC-switch nodes);
\end{itemize}

In general, we can describe the mechanism of the modified MSE loss function in this way: internal mesh nodes and only-Neumann nodes fully contribute to the loss value, only-Dirichlet nodes do not contribute to the loss value, and BC-switch nodes contribute to the loss value is rescaled.
Mathematically, for each batch $\mathcal{B}$ of input-target pairs, the modified MSE loss function can be described as:
\begin{equation}\label{eq:weightedMSEloss}
    \mathrm{MSE}_{\rm BC}(\mathcal{B}):= \frac{1}{|\mathcal{B}|}\sum_{(\v{\mu},\solh)\in\mathcal{B}}
    \left(
    \frac{1}{N_h}\sum_{i=1}^{N_h}
    \lambda_i\big((\solh)_i - (\nnsolh(\v{\mu}))_i\big)^2
    \right)
    \,,
\end{equation}
where $\lambda_i=1$ if the mesh node $\v{x}_i$ is internal or only-Neumann, $\lambda_i=0$ if the mesh node $\v{x}_i$ is only-Dirichlet, and $\lambda_i=10^{-1}$ if the mesh node $\v{x}_i$ is BC-switch.

Concerning the boundary nodes, the idea behind the ${\rm MSE_{BC}}$ loss is to focus on predictions for the only-Neumann and BC-switch nodes, because Dirichlet BCs are imposed. For this reason, we set the values of $\lambda_i$ following this criterion: we have $\lambda_i\in (0,1]$ for BC-switch nodes (specifically, $\lambda_i=10^{-1}$ after preliminary analyses) because these nodes during the inference phase can be either Dirichlet nodes (with value imposed by reading the input) or Neumann nodes (with value predicted by the NN).

%%%%%%%%%%%%%%%%%%%%%%%%%%%%%%%%%%%%%%%%%%%%%%%%%%%%%%%%%%%%%%%%%%%%%%%%%%%%%%%%%%%%%%%%%%%

\section{Numerical Results}\label{sec:results}

In this section, we test the proposed \pGINN architecture in three different settings, comparing the results with a \pFCNN with a number of trainable weights of the same order of magnitude. The three settings are:
\begin{itemize}
    \item Experiment 1, where the boundary parameters define only a Neumann type boundary for a linear diffusion PDE;
    \item Experiment 2, where the boundary parameters define a mixed Dirichlet-Neumann type boundary for an advection-diffusion PDE;
    \item Experiment 3, where the boundary parameters define a mixed Dirichlet-Neumann type boundary for the Navier-Stokes Equations;
\end{itemize}
The data needed for the surrogate model training have been obtained through the employment of FEniCS \cite{fenics} {and its dependencies, running the numerical simulations on the node of a remote cluster endowed with up to 60GB of RAM, 1CPU used (Intel Xeon Silver 4110, 2.10GHz)}. The \pNNs have been developed using the Keras-TensorFlow module of Python \cite{keras2015,tensorflow2015-whitepaper}, and trained on node of a remote cluster endowed with up to 180GB of RAM, 10CPUs used (Intel Xeon E5-2640 v3, 2.60GHz), 1GPU used (GTX980 6G).

In what follows, we list the metric used to assess the accuracy and the robustness of the proposed approach based on \pGINNs, compared to \pFCNNs (representative of a generic NN-based approach).
For a given solution $\solh=\solh(\v{\mu})\in\R^{N_h}$ of the test set $\Xi$ and the {\pNN prediction $\nnsolh=\nnsolh(\v{\mu})\in\R^{N_h}$, we define the following absolute errors}

\begin{equation}\label{eq:abserr_onesol}
            \begin{aligned}
                \mathrm{e}_{\ell_2}(\v{\mu},\solh) = \Norm{\solh - \nnsolh}_2\,, 
                \quad
                & 
                \mathrm{e}_{L^2}(\v{\mu},\solh) = \Norm{\solh - \nnsolh}_{\mathsf M}\,,
                & % \\
                \quad
                \mathrm{e}_{H^1}(\v{\mu},\solh) =  \Norm{\solh - \nnsolh}_{\mathsf A}\,,
            \end{aligned}
\end{equation}

and relative errors

\begin{equation}\label{eq:relerr_onesol}
            \begin{aligned}
                \mathrm{re}_{\ell_2}(\v{\mu},\solh) = \frac{\Norm{\solh - \nnsolh}_2}{\Norm{\solh}_2}\,, 
                \quad
                & 
                \mathrm{re}_{L^2}(\v{\mu},\solh) = \frac{\Norm{\solh - \nnsolh}_{\mathsf M}}{\Norm{\solh}_{\mathsf M}}\,,
                & % \\
                \quad
                \mathrm{re}_{H^1}(\v{\mu},\solh) =  \frac{\Norm{\solh - \nnsolh}_{\mathsf A}}{\Norm{\solh}_{\mathsf A}}\,,
            \end{aligned}
\end{equation}

where $\mathsf M$ and $\mathsf A$ represent the mass and stiffness matrices of the $\mathbb P^1$ finite element PDE discretization, respectively. Namely, $\mathrm{e}_{L^2}$ and $\mathrm{e}_{H^1}$ provide the numerical approximation of the $L^2$-norm and the $H^1$-seminorm. Of these quantities, we will consider the mean over the whole test set $\Xi$; i.e., we consider the mean absolute, or relative, errors:

\begin{equation}\label{eq:err_mean}
            \begin{aligned}
                \MlE(\Xi) & = &\frac{1}{|\Xi|}\sum_{(\v{\mu}, \solh)\in\Xi}\mathrm{e}_{\ell_2}(\v{\mu}, \solh)\,,
                \qquad
                \MlRE(\Xi) & = &\frac{1}{|\Xi|}\sum_{(\v{\mu}, \solh)\in\Xi}\mathrm{re}_{\ell_2}(\v{\mu}, \solh)\,,
                \\
                \MLE(\Xi) & = &\frac{1}{|\Xi|}\sum_{(\v{\mu}, \solh)\in\Xi}\mathrm{e}_{L^2}(\v{\mu}, \solh)\,,
                \qquad
                \MLRE(\Xi) & = &\frac{1}{|\Xi|}\sum_{(\v{\mu}, \solh)\in\Xi}\mathrm{re}_{L^2}(\v{\mu}, \solh)\,,
                \\
                \MHE(\Xi) & = &\frac{1}{|\Xi|}\sum_{(\v{\mu}, \solh)\in\Xi}\mathrm{e}_{H^1}(\v{\mu}, \solh)\,,
                \qquad
                \MHRE(\Xi) & = &\frac{1}{|\Xi|}\sum_{(\v{\mu}, \solh)\in\Xi}\mathrm{re}_{H^1}(\v{\mu}, \solh)\,.
            \end{aligned}
\end{equation}

In particular, for each configuration of training options and hyper-parameters, the \pNNs models are trained 5 times, varying the random seed, and the statistics of the mean errors \eqref{eq:err_mean} are analyzed and reported; this operation is necessary to better understand the robustness of a model with respect to its initialization.

Concerning the different configurations considered, they are characterized by a decreasing training set size and validation set size, while the test set is kept fixed (also for maintaining a rigorous comparison of the performances). The training options adopted for all the experiments are the following:
\begin{itemize}
    \item mini-batch size of 16 samples;
    \item early stopping patience of $250$ epochs returning the best weights \cite{EarlyStopping_TF}, over $10\,000$ maximum number of epochs;
    \item reduce learning rate on plateaus with patience of 100 epochs and factor 0.5 \cite{ReduceLRPlateau_TF};
    \item Adam optimizer \cite{Kingma2015_ADAM} (starting learning rate $0.001$, moment decay rates $\beta_1=0.9$, $\beta_2=0.999$).
\end{itemize}

A valuable set of architecture hyper-parameters has been identified during preliminary analyses, and they are fixed; we postpone to future work an extensive analysis to identify the best combination of hyper-parameters for each model.

Additionally, not only the approximation performance is measured, but also the computational costs of training the models, in terms of average time required by the model for performing a weight-update with respect to a mini-batch of the training set.

In the end, we recall the different notation between continuous formulation and discrete formulation of a problem parametrization. The parameter $\bmu\in\mathcal{P}$ characterizing the PDE problem in the continuous formulation is such that $\bmu=(\bmu_{\phi}, \mu_b, \mu_v)\in\mathcal{P}$, with $\mathcal{P} = \mathbb R^{p_{\phi}} \times \mathcal{F}_{pwc}(\Gamma^{\bmu_b}) \times \mathcal{F}_{pwc}(\Gamma^{\bmu_b})$ (see Section \ref{sec:prob}); on the other hand, the parameter $\bmu\in\mathcal{P}_h$ characterizing the PDE problem in the discrete formulation is such that $\bmu=(\bmu_{\phi}, \bmu_b, \bmu_v)\in\mathcal{P}_h$, with $\mathcal{P}_h = \R^{p_{\phi}}\times \R^{p_b} \times \R^{p_v}$ (see Section \ref{sec:discrete_form}), and it is encoded as illustrated in \eqref{eq:input_pNN}, Section \ref{sec:GINNgen}, for feeding the \pNNs.

\subsection{Experiment 1 - Neumann linear case (Diffusion)}

\begin{figure}[htb!]
    \begin{subfigure}[t]{0.5\textwidth}
        \begin{center}
    \begin{tikzpicture}[scale=5.0]
\filldraw[color=white, fill=gray!10, very thick, loosely dotted](0.5,0.5) circle (0.3);
\filldraw[color=black, very thick, densely  dotted](0,0) -- (0.,1);
\filldraw[color=black, fill=gray!10, very thick](0,1) -- (1,1);
\filldraw[color=black, fill=gray!10, very thick](0,0.) -- (1,0.);
\filldraw[color=black, fill=gray!10, very thick, densely dotted](1,0.) -- (1,1);
\draw [black, very thick,domain=0:45, dashed] plot ({0.5+0.3*cos(\x)}, {0.5+0.3*sin(\x)});
\draw [black, very thick,domain=180:290, dashed] plot ({0.5+0.3*cos(\x)}, {0.5+0.3*sin(\x)});
\draw [black, thick,domain=45:180, dashdotted] plot ({0.5+0.3*cos(\x)}, {0.5+0.3*sin(\x)});
\draw [black, thick,domain=290:360, dashdotted] plot ({0.5+0.3*cos(\x)}, {0.5+0.3*sin(\x)});
\node at (-.17,0.5){\color{black}{$\Gamma_{N}$}};
\node at (0.5,-0.1){\color{black}{$\Gamma_{D}$}}; 
\node at (0.4,0.35){\color{black}{$\Gamma_1^{\boldsymbol{\mu}_b}$}}; 
\node at (0.5,0.65){\color{black}{$\Gamma_0^{\boldsymbol{\mu}_b}$}}; 

%\node at (2.3,0.205){\color{black}{$\Gamma_{N}$}};
%
\node at (-.1,0){\color{black}{$(0,0)$}};

\node at (-.1,1.){\color{black}{$(0,1)$}};
\node at (1.1,1.){\color{black}{$(1,1)$}};
\node at (1.1,-0.){\color{black}{$(1,0)$}};
%\node at (1.1,-.05){\color{black}{$ \color{cyan}{\Gamma_\text{w}}$}};
\end{tikzpicture}
    \caption{}
    \label{subfig:domain}
    \end{center}
        \end{subfigure}%
        ~ 
        \begin{subfigure}[t]{0.5\textwidth}
            \begin{center}
    % \begin{tikzpicture}[scale=5.0]
% \filldraw[color=white, fill=gray!10, very thick, dotted](0.5,0.5) circle (0.3);
% \draw [black, very thick,domain=0:45, dashed] plot ({0.5+0.3*cos(\x)}, {0.5+0.3*sin(\x)});
% \draw [black, very thick,domain=180:290, dashed] plot ({0.5+0.3*cos(\x)}, {0.5+0.3*sin(\x)});
% \draw [black, very thick,domain=45:180, densely dotted] plot ({0.5+0.3*cos(\x)}, {0.5+0.3*sin(\x)});
% \draw [black, very thick,domain=290:360, densely dotted] plot ({0.5+0.3*cos(\x)}, {0.5+0.3*sin(\x)});
% \node at (0.4,0.35){\color{black}{$\Gamma_{N}^{\mu_N^*}$}}; 
% \filldraw[color=white, fill=black](0.32,0.26)circle (0.015);
% %\node at (1.1,-.05){\color{black}{$ \color{cyan}{\Gamma_\text{w}}$}};
% \end{tikzpicture}
\tikzset{cross/.style={cross out, draw=black, minimum size=2*(#1-\pgflinewidth), inner sep=0pt, outer sep=0pt},
%default radius will be 1pt. 
cross/.default={3pt}}
\begin{tikzpicture}[scale=5.0]
\filldraw[color=white, fill=white, very thick, dotted](0.5,0.5) circle (0.3);
\filldraw[color=white, very thick, densely  dotted](0,0) -- (0.,1);
\filldraw[color=white, fill=gray!10, very thick](0,1) -- (1,1);
\filldraw[color=white, fill=gray!10, very thick](0,0.) -- (1,0.);
\filldraw[color=white, fill=gray!10, very thick, densely dotted](1,0.) -- (1,1);
\draw [black, very thick,domain=0:45, dashed] plot ({0.5+0.3*cos(\x)}, {0.5+0.3*sin(\x)});
\draw [black, very thick,domain=180:290, dashed] plot ({0.5+0.3*cos(\x)}, {0.5+0.3*sin(\x)});
\draw [black, thick,domain=45:180, dashdotted] plot ({0.5+0.3*cos(\x)}, {0.5+0.3*sin(\x)});
\draw [black, thick,domain=290:360, dashdotted] plot ({0.5+0.3*cos(\x)}, {0.5+0.3*sin(\x)});
\node at (0.4,0.35){\color{white}{$\Gamma_{N}^{\mu_N^*}$}}; 
%\node at (1.1,-.05){\color{black}{$ \color{cyan}{\Gamma_\text{w}}$}};
\node at (-.17,0.5){\color{white}{$\Gamma_{N}\setminus \Gamma_N^{\mu_N^*}$}};
\node at (0.5,-0.1){\color{white}{$\Gamma_{D}$}}; 
\node at (0.4,0.35){\color{white}{$\Gamma_{N}^{\mu_N^*}$}}; 

%\node at (2.3,0.205){\color{black}{$\Gamma_{N}$}};
%
\node at (-.1,0){\color{white}{$(0,0)$}};

\node at (-.1,1.){\color{white}{$(0,1)$}};
\node at (1.1,1.){\color{white}{$(1,1)$}};
\node at (1.1,-0.){\color{white}{$(1,0)$}};
\filldraw[color=white, fill=black](0.32,0.26)circle (0.015);
\draw (0.6,0.22) node[cross,rotate=10, very thick] {};
\draw (0.2,0.5) node[cross,rotate=10, very thick] {};
\node at (0.07,0.5){\color{black}{$c_1 + \frac{l_1}{2}$}};
\node at (0.6,0.15){\color{black}{$c_1 - \frac{l_1}{2}$}};
\node at (0.32,0.19){\color{black}{$c_1$}};
\filldraw[color=white, fill=black](0.78,0.61)circle (0.015);
\draw (0.8,0.5) node[cross,rotate=10, very thick] {};
\draw (0.73,0.7) node[cross,rotate=10, very thick] {};
\node at (0.86,0.61){\color{black}{$c_2$}};
\node at (0.87,0.74){\color{black}{$c_2 - \frac{l_2}{2}$}};
\node at (0.95,0.5){\color{black}{$c_2 + \frac{l_2}{2}$}};
\end{tikzpicture}
    \caption{}
    \label{subfig:interval}
\end{center}
    \end{subfigure}%
    \caption{
    Experiment 1. Spatial domain $\Omega$ for a specific parametric instance: schematic representation. (A) The homogeneous Dirichlet boundary, the non-homogeneous and homogeneous Neumann boundaries are applied to $\Gamma_D$ (solid line), $\Gamma_1^{\bmu_b}$ (dashed line), and $\Gamma_0^{\bmu_b} \cup \Gamma_N$ (dash-dotted and dotted lines), respectively. (B) Example of selection of circle intervals. The summation only represents the detection of the start and the end node of the intervals.
    }
\label{fig:domain1}
\end{figure}
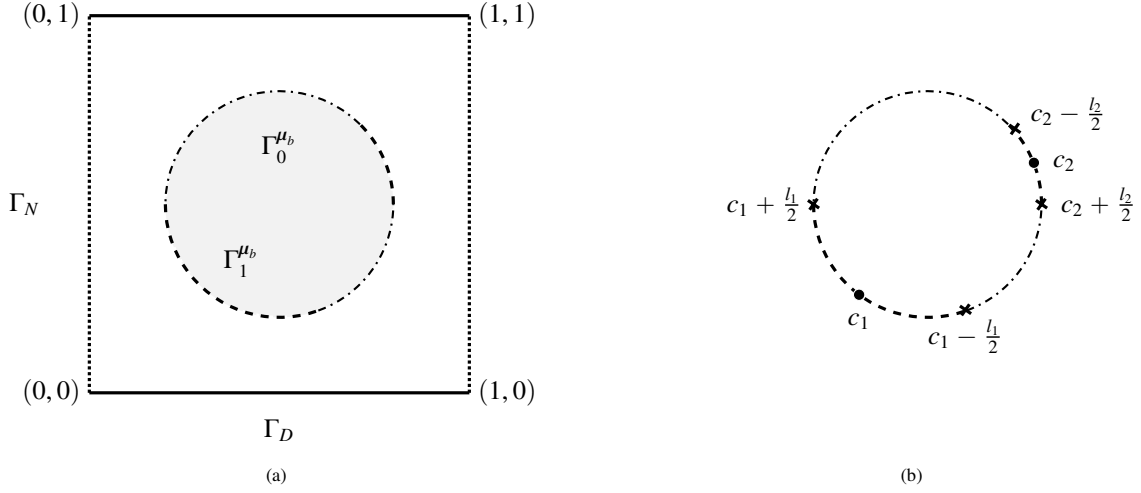

Let us consider $\Omega = [0, 1]^2 \setminus B\left((0.5, 0.5), 0.3\right)$, as depicted in Figure \ref{subfig:domain}; $B\left((0.5, 0.5), 0.3\right)$ denotes the open ball of center $(0.5, 0.5)$ and radius $0.3$. Specifically, we represent a domain for a generic $\bmu=(\mu_b, \mu_v)\in\mathcal{P}$, where $\mu_b$ defines $\Gamma^{\bmu_b}$ as the inner circular boundary $\partial B\left((0.5, 0.5), 0.3\right)$ and characterizes it into the portions $\Gamma^{\bmu_b}_1$ and $\Gamma^{\bmu_b}_0$ such that $\Gamma^{\bmu_b}_1 \cup \Gamma^{\bmu_b}_0=\Gamma^{\bmu_b}=\partial B\left((0.5, 0.5), 0.3\right)$, where we apply the following Neumann BCs via $\mu_v$: $\mu_v(\Gamma^{\bmu_b}_1)=\mu_1\in(0,1]$, $\mu_v(\Gamma^{\bmu_b}_0)=0$. For the ease of notation, we dropped $\bmu_{\phi}$ in $\bmu$ because the physical parameters are not parametrized in this experiment.

As we stated above, we denote as $\Gamma_1^{\bmu_b}$ the portion of $\Gamma^{\bmu_b}$ with non-homogeneous Neumann BCs of value $\mu_1$, while homogeneous Neumann BCs are applied to $\Gamma^{\bmu_b}_0 = \Gamma^{\bmu_b} \setminus \Gamma_1^{\bmu_b}$. The external boundary is $\Gamma_{\text{fix}} = \Gamma_N \cup \Gamma_D$, with $\Gamma_D = ((0,1) \times \{0\}) \cup ((0,1) \times \{1\})$ denoting homogeneous Dirichlet BC, and $\Gamma_N = (\{0\} \times (0,1)) \cup (\{1\} \times (0,1))$ denoting homogeneous Neumann BC. See Figure \ref{subfig:domain}, for a visual example.
In particular, the problem we are solving reads:
$$
\begin{cases}
    - \Delta u = 0 & \text{ in } \Omega,\\
    \displaystyle \frac{\partial{u}}{\partial \v{n}} = \mu_1 & \text{ on } \Gamma_1^{\bmu_b}, \vspace{1mm}\\
    \displaystyle \frac{\partial{u}}{\partial \v{n}} = 0 & \text{ on } \Gamma_0^{\bmu_b} \text{ and on }  \Gamma_{N}, \\
    u = 0 & \text{ on } \Gamma_D
\end{cases}
$$
where $\v{n}$ is the normal outward vector to $\partial \Omega$.

We now present the rationale behind the $\bmu=(\mu_b, \mu_v)$ variability distribution, specifically used for generating the training, validation, and test data.
The parameter $\mu_v(\Gamma_1^{\bmu_b})=\mu_1$ is randomly generated with uniform distribution between zero and one; i.e., $\mu_1\sim\mathcal{U}(0, 1)$.
Each $\mu_b$ is generated by randomly selecting $N_c\in\N$ points on $\Gamma^{\bmu_b}$, with $N_c$ minimum one and up to five; i.e., $N_c\sim \mathcal{U}(1,5)$. Then, we denote these points as $\{\v{c}_i\}_i^{N_c}\subset\Gamma^{\bmu_b}$. After that, for each point $\v{c}_i$, we randomly generate an interval on $\Gamma^{\bmu_b}$ that is centered on $\v{c}_i$ and that will denote a portion of the boundary with non-homogeneous Neumann BC. For generating these intervals, we generate a length $l_i\sim\mathcal{U}(0, L_c/5)$, for each $i=1, \ldots ,N_c$, where $L_c$ is the length of the circumference of $\Gamma^{\bmu_b}$; then, the $i$-th interval is defined as $I_i=[\v{c}_i - l_i/2, \v{c}_i + l_i/2]$, given a parametrized representation of the circumference, and the union of all the intervals $I_i$ defines $\Gamma_1^{\bmu_b}$. See Figure \ref{subfig:interval} for a visual example. It is worth noting that this generation procedure can generate overlapping intervals or degenerate intervals of null length (i.e., $I_i=\{\v{c}_i\}$).

The above generation procedure for $\mu_b$ can be easily translated to a discrete representation of the domain, after the creation of the mesh, for generating a vector $\bmu_b$. Denoting with $\v{x}_1,\ldots ,\v{x}_M$ the $M$ points of the mesh on the boundary $\Gamma^{\bmu_b}$, the $N_c$ points are randomly sampled among them, while the interval lengths are measured in number of interval nodes, i.e.: $l_i\sim\mathcal{U}(1, [M/5])$. Therefore, in the discrete case (and in practice) the parameter $\bmu_b$ is a vector in $\{0, 1\}^M\subset\R^M$ (i.e., $p_b=M$) such that its elements are $\bmu_b\pe{j}=1$ if $\v{x}_j\in\Gamma_1^{\bmu_b}$, and $\bmu_b\pe{j}=0$ otherwise.

\begin{remark}[ROMs for Neumann varying boundary conditions]
    Projection-based ROMs can be applied in the proposed context, but inefficiently. Indeed, the problem at hand is not affine and would require additional hyper-reduction techniques, such as DEIM or EIM, to recover online efficiency. Moreover, as already experimented in \cite{StrazzulloVicini}, DEIM performs poorly when the variation is located mostly along the boundary. Furthermore, the complex behaviour related to different locations of non-homogeneous Neumann boundaries needs a large number of basis functions to properly represent the phenomenon, even with clustering strategies and local POD. 
\end{remark}

\subsubsection{Experiment 1 - Results}\label{sec:results_exp1}

For this first experiment, we discretized the domain $\Omega$ with a mesh of $N_h=1141$ nodes, $p_b=88$ nodes on the boundary $\Gamma^{\bmu_b}$, generating a mesh graph $G$ with diameter $\mathrm{diam}(G)=54$. The dataset is generated by running $5000$ simulations through a $\mathbb{P}^1$ finite element discretization on the given mesh; in particular, each simulation is done with respect to a randomly generated parameter vector $\bmu=(\bmu_b, \bmu_v)=(\bmu_b, \mu_1, 0)\in\{0, 1\}^{p_b}\times (0, 1]\times \{0\}$. Given these 5000 simulation data, 2048 of them are selected randomly and used as a fixed test set; from the remaining data, we randomly extract {an even number $T\in\N$ of training set data and $T/2$ validation set data}. Each model (\pGINN and \pFCNN) is trained five times, with respect to five different random seeds for weight initialization, varying the amount of training data ${T} = 1024, 512, 256, 128$.

\begin{figure}[htb!]
    \centering
    \subcaptionbox{$\MlE$}{\includegraphics[trim=0.5cm 0.0cm 2.5cm 1.8cm,clip,width=0.45\textwidth]{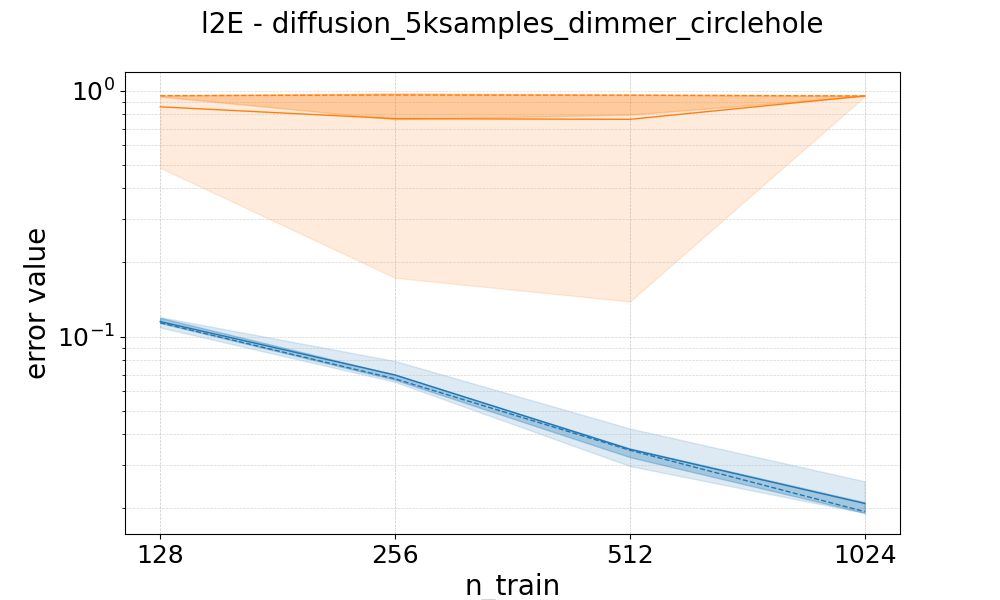}
    }
	\subcaptionbox{$\MLE$}{\includegraphics[trim=0.5cm 0.0cm 2.5cm 1.8cm,clip,width=0.45\textwidth]{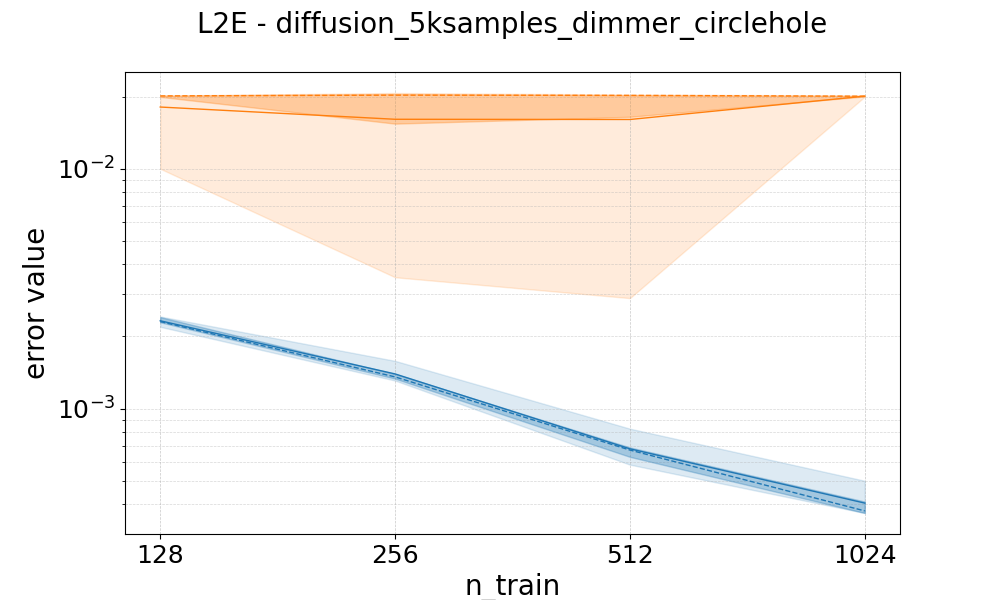}
    }
    \\
    \subcaptionbox{$\MHE$}{\includegraphics[trim=0.5cm 0.0cm 2.5cm 1.8cm,clip,width=0.45\textwidth]{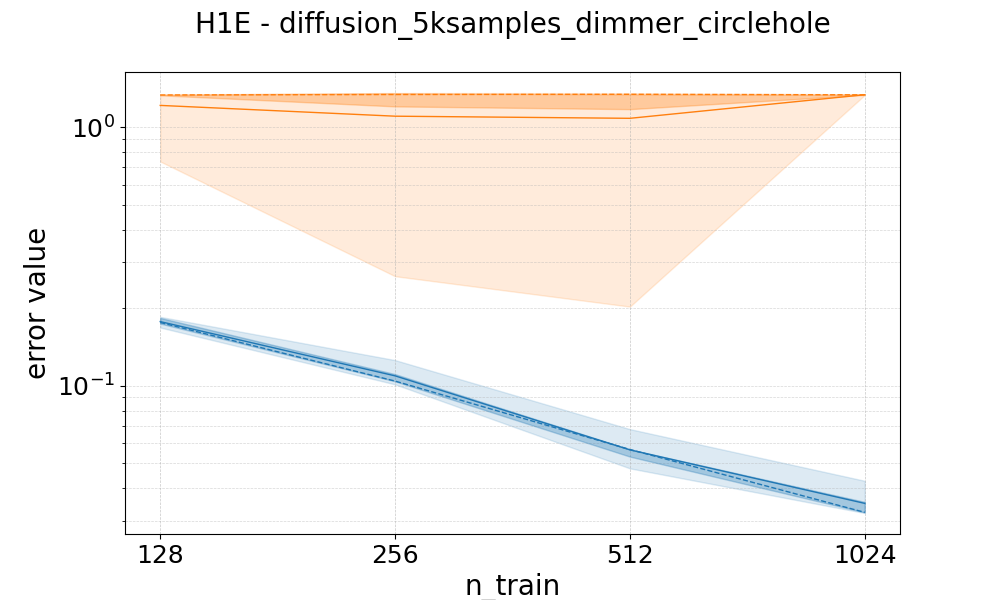}
    }
    \subcaptionbox{legend}{
    \includegraphics[width=0.45\textwidth]{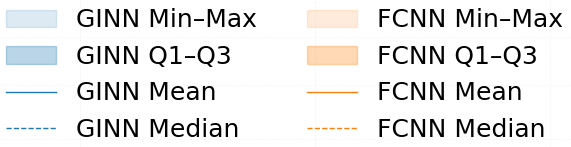}
    \vspace{0.75cm}
    }
    \caption{
    Experiment 1 (Diffusion) - Error statistics with respect to five random seeds for weight initialization, varying the training set size. In blue, the \pGINN performance, in orange the \pFCNN performance. Light colored areas represent the Min-Max range of values, the dark colored areas represent the values between fist and third quartiles. The continuous lines represent the average errors (same as Table \ref{tab:results_exp1}), the dotted lines represent the medians.
    }
    \label{fig:results_exp1}
\end{figure}

\begin{figure}[htb!]
    \centering
    \subcaptionbox{Best case - $\solh$}{\includegraphics[trim=1.95cm 4.cm 7cm 5cm,clip,width=0.28\textwidth]{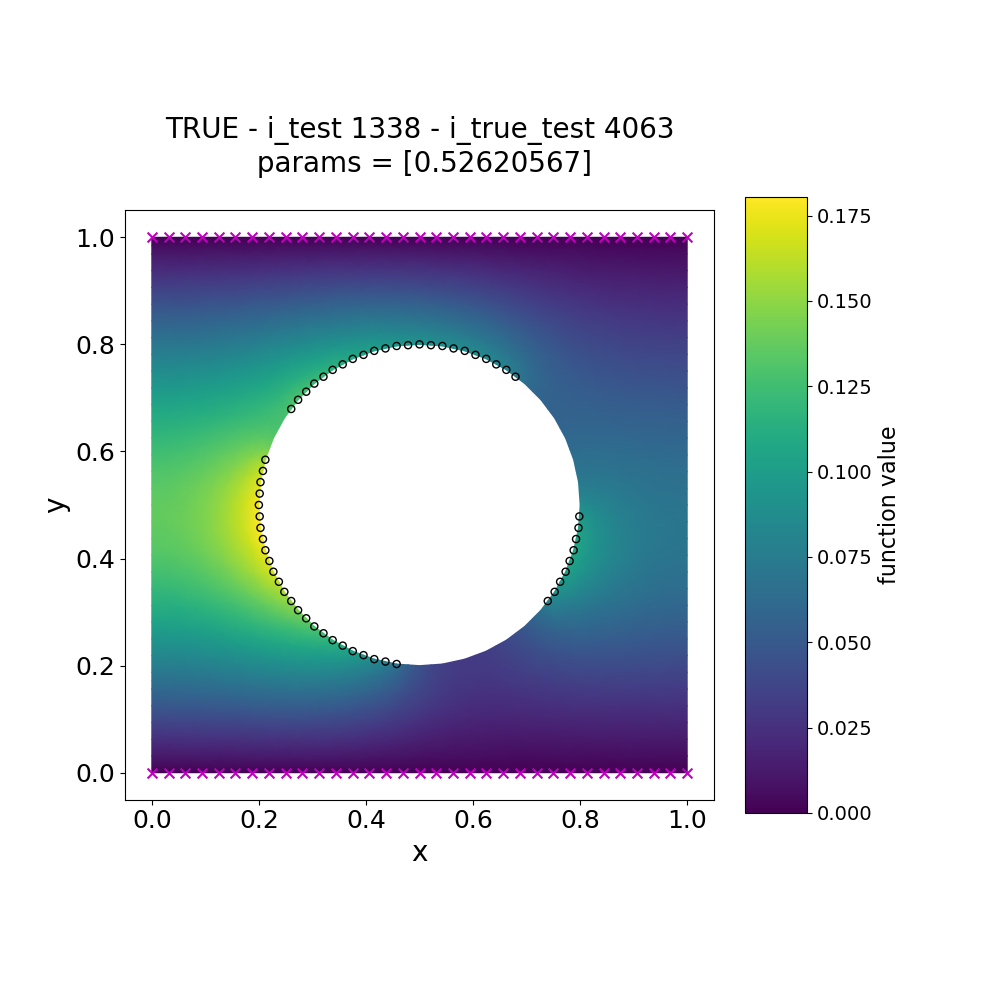}
    }
    \subcaptionbox{Best case - $\nnsolh$}{\includegraphics[trim=1.95cm 4.cm 3.55cm 4.75cm,clip,width=0.335\textwidth]{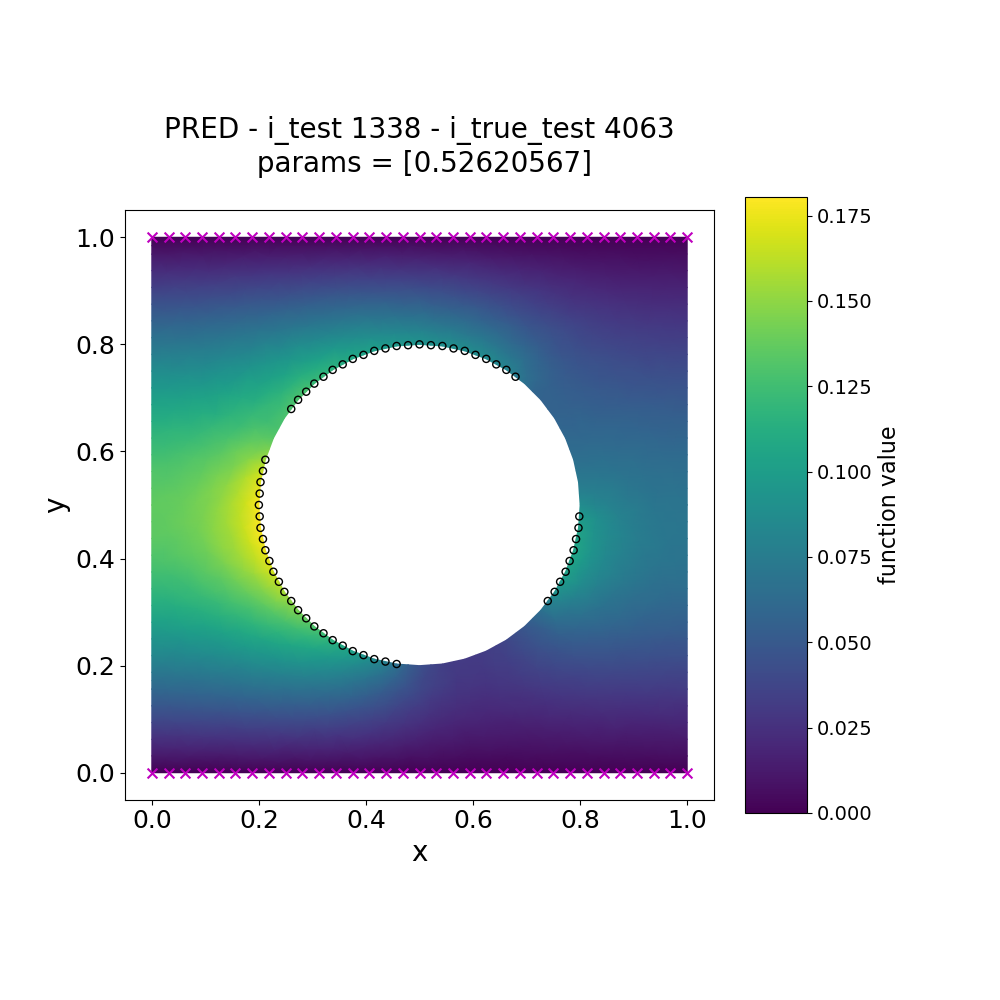}
    }
    \subcaptionbox{Best case - $|\solh - \nnsolh|$}{\includegraphics[trim=1.95cm 4.cm 3.25cm 4.75cm,clip,width=0.345\textwidth]{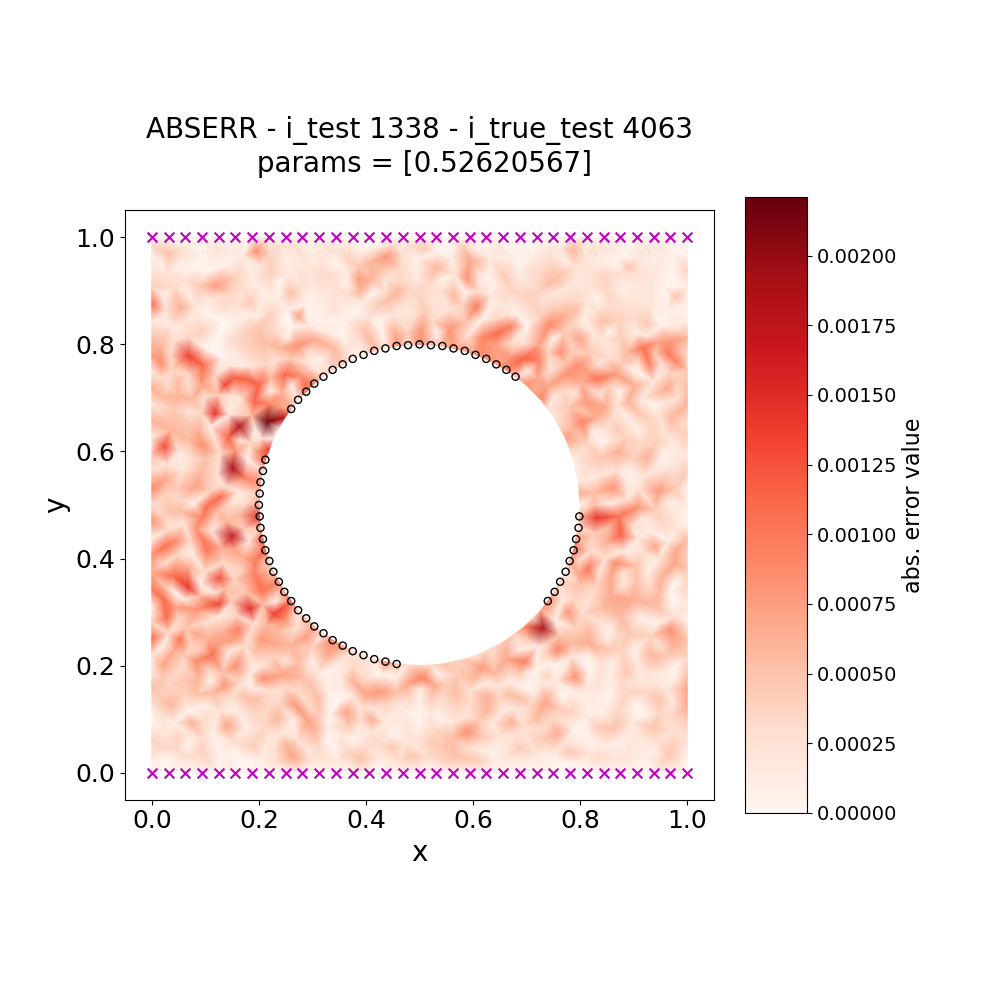}
    }
    \\
    \subcaptionbox{Median case - $\solh$}{\includegraphics[trim=1.95cm 4.cm 7cm 5cm,clip,width=0.28\textwidth]{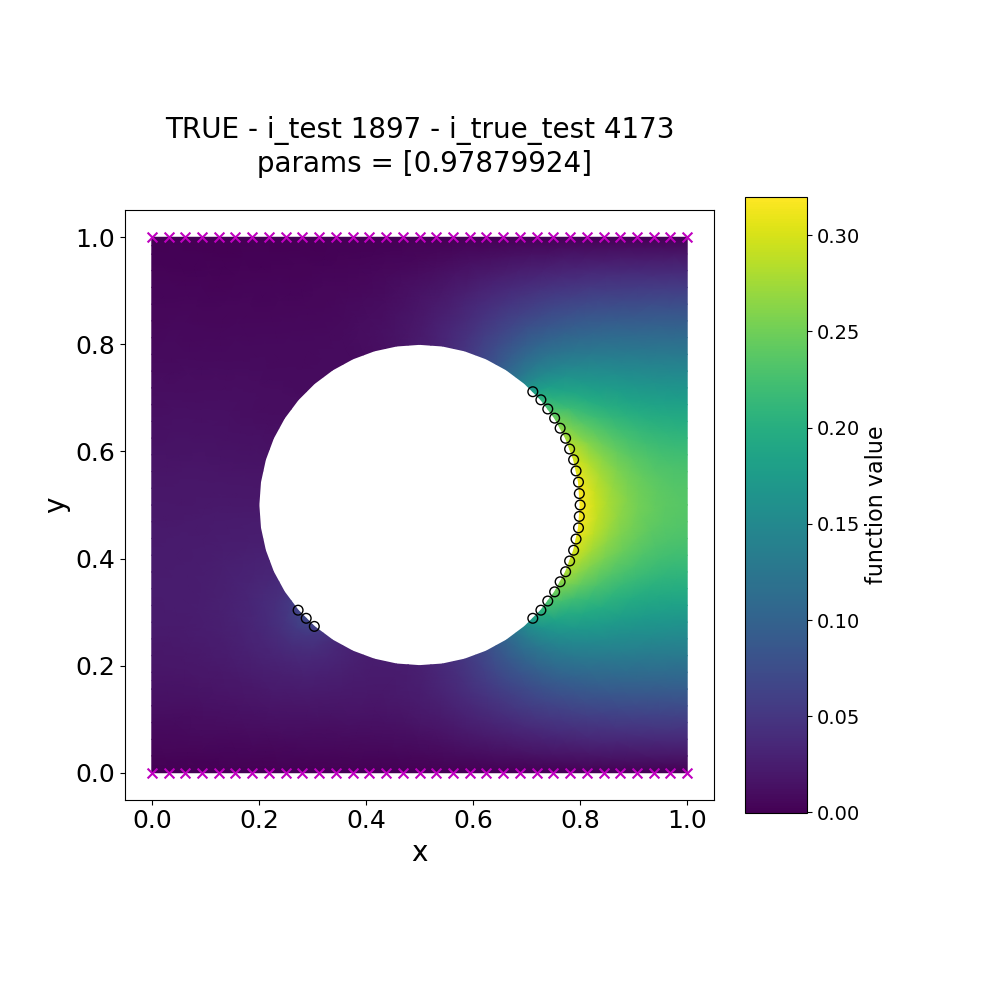}
    }
    \subcaptionbox{Median case - $\nnsolh$}{\includegraphics[trim=1.95cm 4.cm 3.55cm 4.75cm,clip,width=0.335\textwidth]{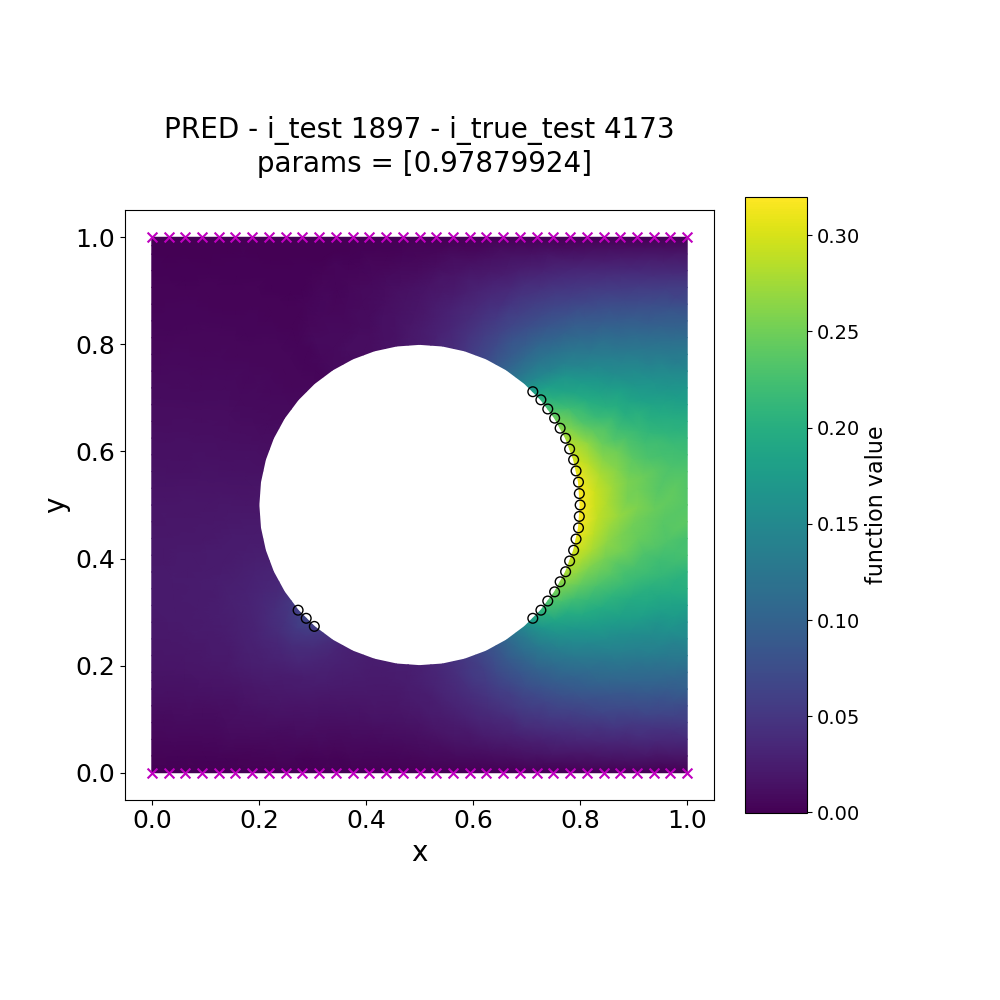}
    }
    \subcaptionbox{Median case - $|\solh - \nnsolh|$}{\includegraphics[trim=1.95cm 4.cm 3.25cm 4.75cm,clip,width=0.345\textwidth]{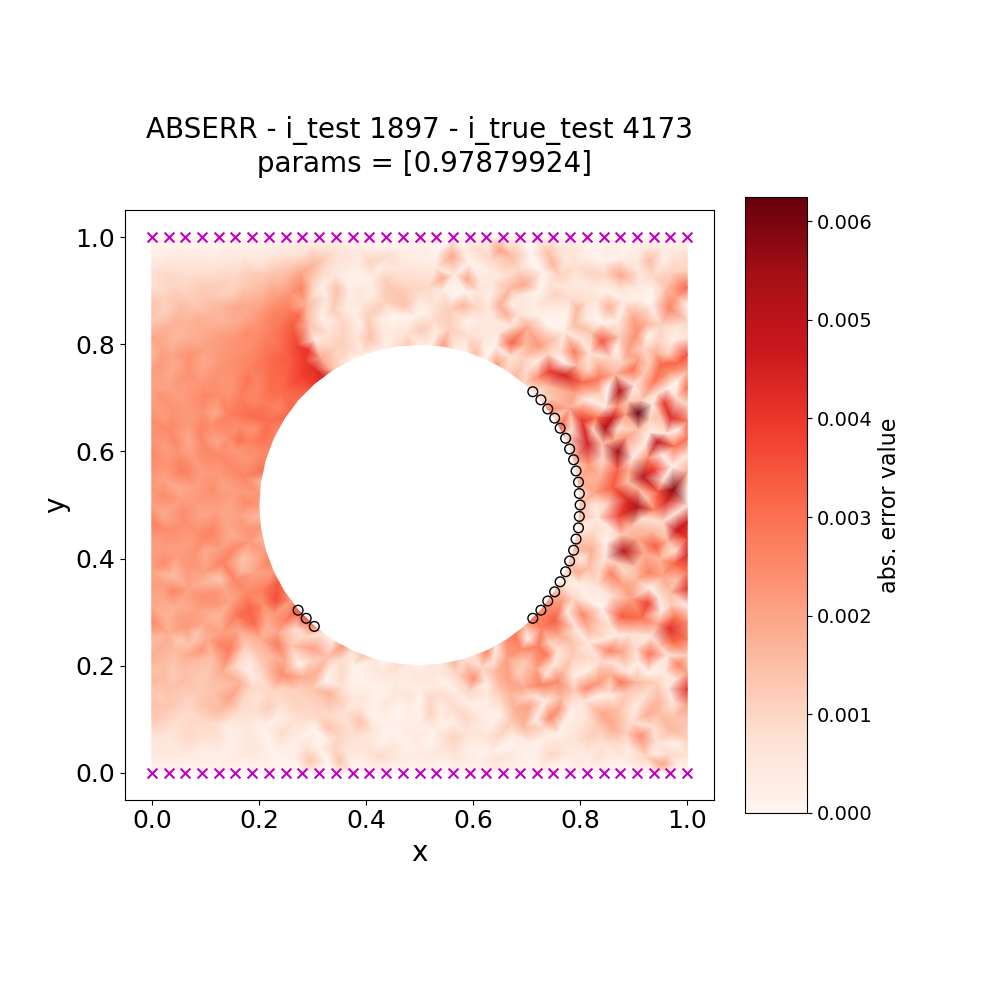}
    }
    \\
    \subcaptionbox{Worst case - $\solh$}{\includegraphics[trim=1.95cm 4.cm 7cm 5cm,clip,width=0.28\textwidth]{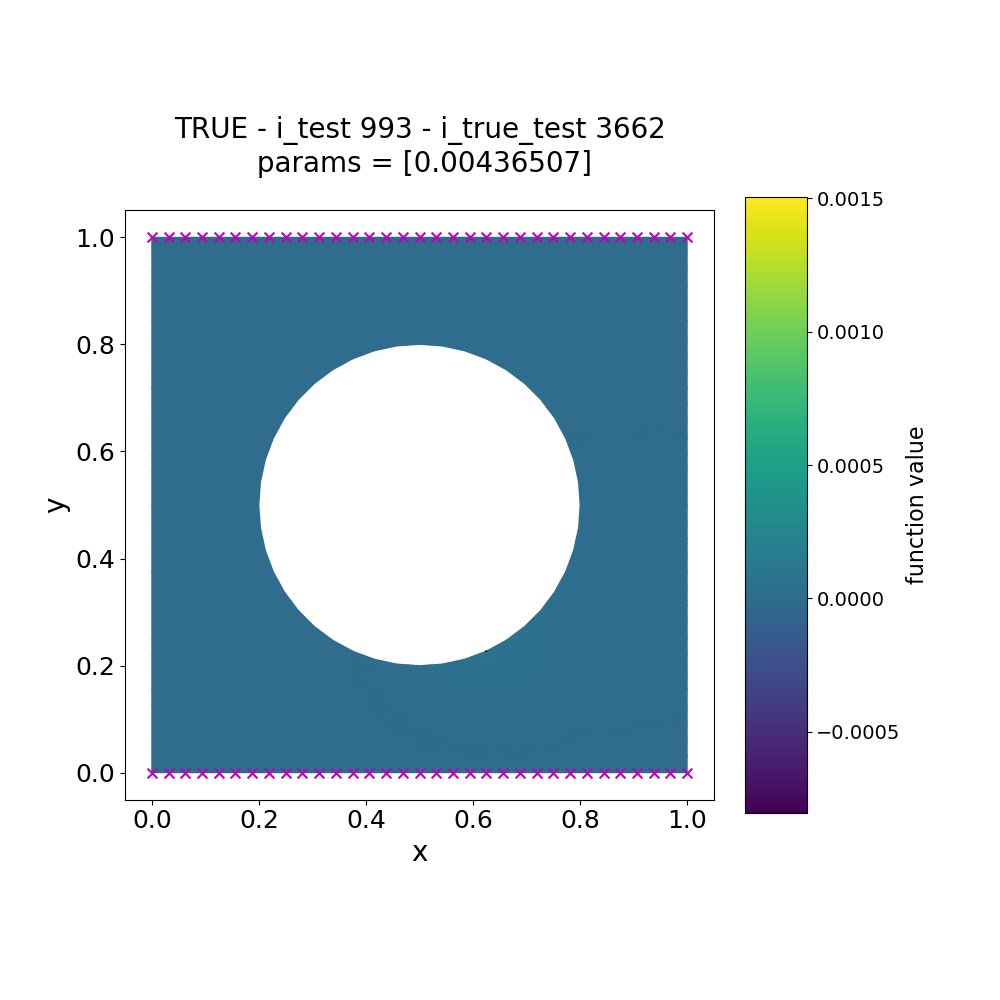}
    }
    \subcaptionbox{Worst case - $\nnsolh$}{\includegraphics[trim=1.95cm 4.cm 3.55cm 4.75cm,clip,width=0.335\textwidth]{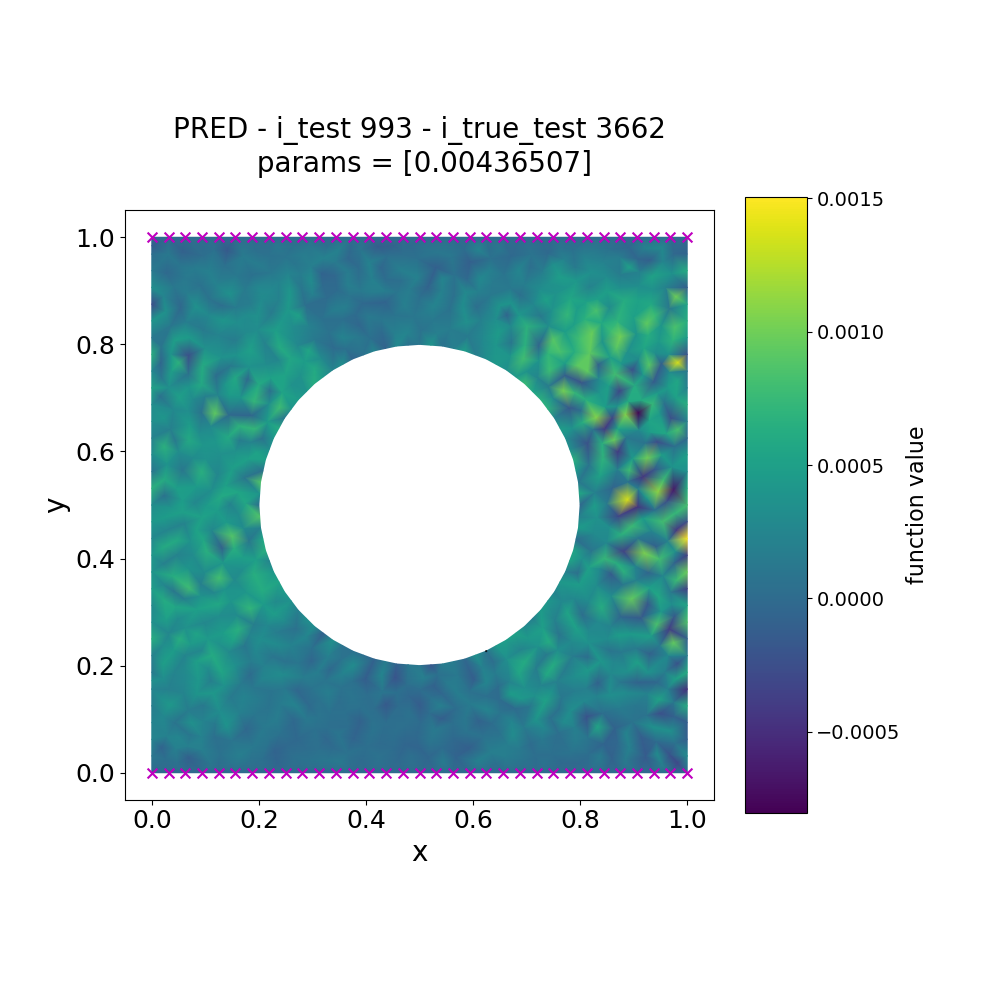}
    }
    \subcaptionbox{Worst case - $|\solh - \nnsolh|$}{\includegraphics[trim=1.95cm 4.cm 3.25cm 4.75cm,clip,width=0.345\textwidth]{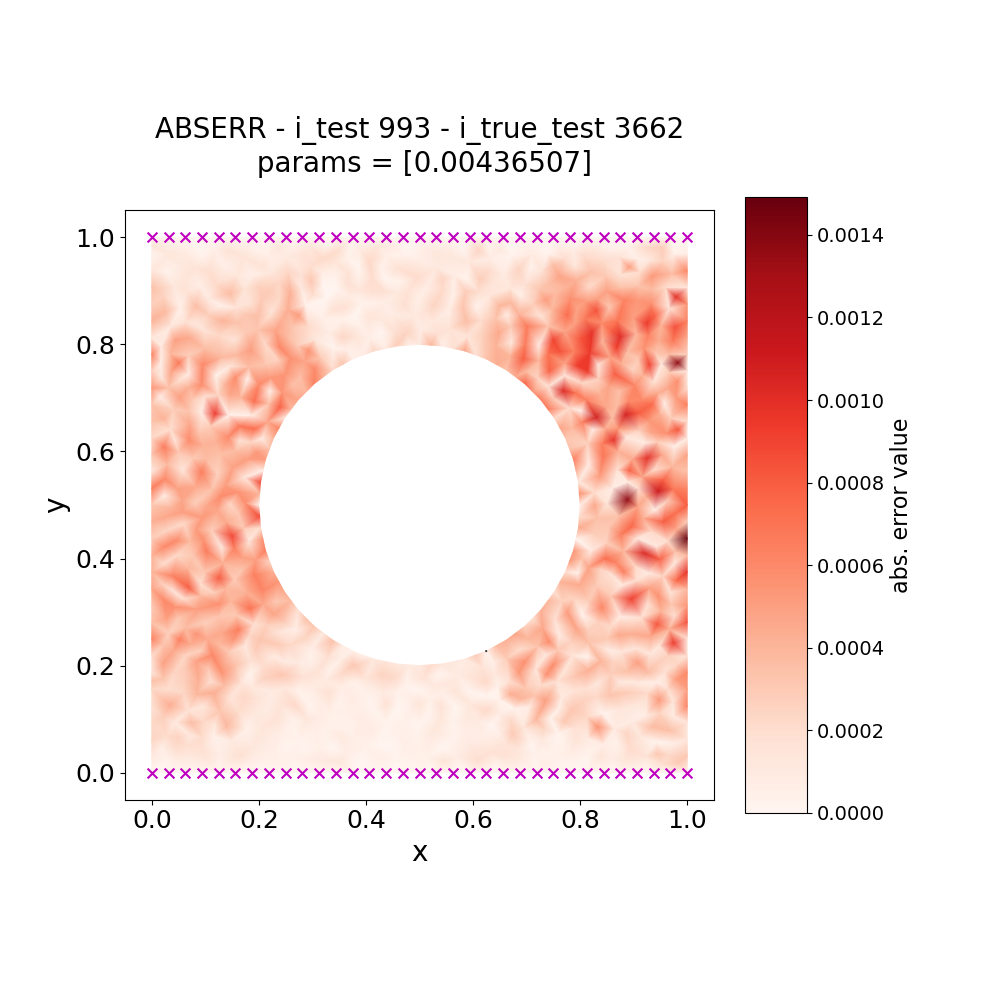}
    }
    \caption{
    Experiment 1 (Diffusion) - Best, median, and worst prediction cases with respect to the test set $\Xi$, ranked with respect to $\MLRE$. Black, empty circle dots denote the boundary nodes with Neumann BC of value $\mu_1\in(0, 1]$; magenta crosses denote boundary nodes with fixed homogeneous Dirichlet BCs.
    }
    \label{fig:examples_exp1}
\end{figure}

The model performance is evaluated by considering the average errors $\MlE$, $\MLE$, and $\MHE$. These errors, together with information about training epochs and training time, are reported in Table \ref{tab:results_exp1} and illustrated in Figure \ref{fig:results_exp1}. {In Figure \ref{fig:examples_exp1}, we report some prediction examples taken from the test set $\Xi$.} We do not report the relative errors for this test case because many solutions have nearly zero norm, making the relative errors meaningless {(e.g., see Figure \ref{fig:examples_exp1}, bottom row)}.

\begin{table}[htbp]
\centering
\begin{tabular}{l|l|l||r|r|r||r|r}
Model & n. weights & $T$ & $\MlE$ & $\MLE$ & $\MHE$ & batch tr.time (s) & tr. epochs \\
\hline
\hline
\multirow{4}{0pt}{\pFCNN} & \multirow{4}{0pt}{6.576e+06} & 128  & 8.593e-01 & 1.809e-02 & 1.215e+00 & 2.166e-02 & 471.4 
\\
 & & 256  & 7.692e-01 & 1.609e-02 & 1.103e+00 & 1.931e-02 & 622.4 
\\
 & & 512  & 7.643e-01 & 1.605e-02 & 1.083e+00 & 1.818e-02 & 580.0 
\\
 & & 1024 & 9.520e-01 & 2.009e-02 & 1.336e+00 & 1.747e-02 & 411.0 
\\
\hline
\multirow{4}{0pt}{\pGINN} & \multirow{4}{0pt}{1.898e+06} & 128  & 1.150e-01 & 2.325e-03 & 1.769e-01 & 1.250e+00 & 1745.4 
\\
 & & 256  & 6.973e-02 & 1.396e-03 & 1.091e-01 & 8.090e-01 & 1535.8 
\\
 & & 512  & 3.476e-02 & 6.800e-04 & 5.643e-02 & 6.684e-01 & 1462.6 
\\
 & & 1024 & 2.093e-02 & 4.040e-04 & 3.503e-02 & 5.393e-01 & 1742.2 \\
\end{tabular}
\caption{
Experiment 1 (Diffusion) - Average performance of the models, with respect to five random seeds for weight initialization, varying the training set size.
}
\label{tab:results_exp1}
\end{table}

Looking at the behaviors and values of the average errors of the trained models, the advantage of using a \pGINN model instead of a \pFCNN becomes evident. Indeed, \pGINN models exhibit better predictive performance on the test set; moreover, this performance is stable with respect to the random initialization of the weights. In contrast, \pFCNNs only rarely achieve slightly better results for some initializations (and still not better than those of \pGINNs).

In general, we observe that \pFCNNs suffer from poor generalization properties, even when the size of the training set is increased; indeed, the average value of the errors is almost constant. On the contrary, \pGINNs exhibit a clear error reduction as the training set size increases and, moreover, they achieve remarkably low error values even when only a few hundred training samples are used.

Regarding computational costs, we observe that \pGINNs have approximately one-sixth the number of trainable parameters of \pFCNNs (1.898e+06 weights instead of 6.576e+06, see Section \ref{sec:architectures} for details on the NN architectures); however, their training time is longer (approximately one order of magnitude greater), because the GI layers are based on sparse tensor/matrix operations, for a wide larger number of layers ($H=81$ GI layers for \pGINNs and $H=3$ FC layers for \pFCNNs, see Section \ref{sec:architectures}).

\subsection{Experiment 2 - Neumann and Dirichlet linear case (Advection-Diffusion)}\label{sec:exp2}

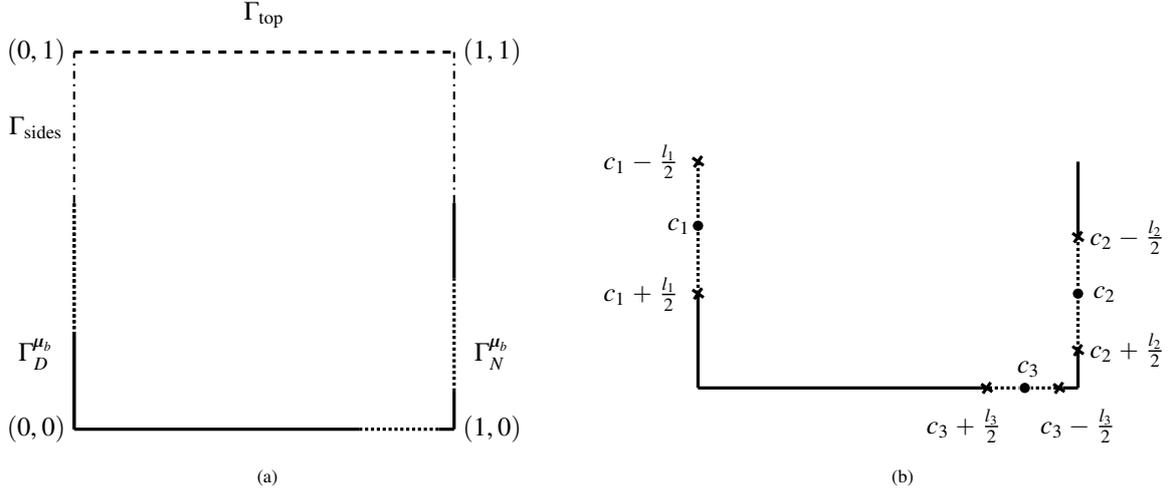
\begin{figure}[H]
    \begin{subfigure}[t]{0.5\textwidth}
        \begin{center}
\begin{tikzpicture}[scale=5]

\filldraw[color=black, very thick](0,0) -- (0.,0.25);
\filldraw[color=black, very thick, densely  dotted](0,0.25) -- (0.,.6);
\filldraw[color=black, thick, dashdotted](0,0.6) -- (0.,1.);
%\filldraw[color=black, fill=gray!10, very thick, dashed](0,1) -- (1,1);
%\filldraw[color=black, fill=gray!10, very thick](0.4,0.) -- (1,0.);
\filldraw[color=black, very thick](0,0) -- (0.,0.25);
\filldraw[color=black, very thick, densely  dotted](0,0.25) -- (0.,.6);
\filldraw[color=black, fill=gray!10, very thick, dashed](0,1) -- (1,1);
\filldraw[color=black, fill=gray!10, very thick](0,0.) -- (0.75,0.);
\filldraw[color=black, fill=gray!10, very thick,densely  dotted](0.75,0) -- (.96,0.);
\filldraw[color=black, fill=gray!10, very thick](0.96,0) -- (1,0.);
\filldraw[color=black, fill=gray!10, very thick](1,0.) -- (1,0.1);
\filldraw[color=black, fill=gray!10, very thick, densely dotted](1,0.1) -- (1,0.4);

%\filldraw[color=black, fill=gray!10, very thick](1,0.2) -- (1,0.1);
\filldraw[color=black, fill=gray!10, very thick](1,0.4) -- (1,0.6);

\filldraw[color=black, fill=gray!10,thick, dashdotted](1,0.6) -- (1,1);

\node at (1.1,0.2){\color{black}{$\Gamma_{N}^{\boldsymbol{\mu}_b}$}};

\node at (-.1,0.2){\color{black}{$\Gamma_D^{\boldsymbol{\mu}_b}$}};
\node at (-.1,0.8){\color{black}{$\Gamma_{\text{sides}}$}};
\node at (.5,1.1){\color{black}{$\Gamma_{\text{top}}$}};
\node at (-.1,0){\color{black}{$(0,0)$}};

\node at (-.1,1.){\color{black}{$(0,1)$}};
\node at (1.1,1.){\color{black}{$(1,1)$}};
\node at (1.1,-0.){\color{black}{$(1,0)$}};
%\node at (1.1,-.05){\color{black}{$ \color{cyan}{\Gamma_\text{w}}$}};

% \node at (-.5,0.5){\color{black}{$\Gamma_{In}^{\mu_D^*}$}};
% \node at (8.5,0.5){\color{black}{$\Gamma_{N}$}};

% \node at (2.5,-1.5){\color{black}{$\Gamma_{D}^{\mu_D^*}$}};

\end{tikzpicture}
\caption{}
\label{subfig:domain3}
\end{center}
    \end{subfigure}%
    ~ 
    \begin{subfigure}[t]{0.5\textwidth}
        \begin{center}
\tikzset{cross/.style={cross out, draw=black, minimum size=2*(#1-\pgflinewidth), inner sep=0pt, outer sep=0pt},
%default radius will be 1pt. 
cross/.default={3pt}}
\begin{tikzpicture}[scale=5]

\filldraw[color=black, very thick](0,0) -- (0.,0.25);
\filldraw[color=black, very thick, densely  dotted](0,0.25) -- (0.,.6);
\filldraw[color=white, very thick](0,0.6) -- (0.,1.);
\filldraw[color=white, fill=gray!10, very thick, dashed](0,1) -- (1,1);
\filldraw[color=black, fill=gray!10, very thick](0,0.) -- (0.75,0.);
\filldraw[color=black, fill=gray!10, very thick,densely  dotted](0.75,0) -- (.96,0.);
\filldraw[color=black, fill=gray!10, very thick](0.96,0) -- (1,0.);
\filldraw[color=black, fill=gray!10, very thick](1,0.) -- (1,0.1);
\filldraw[color=black, fill=gray!10, very thick, densely dotted](1,0.1) -- (1,0.4);

%\filldraw[color=black, fill=gray!10, very thick](1,0.2) -- (1,0.1);
\filldraw[color=black, fill=gray!10, very thick](1,0.4) -- (1,0.6);

\node at (1.2,0.2){\color{white}{$\partial \Omega \setminus \Gamma_{D}^{\mu_D^*}$}};

\node at (-.2,0.2){\color{white}{$\Gamma_D^{\mu_N^*}$}};
\node at (.5,1.1){\color{white}{$\Gamma_{\text{in}}$}};
\node at (-.1,0){\color{white}{$(0,0)$}};

\node at (-.1,1.){\color{white}{$(0,1)$}};
\node at (1.1,1.){\color{white}{$(1,1)$}};
\node at (1.1,-0.){\color{white}{$(1,0)$}};
\filldraw[color=white, fill=black](1,0.25)circle (0.015);
\draw (0.76,-0.) node[cross,rotate=10, very thick] {};
\draw (1,0.4) node[cross,rotate=10, very thick] {};
\draw (1,0.1) node[cross,rotate=10, very thick] {};
\draw (0.95,-0.) node[cross,rotate=10, very thick] {};

\node at (1.07,0.25){\color{black}{$c_2$}};
\node at (1.13,0.1){\color{black}{$c_2 + \frac{l_2}{2}$}};
\node at (1.13,0.4){\color{black}{$c_2 - \frac{l_2}{2}$}};
\node at (0.7,-0.1){\color{black}{$c_3 + \frac{l_3}{2}$}};
\node at (1,-0.1){\color{black}{$c_3 - \frac{l_3}{2}$}};
\node at (0.87,0.05){\color{black}{$c_3$}};
\filldraw[color=white, fill=black](0.86,0.)circle (0.015);
\filldraw[color=white, fill=black](1,0.25)circle (0.015);
\draw (0,0.6) node[cross,rotate=10, very thick] {};
\draw (0,0.25) node[cross,rotate=10, very thick] {};
\node at (-.05,0.43){\color{black}{$c_1$}};
\filldraw[color=white, fill=black](0.,0.43)circle (0.015);
\node at (-.15,0.6){\color{black}{$c_1 - \frac{l_1}{2}$}};
\node at (-.15,0.25){\color{black}{$c_1 + \frac{l_1}{2}$}};

\end{tikzpicture}
\caption{}
\label{subfig:interval_square}
\end{center}
    \end{subfigure}%
    \caption{(Experiments 2 and 3) Spatial domain $\Omega$ for a specific parametric instance: schematic representation. (A) ${\Gamma_{\text{sides}}}$ (dash-dotted line) and ${\Gamma_{\text{top}}}$ (dashed line) form the fixed boundary ${\Gamma_{\text{fix}}}$. The boundary parameter $\bmu_b$ characterizes $\Gamma_N^{\bmu_b}$ (dotted line) and $\Gamma_D^{\bmu_b} = \Gamma^{\bmu_b}\setminus \Gamma_N^{\bmu_b}$ (solid line), respectively. (B) Example of selection of circle intervals. The summation only represents the detection of the start and the end node of the intervals.}
\label{fig:domain3}
\end{figure}

We consider the physical domain $\Omega$ as the unit square. In this specific test case, the geometrical parameter affects both the Dirichlet and the Neumann boundaries. Specifically, we fix a boundary $\Gamma_{\text{fix}} = {\Gamma_{\text{sides}}} \cup {\Gamma_{\text{top}}}$, with ${\Gamma_{\text{sides}}} = (\{0\} \times [0.6, 1]) \cup  (\{1\} \times [0.6, 1])$ and ${\Gamma_{\text{top}}} = [0,1] \times \{1\}$. On $\Gamma_{\text{fix}}$, homogeneous Dirichlet boundary conditions are applied. In terms of parameters, we define $\bmu = (\bmu_{\phi}, \mu_b) = (\mu_1, \mu_2, \mu_b)$, where $\mu_1\in [10^{-3},100]$ denotes the diffusion parameter, $\mu_2\in [(-3/4)\pi,(-1/4) \pi]$ denotes the angle of the unitary advection vector $\v{\alpha}\in\R^2$ with respect to the horizontal axis, and $\mu_b$ defines $\Gamma^{\bmu_b}:=\partial\Omega \setminus \Gamma_{\text{fix}}$ and characterizes the portions where homogeneous Dirichlet BCs or homogeneous Neumann BCs hold (denoted as $\Gamma^{\bmu_b}_D$ and $\Gamma^{\bmu_b}_N$, respectively); for the ease of notation, we dropped the function $\mu_v$ in $\bmu$ because it is constantly equal to zero (only homogeneous BCs). See Figure \ref{subfig:domain3} for a visual example of the domain $\Omega$ and with BCs defined by $\Gamma_{\text{fix}}$ and $\Gamma^{\bmu_b}$, respectively.
Then, the problem we are solving is:
$$
\begin{cases}
    - \mu_1 \Delta u + {\v{\alpha}}(\mu_2) \cdot \nabla u = 1 & \text{ in } \Omega,\\
    \displaystyle \frac{\partial{u}}{\partial {\v{n}}} = {0} & \text{ on } \Gamma_N^{\bmu_b},\\
    \displaystyle u = 0 & \text{ on } \Gamma_D^{\bmu_b}\text{ and on } \Gamma_{\text{fix}}, \\
\end{cases}
$$
where ${\v{\alpha}}(\mu_2) = (\cos(\mu_2), \sin(\mu_2))$.

We now present the rationale behind the variability distribution of $\bmu=(\mu_1, \mu_2, \mu_b)$, specifically used for generating the training, validation, and test data.
Each $\mu_1$ and $\mu_2$ is just randomly generated with uniform distribution with respect to the ranges defined above; i.e., $\mu_1\sim\mathcal{U}(10^{-3}, 100)$ and $\mu_2\sim\mathcal{U}((-3/4)\pi, (-1/4)\pi)$.
Now, for generating $\mu_b$, we split $\Gamma^{\bmu_b}$ into three segments $\Gamma^{\bmu_b}_1=\{0\}\times[0, 0.6]$, $\Gamma^{\bmu_b}_2=\{1\}\times[0, 0.6]$, and $\Gamma^{\bmu_b}_3=[0, 1]\times\{0\}$; then, we randomly sample three points $\v{c}_1=(0, y_1), \v{c}_2=(1, y_2), \v{c}_3=(x_3, 0)$ on $\Gamma^{\bmu_b}_1$, $\Gamma^{\bmu_b}_2$, $\Gamma^{\bmu_b}_3$, respectively, and we generate three interval lengths $l_1,l_2,l_3$ such that $l_1,l_2\sim \mathcal{U}(0.4, 0.6)$ and $l_3\sim\mathcal{U}(2/3, 1)$. Therefore, for each sampling of $\bmu$, we define $\Gamma^{\bmu_b}_N:=I_1\cup I_2\cup I_3$ and $\Gamma^{\bmu_b}_D=\Gamma^{\bmu_b}\setminus\Gamma^{\bmu_b}_N$, where
\begin{equation*}
    \begin{aligned}
        I_1 &= \{0\}\times \left[\max
        \left \{ 0, y_1 - \frac{l_1}{2}\right \}, \min\left \{0.6, y_1 + \frac{l_1}{2} \right \} \right]\,,\\
        I_2 &= \{1\}\times \left[\max
        \left \{0, y_2 - \frac{l_2}{2} \right \}, \min
        \left \{0.6, y_2 + \frac{l_2}{2} \right \} \right]\,,\\
        I_3 &= \left[\max \left \{0, x_3 - \frac{l_3}{2} \right \}, \min \left \{1, x_3 + \frac{l_3}{2} \right \} \right]\times \{0\}\,.
    \end{aligned}
\end{equation*}

The above generation procedure for $\mu_b$ can be easily translated to a discrete representation of the domain, after the creation of the mesh, for generating a vector $\bmu_b$. Denoting with $\v{x}_1,\ldots ,\v{x}_M$ the $M$ points of the mesh on the boundary $\Gamma^{\bmu_b}$, the points $\v{c}_1,\v{c}_2,\v{c}_3$ are randomly sampled among them, while the interval lengths are measured in number of interval nodes. Therefore, in the discrete case (and in practice) the parameter $\bmu_b$ is a vector in $\{0, 1\}^M\subset\R^M$ (i.e., $p_b=M$) such that its elements are $\bmu_b\pe{j}=1$ if $\v{x}_j\in\Gamma_N^{\bmu_b}$, and $\bmu_b\pe{j}=0$ otherwise. See Figure \ref{subfig:interval_square} for a visual example. It is worth noting that this generation procedure can generate intervals connected at the corners.

\subsubsection{{Experiment 2 - Results}}\label{sec:results_exp2}

For this second experiment, we discretized the domain $\Omega$ with a mesh of $N_h=3993$ nodes, $p_b=256$ nodes on the boundary $\Gamma^{\bmu_b}$, generating a mesh graph $G$ with diameter $\mathrm{diam}(G)=74$. The dataset is generated by running $10\,000$ simulations through a $\mathbb{P}^1$ finite element discretization on the given mesh; in particular, each simulation is done with respect to a randomly generated parameter vector $\bmu=(\mu_1,\mu_2,\bmu_b)\in[10^{-3}, 100] \times [(-3/4)\pi,(-1/4)\pi] \times \{0, 1\}^{p_b}$. Given these simulation data, 2048 of them are selected randomly and used as a fixed test set; from the remaining data, we randomly extract an even number $T\in\N$ of training set data and $T/2$ validation set data. Each model (\pGINN and \pFCNN) is trained five times, with respect to five different random seeds for weight initialization, varying the amount of training data $T = 1024, 512, 256, 128$.

\begin{figure}[htbp!]
    \centering
    \subcaptionbox{legend}{\includegraphics[width=0.45\textwidth]{Figures/legend.png}
    }
    \\
    \subcaptionbox{$\MlE$}{\includegraphics[trim=0.5cm 0.0cm 2.5cm 1.8cm,clip,width=0.45\textwidth]{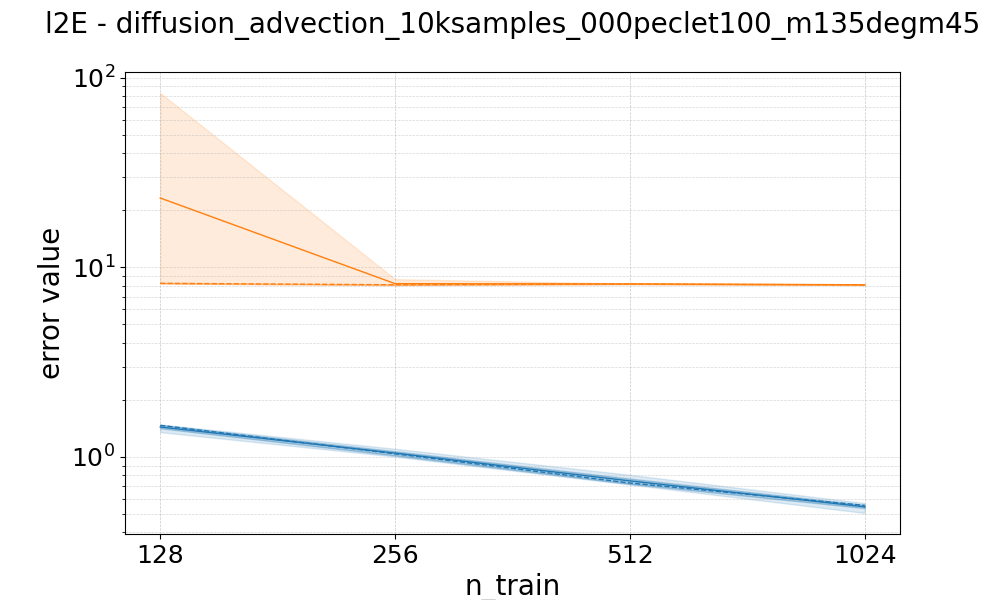}
    }
    \subcaptionbox{$\MlRE$}{\includegraphics[trim=0.5cm 0.0cm 2.5cm 1.8cm,clip,width=0.45\textwidth]{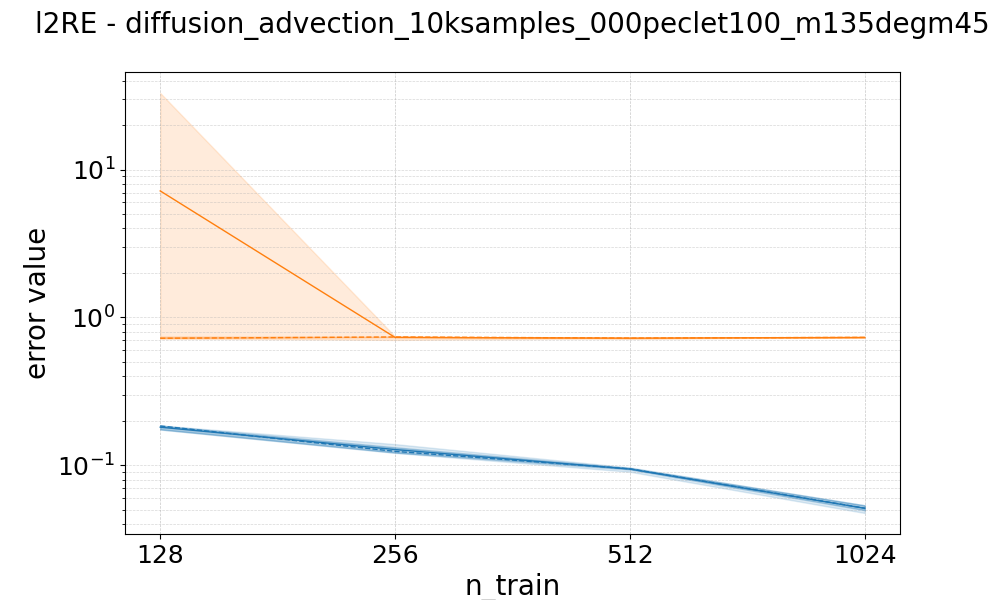}
    }
    \\
	\subcaptionbox{$\MLE$}{\includegraphics[trim=0.5cm 0.0cm 2.5cm 1.8cm,clip,width=0.45\textwidth]{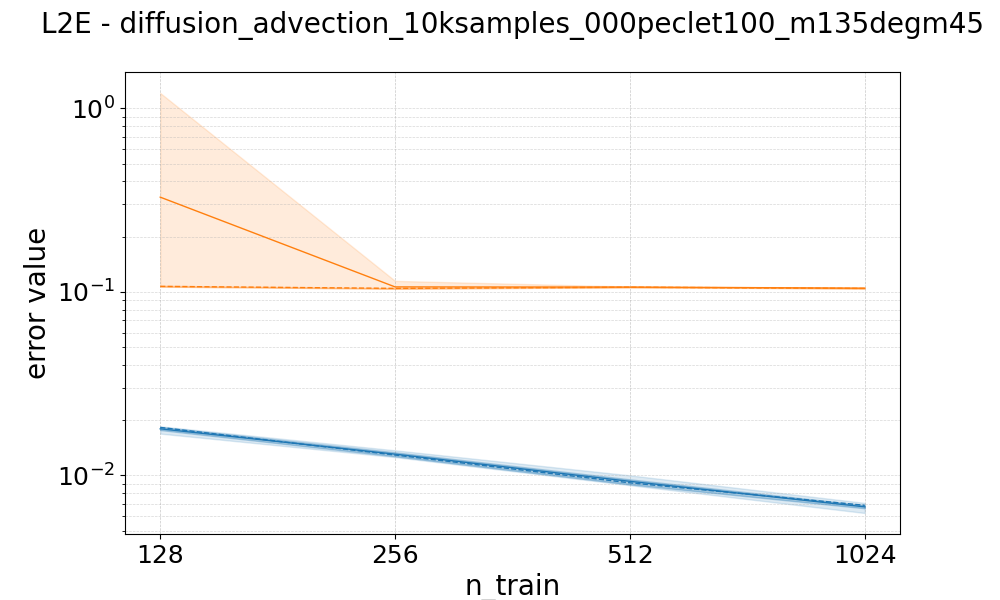}
    }
    \subcaptionbox{$\MLRE$}{\includegraphics[trim=0.5cm 0.0cm 2.5cm 1.8cm,clip,width=0.45\textwidth]{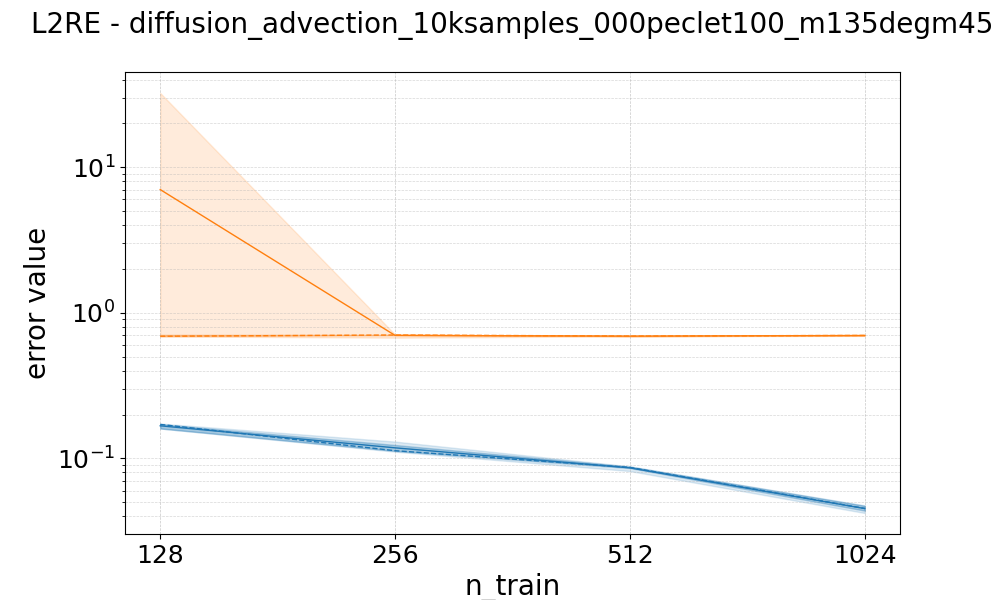}
    }
    \\
    \subcaptionbox{$\MHE$}{\includegraphics[trim=0.5cm 0.0cm 2.5cm 1.8cm,clip,width=0.45\textwidth]{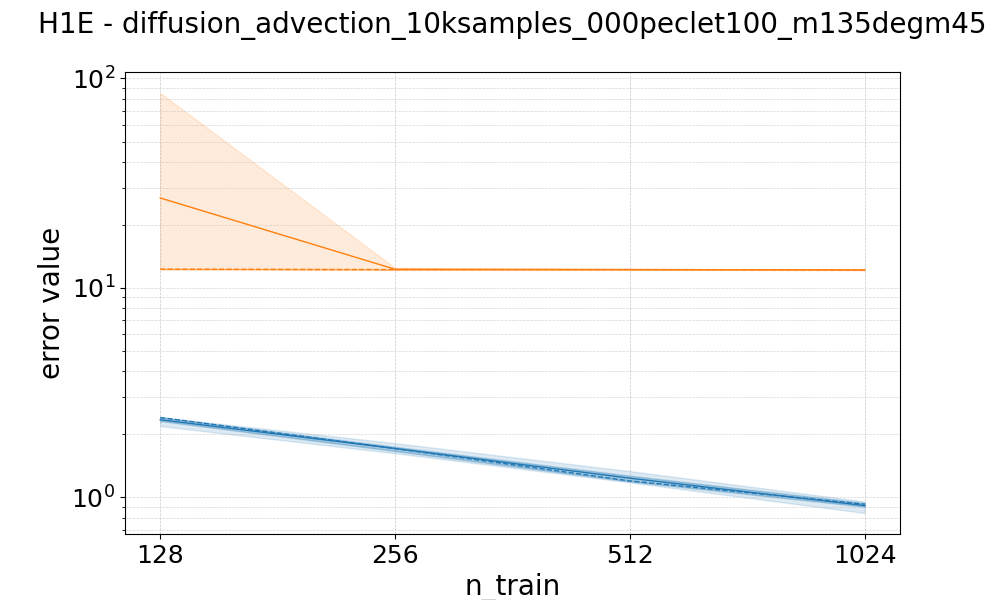}
    }
    \subcaptionbox{$\MHRE$}{\includegraphics[trim=0.5cm 0.0cm 2.5cm 1.8cm,clip,width=0.45\textwidth]{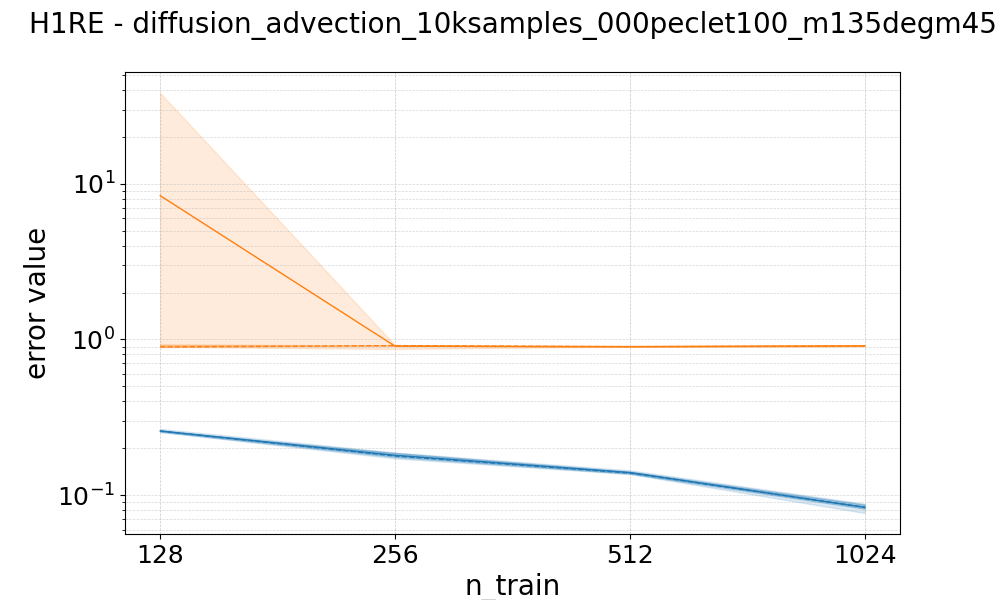}
    }
    \caption{
    Experiment 2 (Advection-Diffusion) - Error statistics with respect to five random seeds for weight initialization, varying the training set size. In blue, the \pGINN performance, in orange the \pFCNN performance. Light colored areas represent the Min-Max range of values, the dark colored areas represent the values between fist and third quartiles. The continuous lines represent the average errors (same as Table \ref{tab:results_exp2}), the dotted lines represent the medians.
    }
    \label{fig:results_exp2}
\end{figure}

\begin{figure}[htb!]
    \centering
    \subcaptionbox{Best case - $\solh$}{\includegraphics[trim=1.95cm 4.cm 7cm 5cm,clip,width=0.28\textwidth]{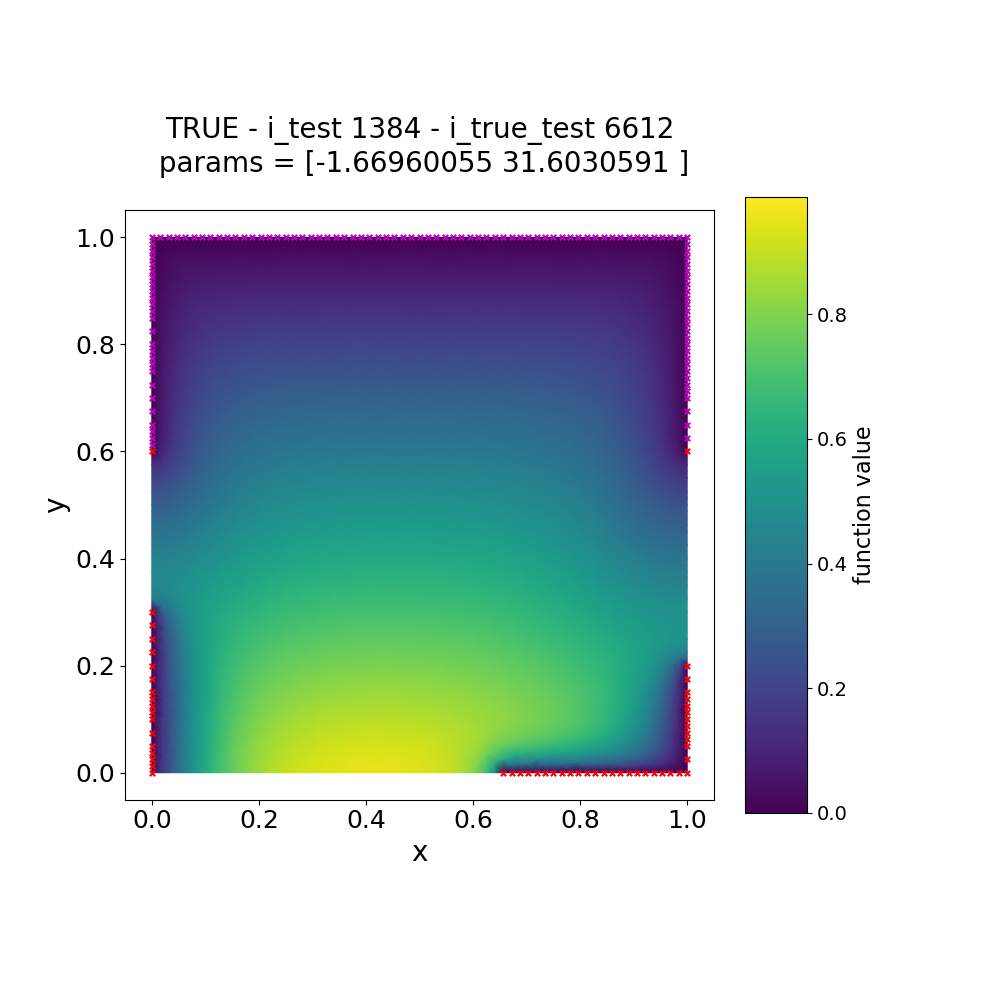}
    }
    \subcaptionbox{Best case - $\nnsolh$}{\includegraphics[trim=1.95cm 4.cm 3.75cm 4.75cm,clip,width=0.335\textwidth]{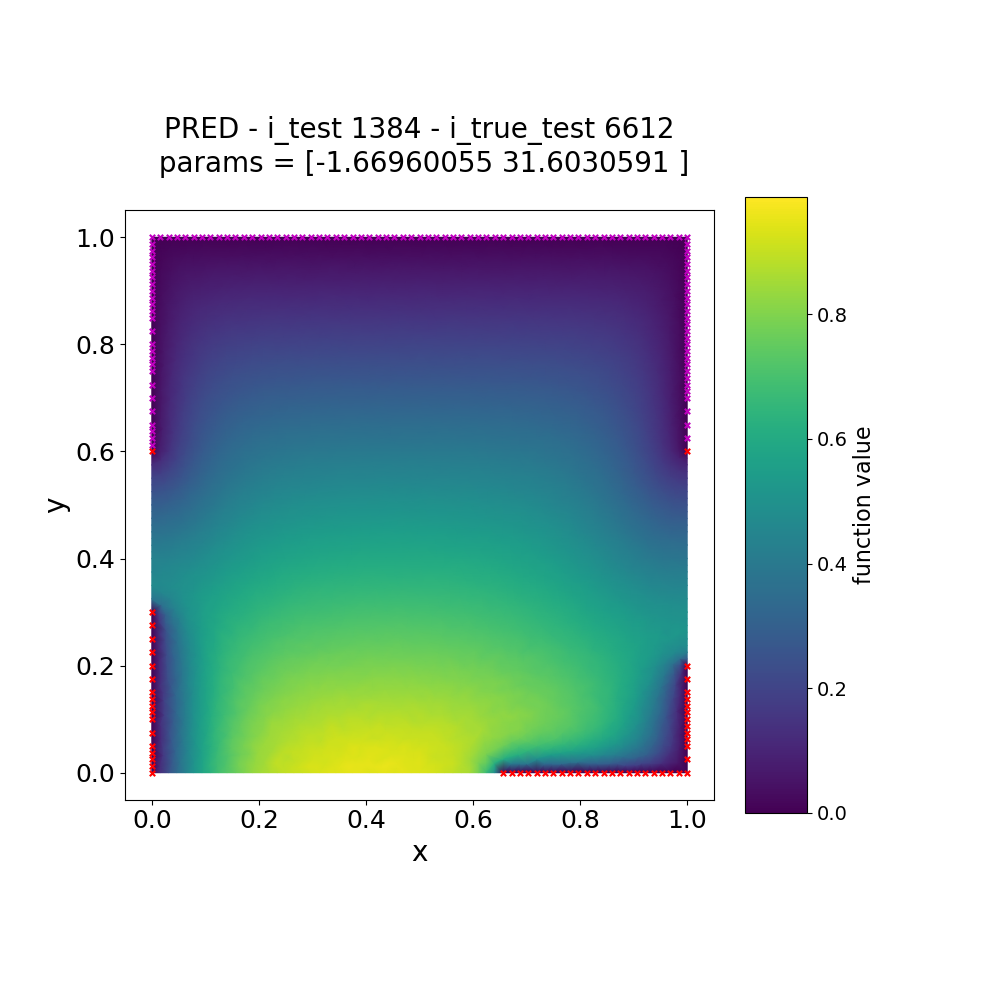}
    }
    \subcaptionbox{Best case - $|\solh - \nnsolh|$}{\includegraphics[trim=1.95cm 4.cm 3.25cm 4.75cm,clip,width=0.345\textwidth]{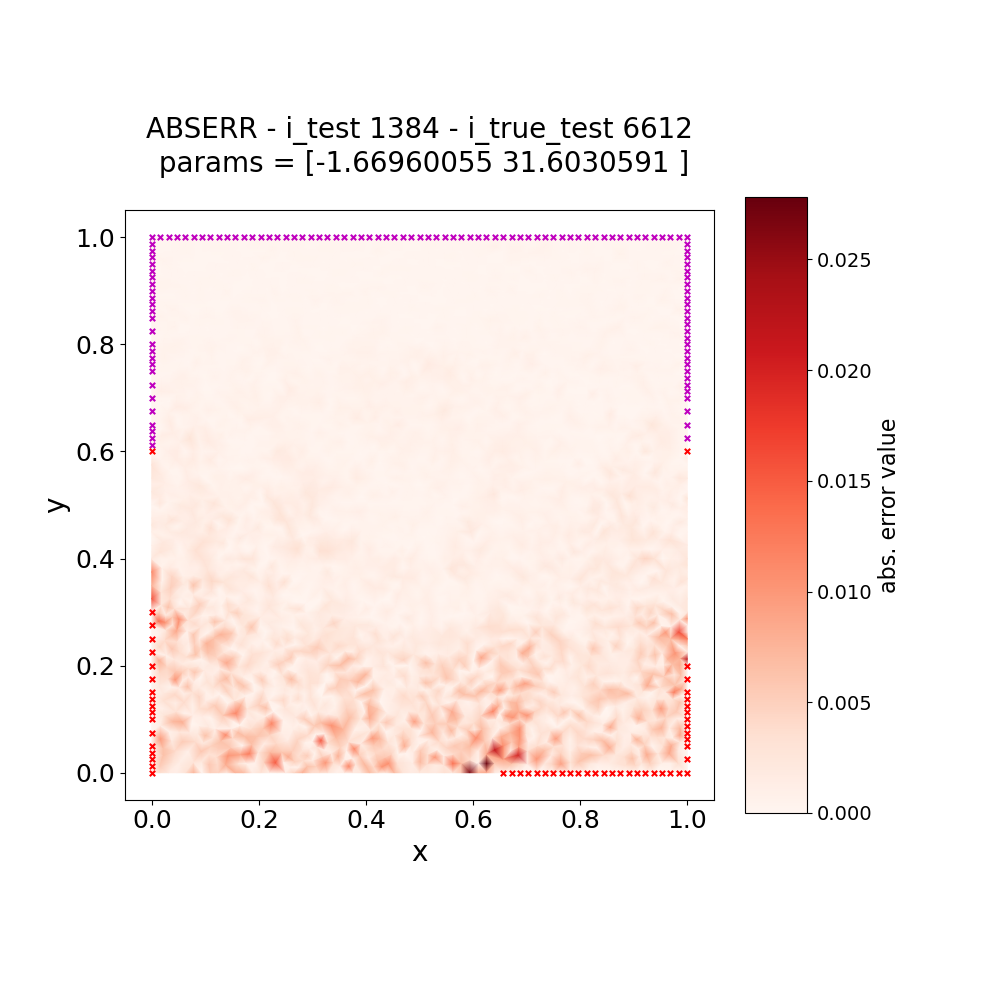}
    }
    \\
    \subcaptionbox{Median case - $\solh$}{\includegraphics[trim=1.95cm 4.cm 7cm 5cm,clip,width=0.28\textwidth]{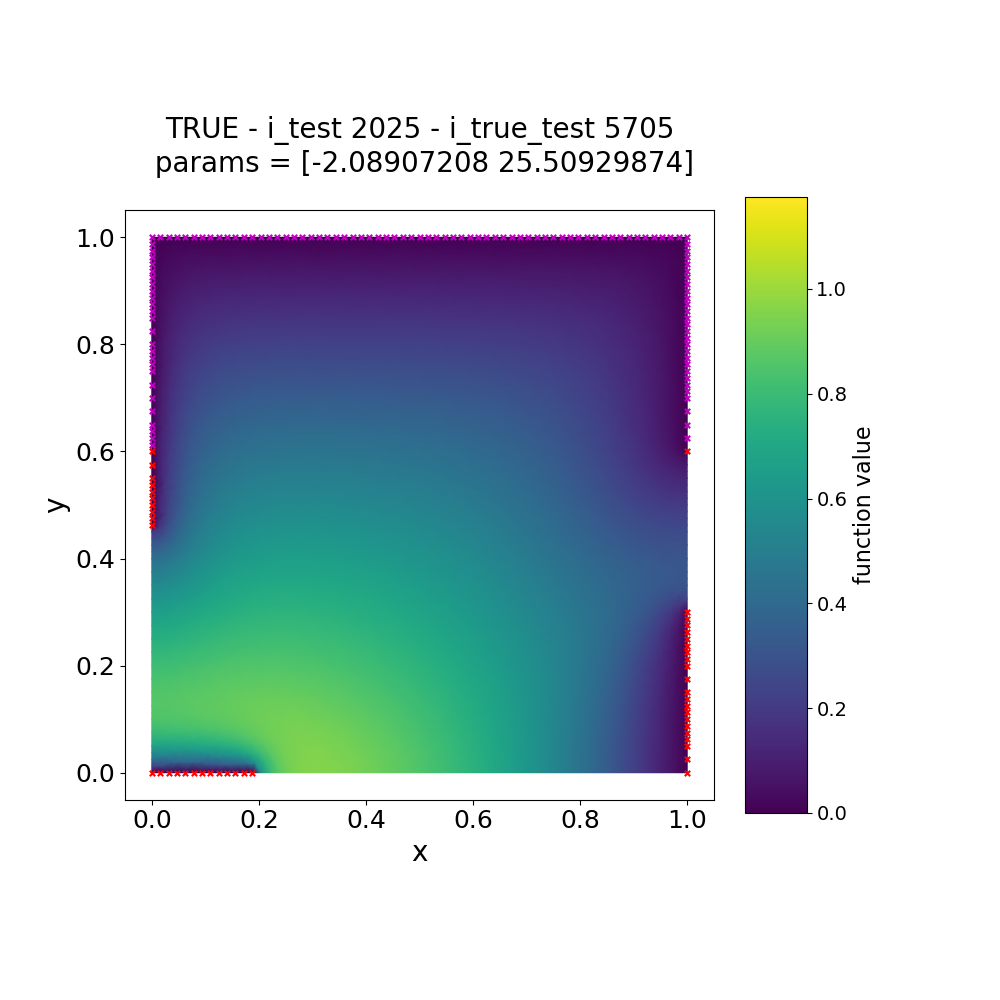}
    }
    \subcaptionbox{Median case - $\nnsolh$}{\includegraphics[trim=1.95cm 4.cm 3.75cm 4.75cm,clip,width=0.335\textwidth]{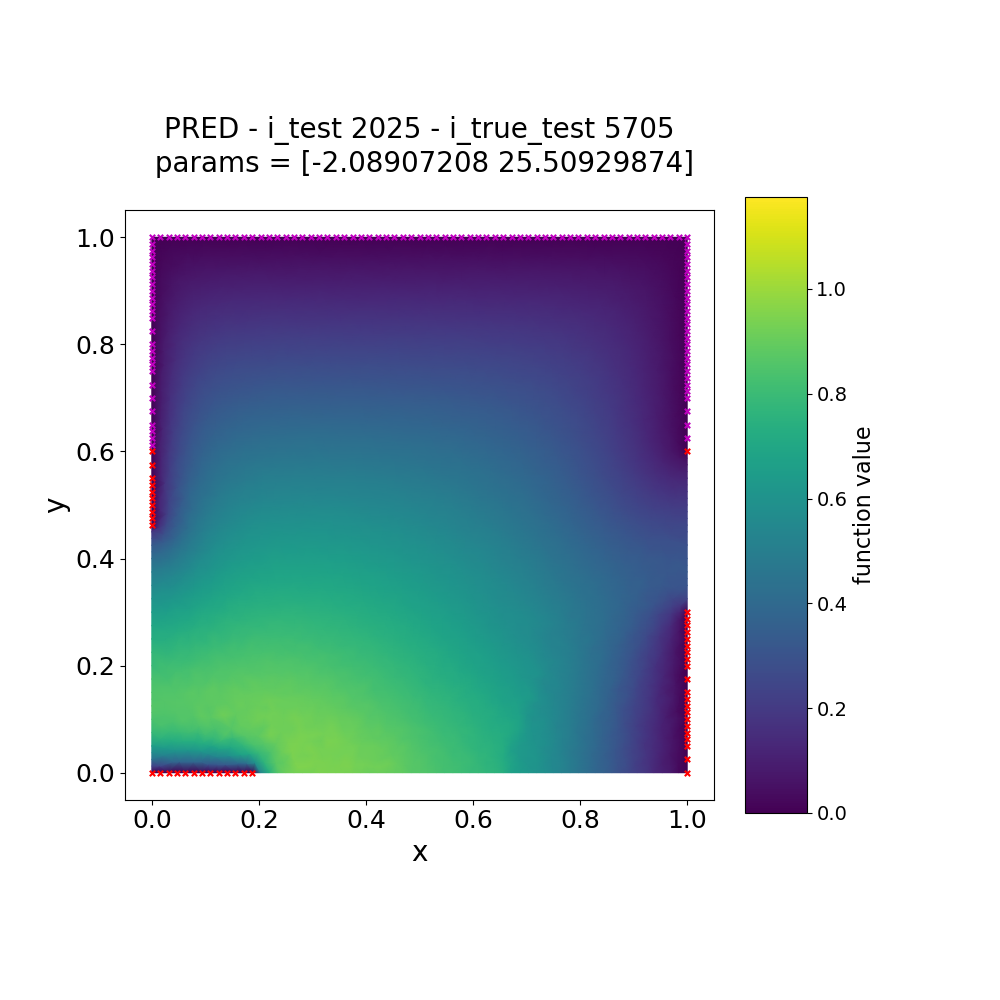}
    }
    \subcaptionbox{Median case - $|\solh - \nnsolh|$}{\includegraphics[trim=1.95cm 4.cm 3.35cm 4.75cm,clip,width=0.345\textwidth]{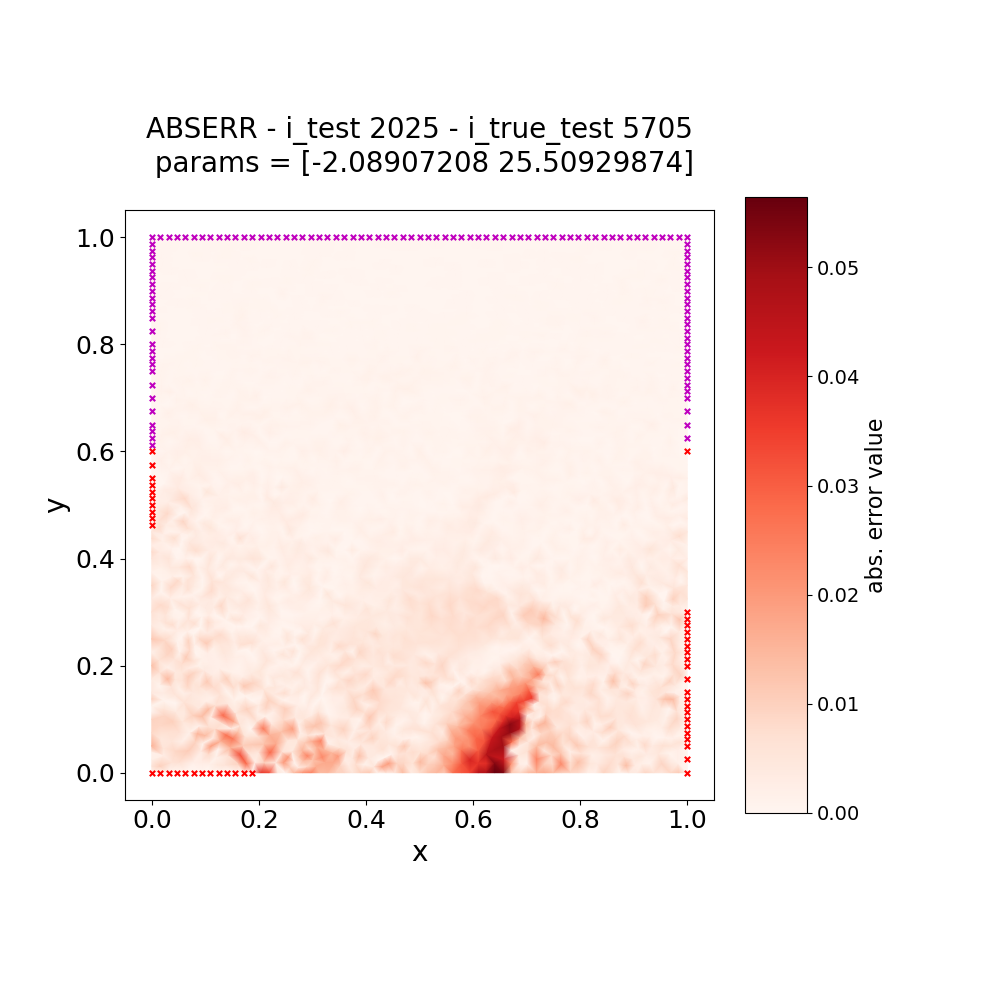}
    }
    \\
    \subcaptionbox{Worst case - $\solh$}{\includegraphics[trim=1.95cm 4.cm 7cm 5cm,clip,width=0.28\textwidth]{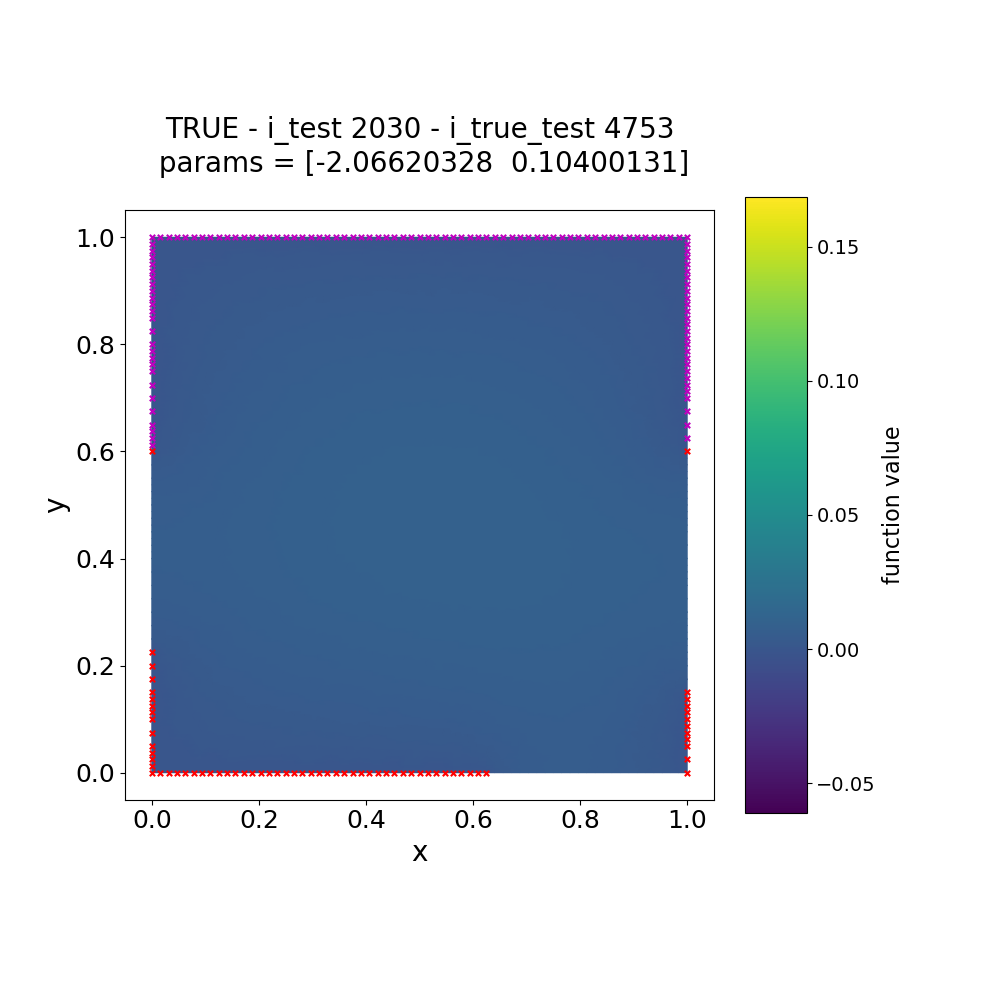}
    }
    \subcaptionbox{Worst case - $\nnsolh$}{\includegraphics[trim=1.95cm 4.cm 3.55cm 4.75cm,clip,width=0.335\textwidth]{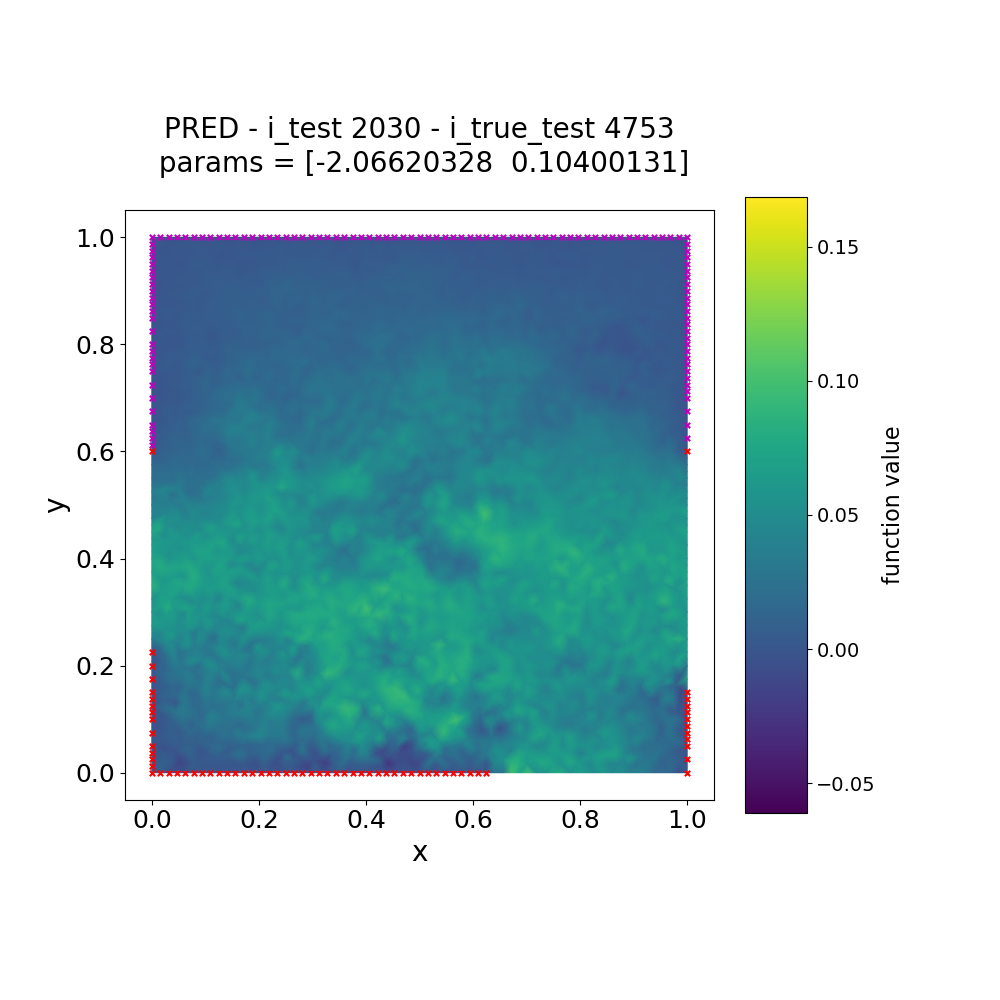}
    }
    \subcaptionbox{Worst case - $|\solh - \nnsolh|$}{\includegraphics[trim=1.95cm 4.cm 3.35cm 4.75cm,clip,width=0.345\textwidth]{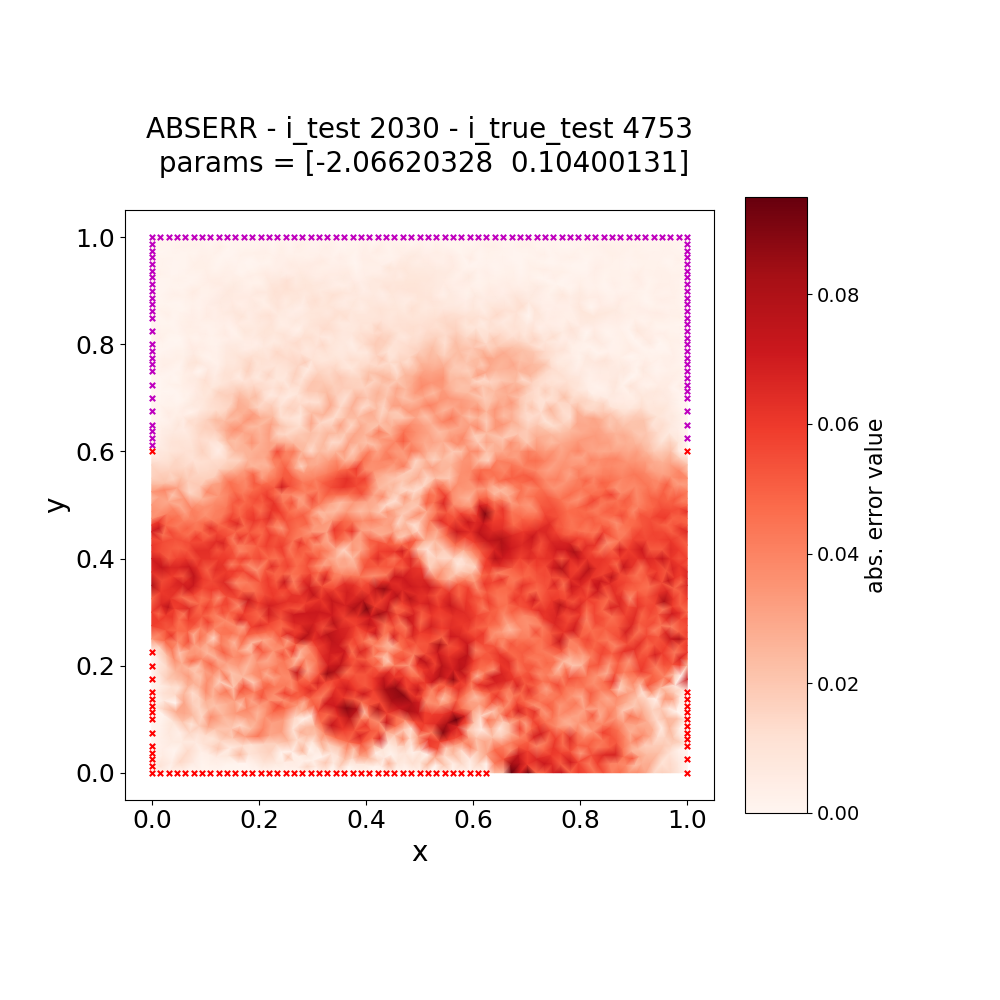}
    }
    \caption{
    Experiment 2 (Advection-Diffusion) - Best, median, and worst prediction cases with respect to the test set $\Xi$, ranked with respect to $\MLRE$. Magenta dots denote fixed homogeneous Dirichlet BCs, while red dots denote the homogeneous Dirichlet BC defined by the boundary parametrization of the problem.
    }
    \label{fig:examples_exp2}
\end{figure}

The model performance is evaluated by considering the average errors illustrated in \eqref{eq:err_mean}. These errors, together with information about training epochs and training time, are reported in Table \ref{tab:results_exp2} and illustrated in Figure \ref{fig:results_exp2}. In Figure \ref{fig:examples_exp2}, we report some prediction examples taken from the test set $\Xi$.

\begin{table}[htbp]
\centering
\resizebox{0.99\textwidth}{!}{
{
\begin{tabular}{l|l|l||r|r|r||r|r|r||r|r}
Model & n. weights & $T$ 
& $\MlE$ & $\MLE$ & $\MHE$ 
& $\MlRE$ & $\MLRE$ & $\MHRE$ 
& batch tr.time (s) & tr. epochs \\
\hline
\hline
\multirow{4}{0pt}{\pFCNN} & \multirow{4}{0pt}{8.017e+07} & 128  
& 2.322e+01 & 3.283e-01 & 2.689e+01 
& 7.179e+00 & 7.018e+00 & 8.385e+00 
& 2.235e+01 & 431.0 \\
 & & 256  
& 8.206e+00 & 1.065e-01 & 1.227e+01 
& 7.329e-01 & 6.991e-01 & 9.041e-01 
& 2.232e+01 & 460.8 \\
 & & 512  
& 8.180e+00 & 1.063e-01 & 1.223e+01 
& 7.251e-01 & 6.913e-01 & 8.985e-01 
& 2.234e+01 & 632.8 \\
 & & 1024 
& 8.096e+00 & 1.048e-01 & 1.219e+01 
& 7.327e-01 & 6.981e-01 & 9.092e-01 
& 2.228e+01 & 628.4 \\
\hline
\multirow{4}{0pt}{\pGINN} & \multirow{4}{0pt}{9.248e+06} & 128  
& 1.439e+00 & 1.795e-02 & 2.347e+00 
& 1.812e-01 & 1.674e-01 & 2.582e-01 
& 5.835e+00 & 1161.2 \\
 & & 256  
& 1.048e+00 & 1.306e-02 & 1.711e+00 
& 1.279e-01 & 1.182e-01 & 1.800e-01 
& 4.015e+00 & 1152.0 \\
 & & 512  
& 7.477e-01 & 9.286e-03 & 1.234e+00 
& 9.436e-02 & 8.607e-02 & 1.395e-01 
& 3.040e+00 & 1079.2 \\
 & & 1024 
& 5.458e-01 & 6.741e-03 & 9.134e-01 
& 5.121e-02 & 4.522e-02 & 8.327e-02 
& 2.557e+00 & 1244.4 \\
\end{tabular}
}
}
\caption{
Experiment 2 (Advection-Diffusion) - Average performance of the models, with respect to five random seeds for weight initialization, varying the training set size.
}
\label{tab:results_exp2}
\end{table}

Also in this second experiment, looking at the behaviors and values of the average errors, it is evident the advantage of using a \pGINN model instead of a \pFCNN one. Indeed, also in this case the \pGINN models exhibit better predictive performance on the test set, that are also stable with respect to the random initialization of the weights. In contrast, \pFCNNs show poor prediction performance, in all the configurations analyzed.

In general, we observe that \pFCNNs suffer from poor generalization properties, with an average value of the errors that is almost constant, independently on the training set size. On the contrary, \pGINNs (like in Experiment 1) exhibit a clear error reduction as the training set size increases, achieving low error values even when only a few hundred training samples are used.

Regarding computational costs, we observe that \pGINNs have a number of trainable parameters that is one order of magnitude smaller than the number of weights of \pFCNNs (9.248e+06 weights instead of 8.017e+07); moreover, in this experiment we notice a completely opposite behavior with respect to Experiment 1, where the training time of \pGINNs is remarkably shorter than the one of \pFCNNs (seconds instead of tens of seconds per mini-batch). The reason of this opposite behavior for the training times is the larger number of mesh nodes; indeed, in this experiment we have $N_h=3993$ nodes, while in Experiment 1 the nodes were $1141$. Therefore, \pFCNN models requires much more computational resources to manage the training of FC layers made of $N_h$ units, because FC layers scale much worse than GI layers (both in number of weights and in training time) with respect to this quantity.

\subsection{Experiment 3 - Neumann and Dirichlet nonlinear case (Navier-Stokes)}

We consider steady and incompressible Navier-Stokes equations for the domain $\Omega = [0,1]^2$ with boundary conditions defined and parametrized as the domain of the second experiment (see Section \ref{sec:exp2} and Figure \ref{fig:domain3}), except for the different values of the BCs. Specifically, the problem reads as
\begin{equation}
\label{eq:NS_eq}
\begin{cases}
-\displaystyle \frac{1}{\mu_{\phi}} \Delta \v{u} + (\v{u}\cdot\nabla) \v{u} + \nabla p = \v{0} \quad &\text{in} \ \Omega, \\
\nabla \cdot \v{u} = 0 \quad &\text{in} \ \Omega, \\
\v{u} = \v{u}_{\text{in}} \quad &\text{on} \ \Gamma_{\text{top}}, \\
\v{u} = \v{0} \quad &\text{on} \ \Gamma_{\text{sides}}\text{ and on } \Gamma_{D}^{\bmu_b}, \\
\displaystyle - p\v{n} + \frac{1}{\mu_1} \frac{\partial \v{u}}{\partial\v{n}} = 0 \quad &\text{on} \ \Gamma_{N}^{\bmu_b},
\end{cases}
\end{equation}
where $\mu_\phi\in\R$ denotes the Reynolds number ($Re$), $\v{u}$ and $p$ are the velocity and the pressure fluid fields, respectively, and 
\begin{equation}
\label{eq:inlet}
\v{u}_{\text{in}}(x_1, x_2) = \begin{cases}
\displaystyle \left  ( \frac{x_1}{0.06}, 0 \right) & \text{ for $x_1 \leq 0.06$ in } \Gamma_{\text{top}}, \\
\displaystyle \left  ( \frac{1 - x_1}{0.06}, 0 \right) & \text{ for $x_1 \geq 0.94$ in } \Gamma_{\text{top}}, \\
\displaystyle (1, 0) & \text{ elsewhere in } \Gamma_{\text{top}}.
\end{cases}
\end{equation}

For the ease of notation, we dropped the function $\mu_v$ in $\bmu$ because we have homogeneous BCs only.

Concerning the rationale behind the $\bmu=(\mu_{\phi}, \mu_b)$ variability distribution used for generating the training, validation, and test data, we have that $\mu_b$ is generated as the $\mu_b$ of the previous experiment (see Section \ref{sec:exp2}), while the parameter $\mu_{\phi}$ is such that $Re$ is randomly chosen between 100 and 500; i.e., $\mu_{\phi}\sim \mathcal{U}(100, 500)$.

\subsubsection{{Experiment 3 - Results}}\label{sec:results_exp3}

For this third experiment, we discretized the domain $\Omega$ with a mesh of $N_h=1444$ nodes, $p_b=144$ nodes on the boundary $\Gamma^{\bmu_b}$, generating a mesh graph $G$ with diameter $\mathrm{diam}(G)=47$. The dataset is generated by running $10\,000$ simulations through a $\mathbb P^2-\mathbb P^1$ finite element discretization on the given mesh; in particular, each simulation is done with respect to a randomly generated parameter vector $\bmu=(\mu_{\phi},\bmu_b)\in[100, 500] \times \{0, 1\}^{p_b}$. Among the $10\,000$ simulations, we removed the ones that did not reach convergence (due to limit-case scenarios randomly generated), obtaining a total of $9\,118$ simulations. Given these simulation data, 2048 of them are selected randomly and used as a fixed test set; from the remaining data, we randomly extract an even number $T\in\N$ of training set data and $T/2$ validation set data. Each model (\pGINN and \pFCNN) is trained five times, with respect to five different random seeds for weight initialization, varying the amount of training data $T = 1024, 512, 256, 128$.

\begin{figure}[htbp!]
    \centering
    \subcaptionbox{legend}{\includegraphics[width=0.45\textwidth]{Figures/legend.png}
    }
    \\
    \subcaptionbox{$\MlE$}{\includegraphics[trim=0.5cm 0.0cm 2.5cm 1.8cm,clip,width=0.45\textwidth]{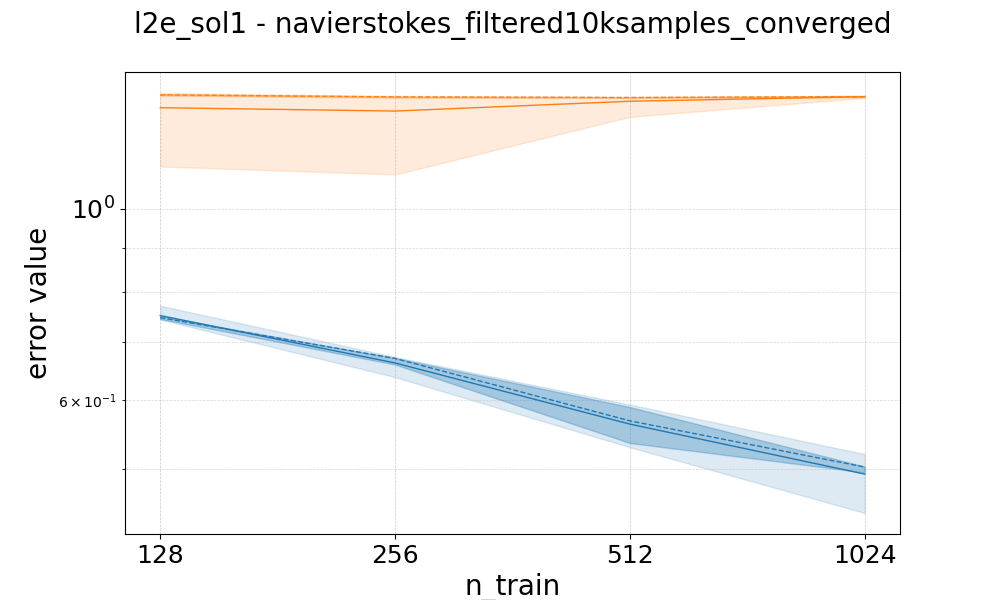}
    }
    \subcaptionbox{$\MlRE$}{\includegraphics[trim=0.5cm 0.0cm 2.5cm 1.8cm,clip,width=0.45\textwidth]{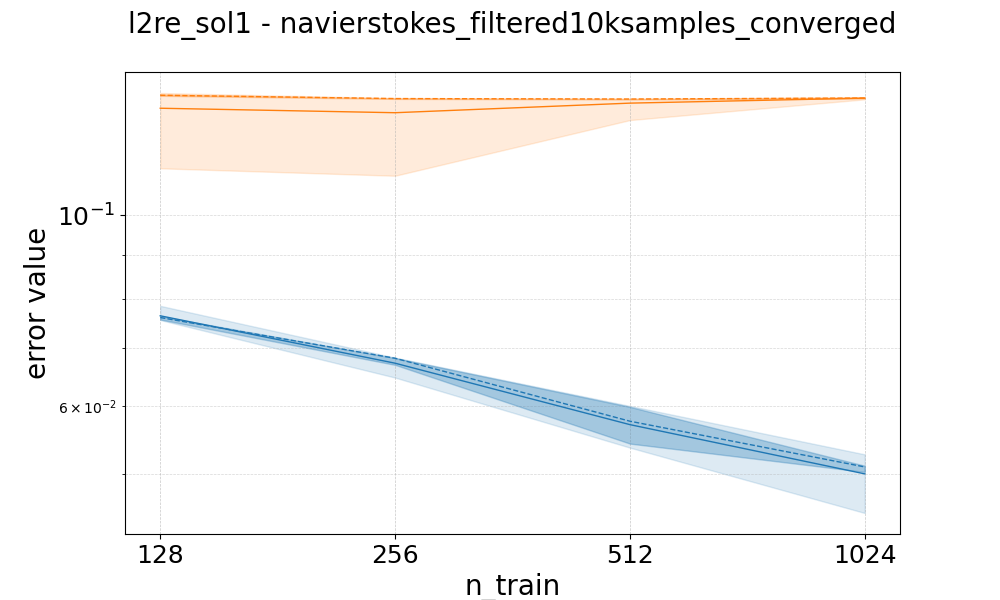}
    }
    \\
	\subcaptionbox{$\MLE$}{\includegraphics[trim=0.5cm 0.0cm 2.5cm 1.8cm,clip,width=0.45\textwidth]{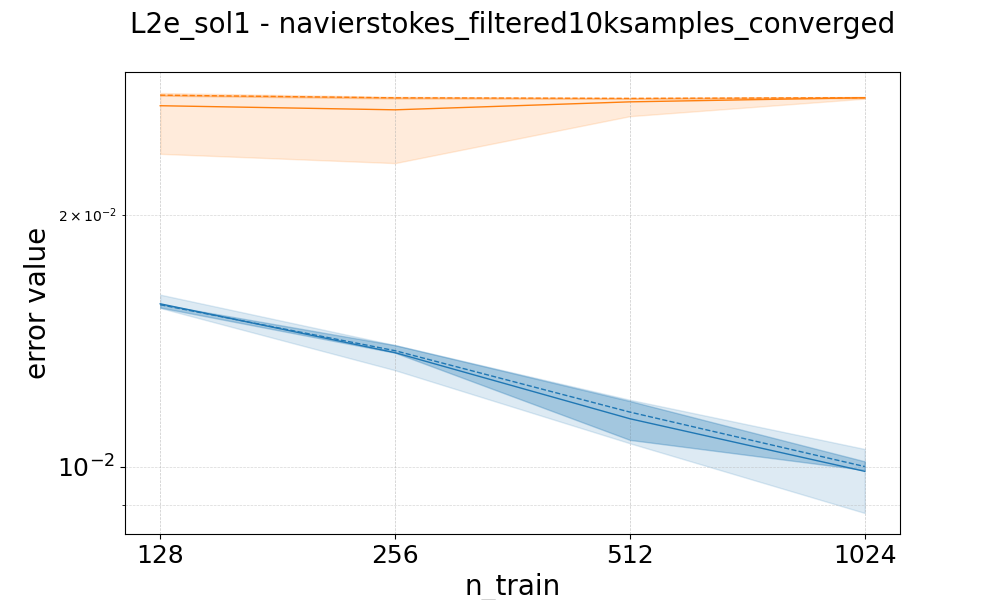}
    }
    \subcaptionbox{$\MLRE$}{\includegraphics[trim=0.5cm 0.0cm 2.5cm 1.8cm,clip,width=0.45\textwidth]{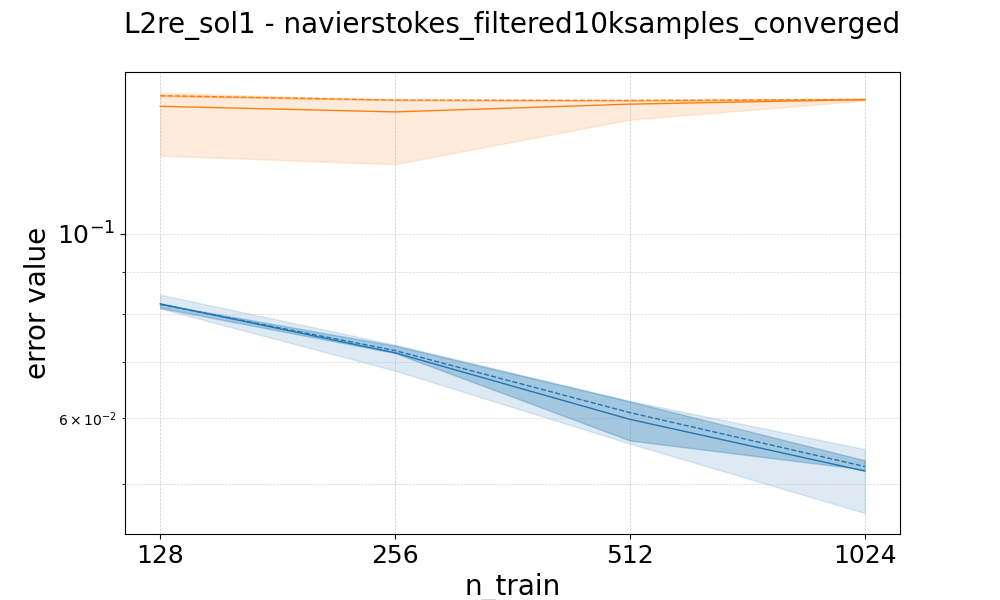}
    }
    \\
    \subcaptionbox{$\MHE$}{\includegraphics[trim=0.5cm 0.0cm 2.5cm 1.8cm,clip,width=0.45\textwidth]{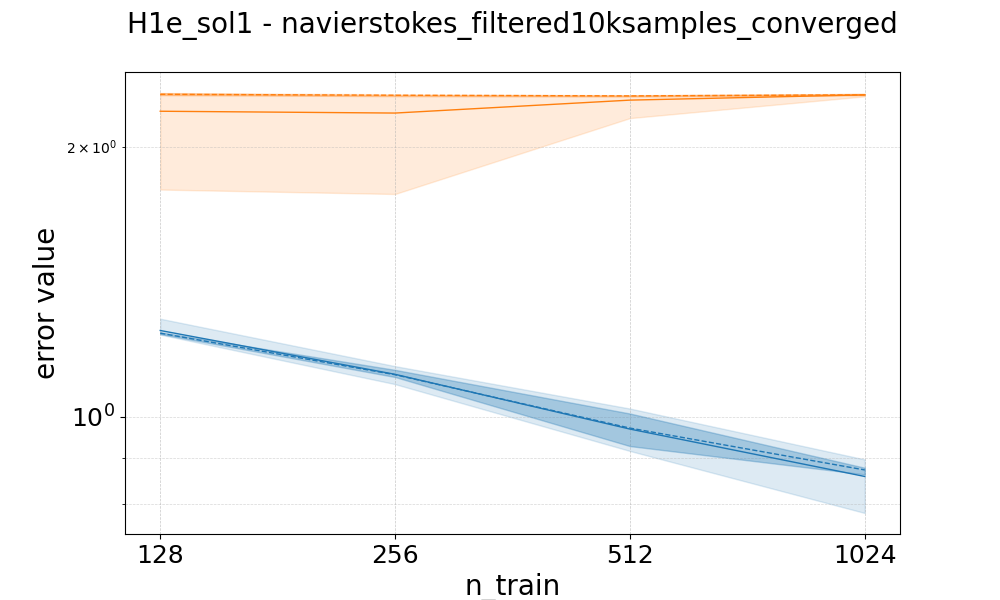}
    }
    \subcaptionbox{$\MHRE$}{\includegraphics[trim=0.5cm 0.0cm 2.5cm 1.8cm,clip,width=0.45\textwidth]{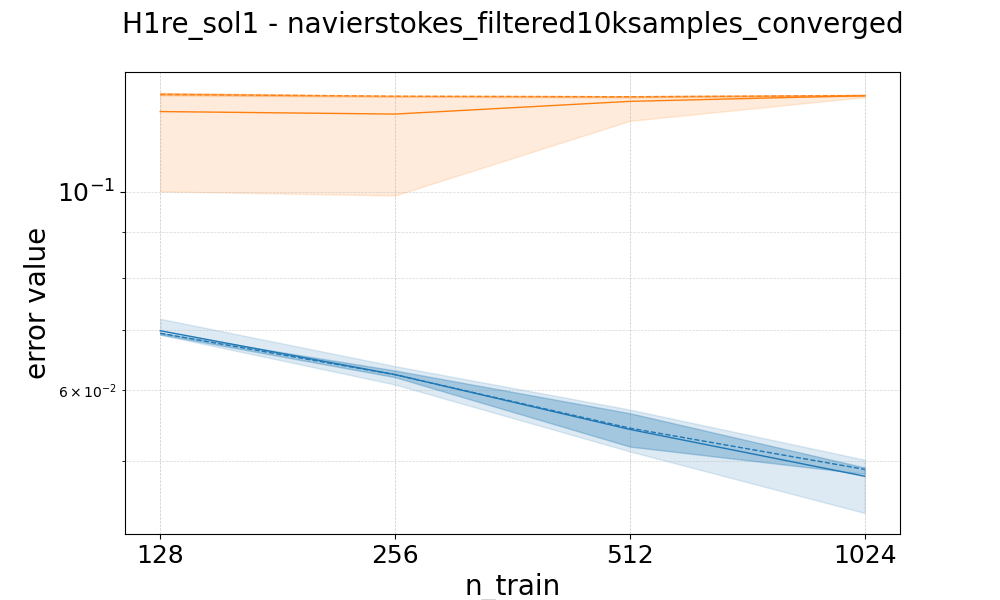}
    }
    \caption{
    Experiment 3 (Navier-Stokes) - Error statistics with respect to five random seeds for weight initialization, varying the training set size. In blue, the \pGINN performance, in orange the \pFCNN performance. Light colored areas represent the Min-Max range of values, the dark colored areas represent the values between fist and third quartiles. The continuous lines represent the average errors (same as Table \ref{tab:results_exp3_sol1}), the dotted lines represent the medians. Results for the first component of solution $\v{u}$.
    }
    \label{fig:results_exp3_sol1}
\end{figure}

\begin{figure}[htbp!]
    \centering
    \subcaptionbox{legend}{\includegraphics[width=0.45\textwidth]{Figures/legend.png}
    }
    \\
    \subcaptionbox{$\MlE$}{\includegraphics[trim=0.5cm 0.0cm 2.5cm 1.8cm,clip,width=0.45\textwidth]{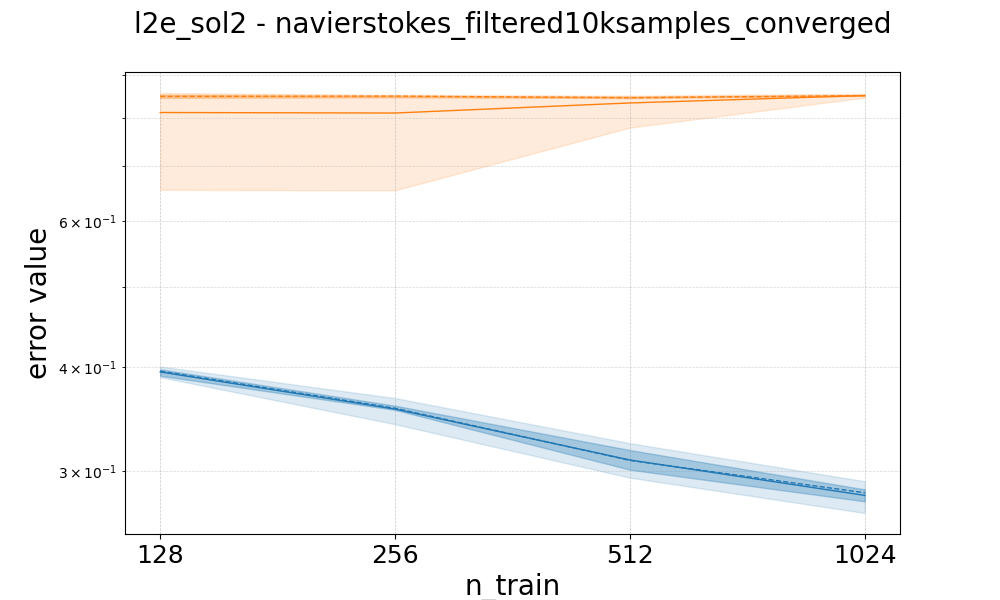}
    }
    \subcaptionbox{$\MlRE$}{\includegraphics[trim=0.5cm 0.0cm 2.5cm 1.8cm,clip,width=0.45\textwidth]{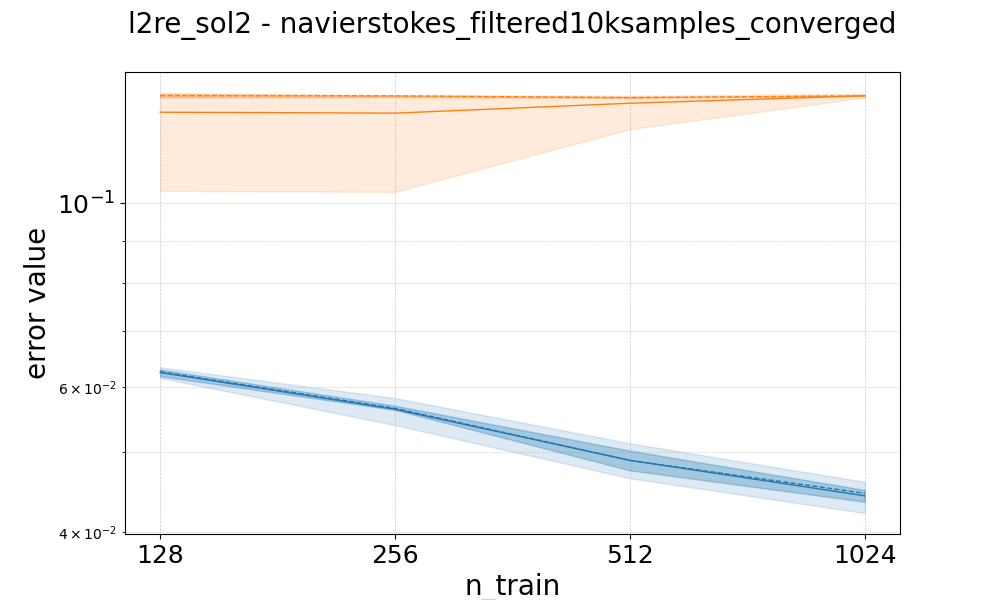}
    }
    \\
	\subcaptionbox{$\MLE$}{\includegraphics[trim=0.5cm 0.0cm 2.5cm 1.8cm,clip,width=0.45\textwidth]{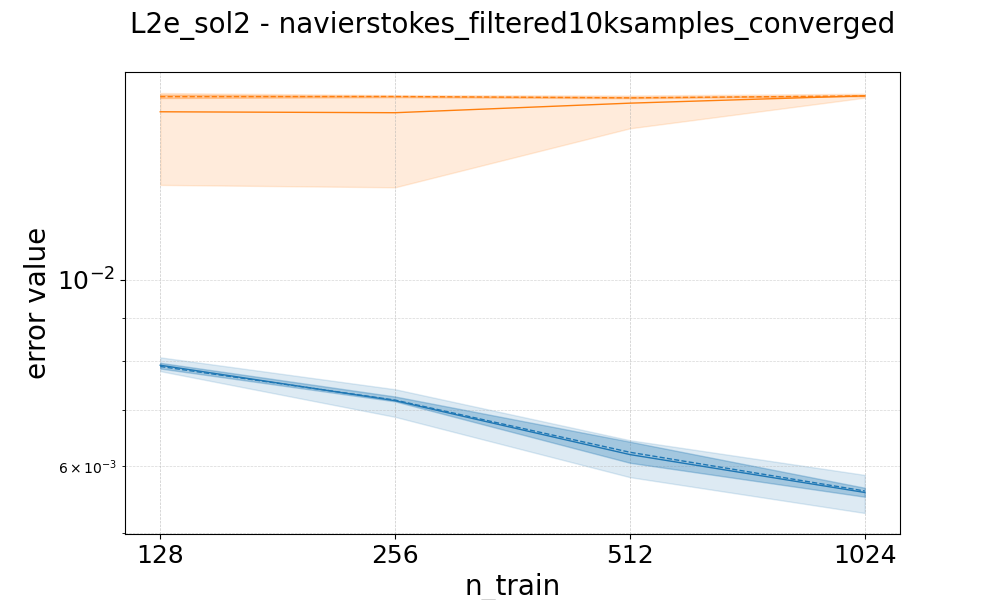}
    }
    \subcaptionbox{$\MLRE$}{\includegraphics[trim=0.5cm 0.0cm 2.5cm 1.8cm,clip,width=0.45\textwidth]{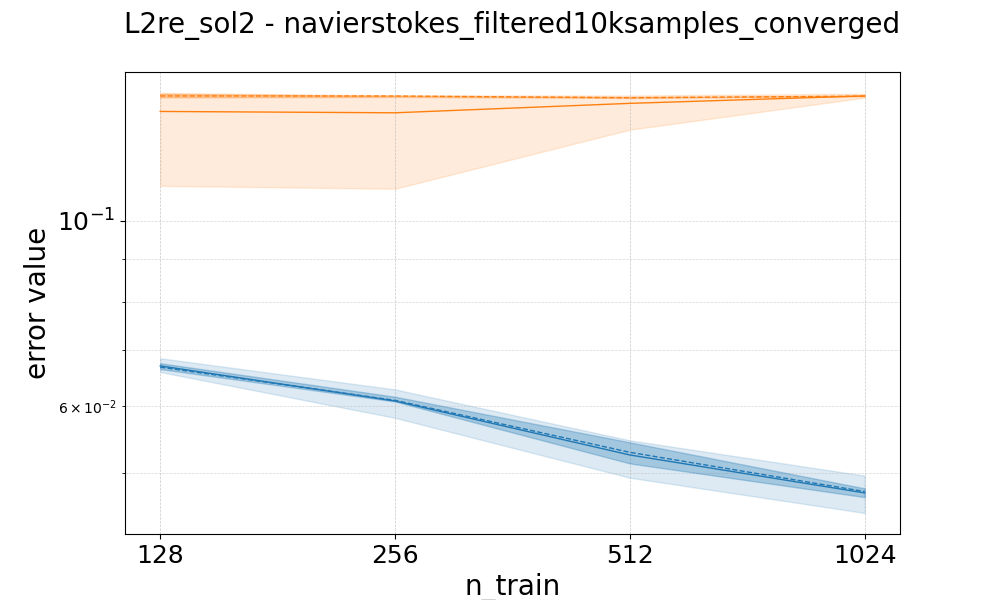}
    }
    \\
    \subcaptionbox{$\MHE$}{\includegraphics[trim=0.5cm 0.0cm 2.5cm 1.8cm,clip,width=0.45\textwidth]{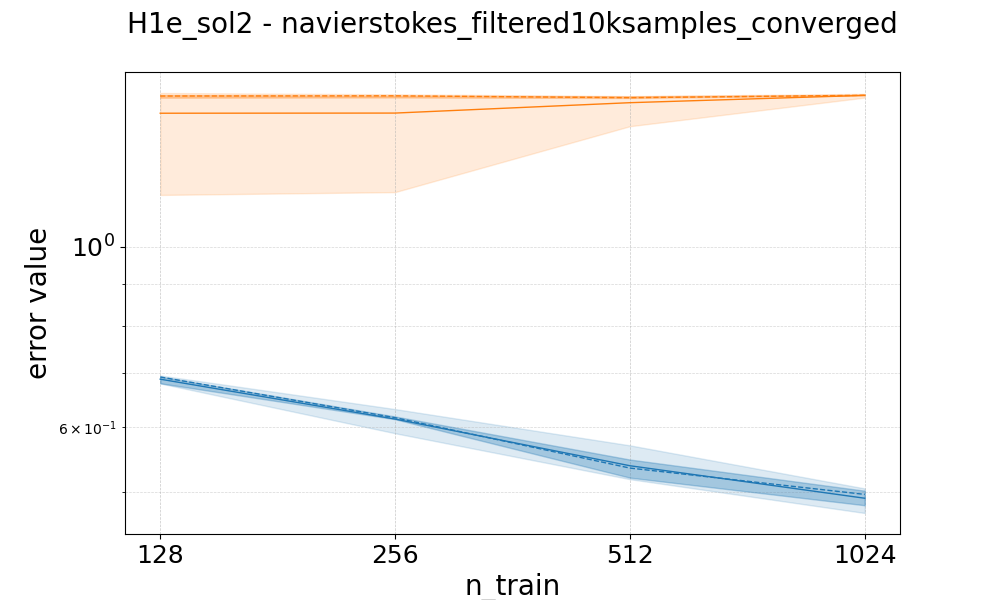}
    }
    \subcaptionbox{$\MHRE$}{\includegraphics[trim=0.5cm 0.0cm 2.5cm 1.8cm,clip,width=0.45\textwidth]{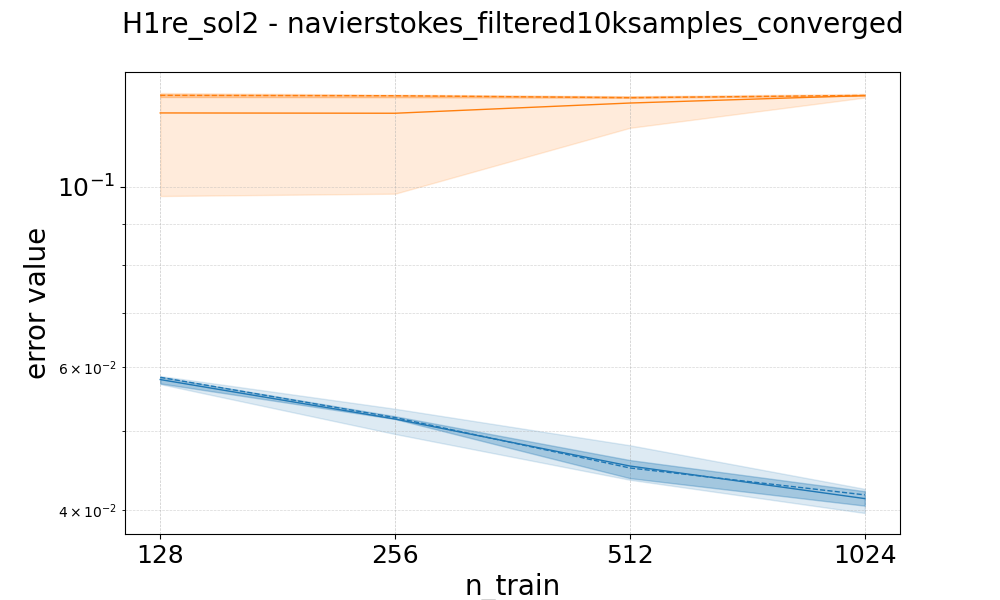}
    }
    \caption{
    Experiment 3 (Navier-Stokes) - Error statistics with respect to five random seeds for weight initialization, varying the training set size. In blue, the \pGINN performance, in orange the \pFCNN performance. Light colored areas represent the Min-Max range of values, the dark colored areas represent the values between fist and third quartiles. The continuous lines represent the average errors (same as Table \ref{tab:results_exp3_sol2}), the dotted lines represent the medians. Results for the second component of solution $\v{u}$.
    }
    \label{fig:results_exp3_sol2}
\end{figure}

\begin{figure}[htbp!]
    \centering
    \subcaptionbox{legend}{\includegraphics[width=0.45\textwidth]{Figures/legend.png}
    }
    \\
    \subcaptionbox{$\MlE$}{\includegraphics[trim=0.5cm 0.0cm 2.5cm 1.8cm,clip,width=0.45\textwidth]{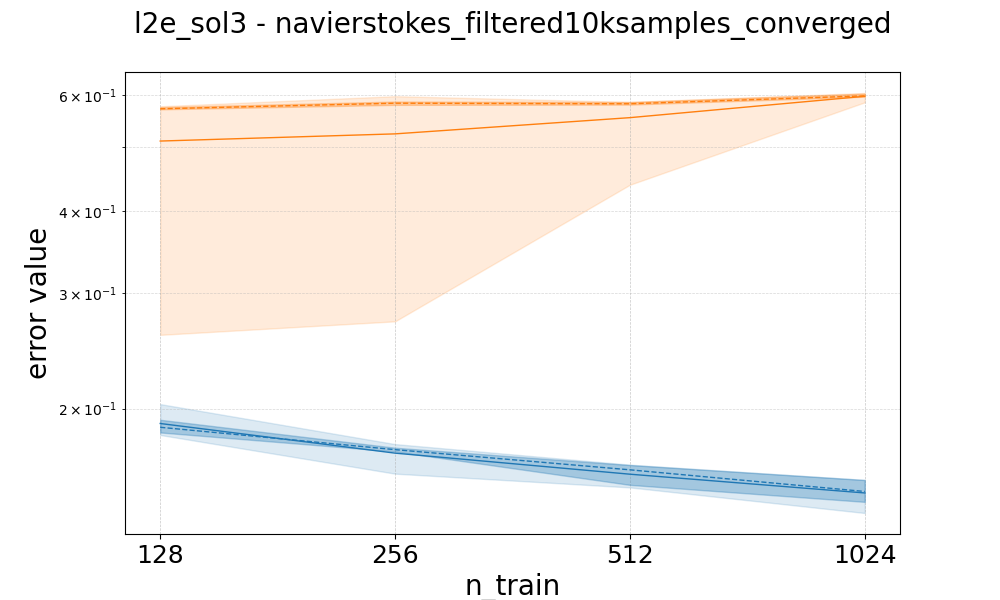}
    }
    \subcaptionbox{$\MlRE$}{\includegraphics[trim=0.5cm 0.0cm 2.5cm 1.8cm,clip,width=0.45\textwidth]{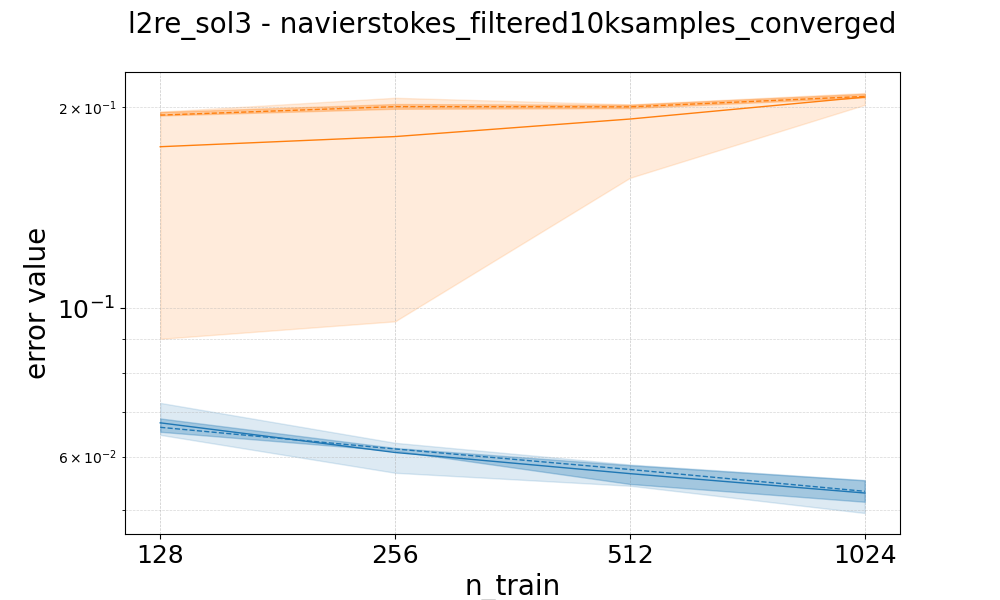}
    }
    \\
	\subcaptionbox{$\MLE$}{\includegraphics[trim=0.5cm 0.0cm 2.5cm 1.8cm,clip,width=0.45\textwidth]{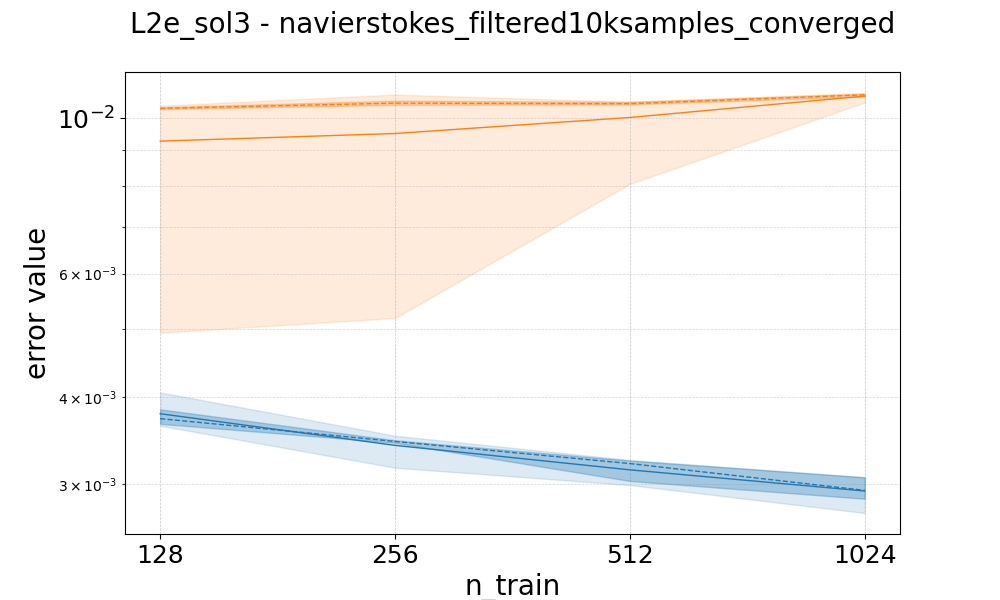}
    }
    \subcaptionbox{$\MLRE$}{\includegraphics[trim=0.5cm 0.0cm 2.5cm 1.8cm,clip,width=0.45\textwidth]{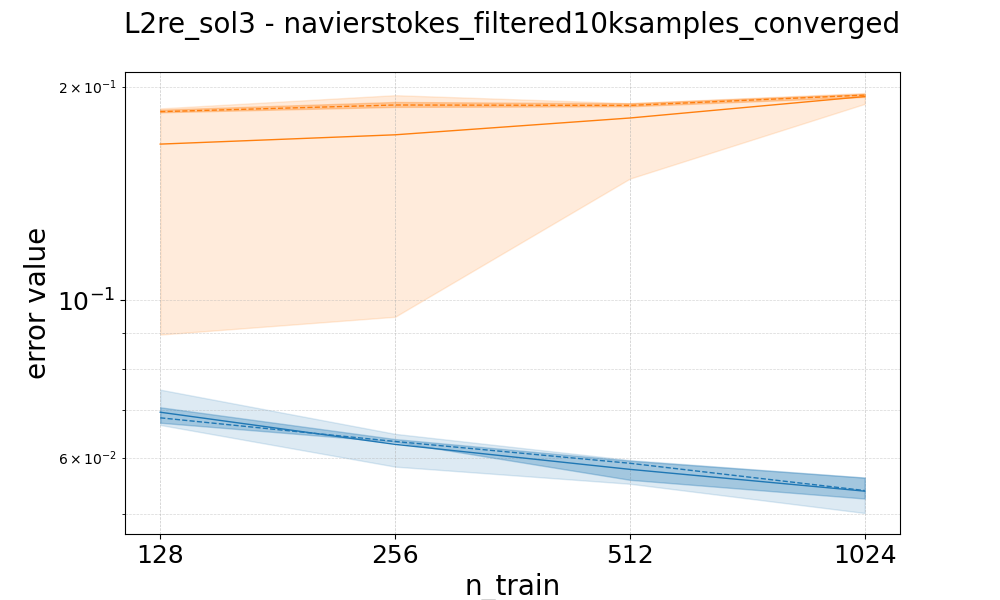}
    }
    \\
    \subcaptionbox{$\MHE$}{\includegraphics[trim=0.5cm 0.0cm 2.5cm 1.8cm,clip,width=0.45\textwidth]{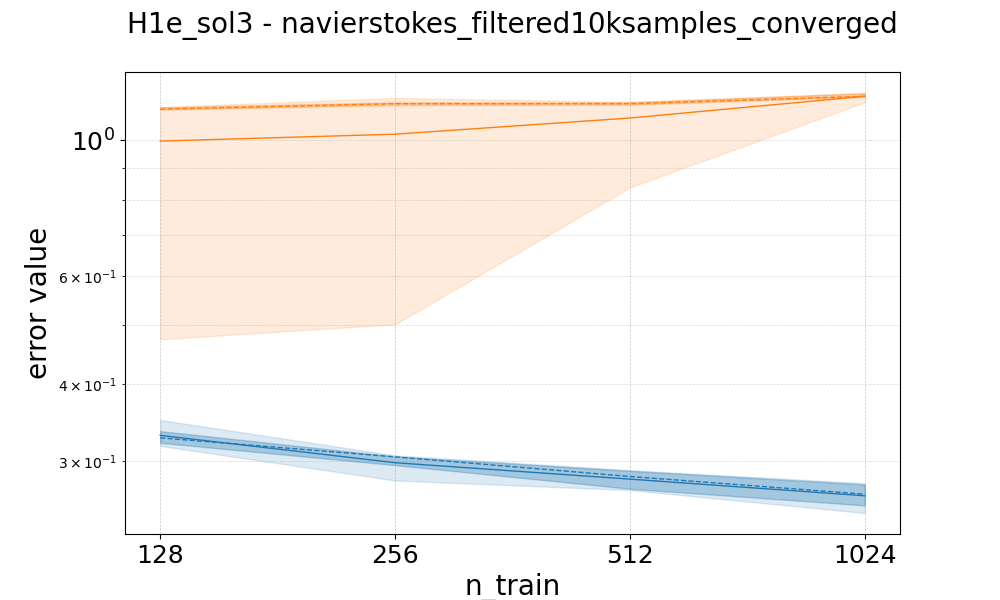}
    }
    \subcaptionbox{$\MHRE$}{\includegraphics[trim=0.5cm 0.0cm 2.5cm 1.8cm,clip,width=0.45\textwidth]{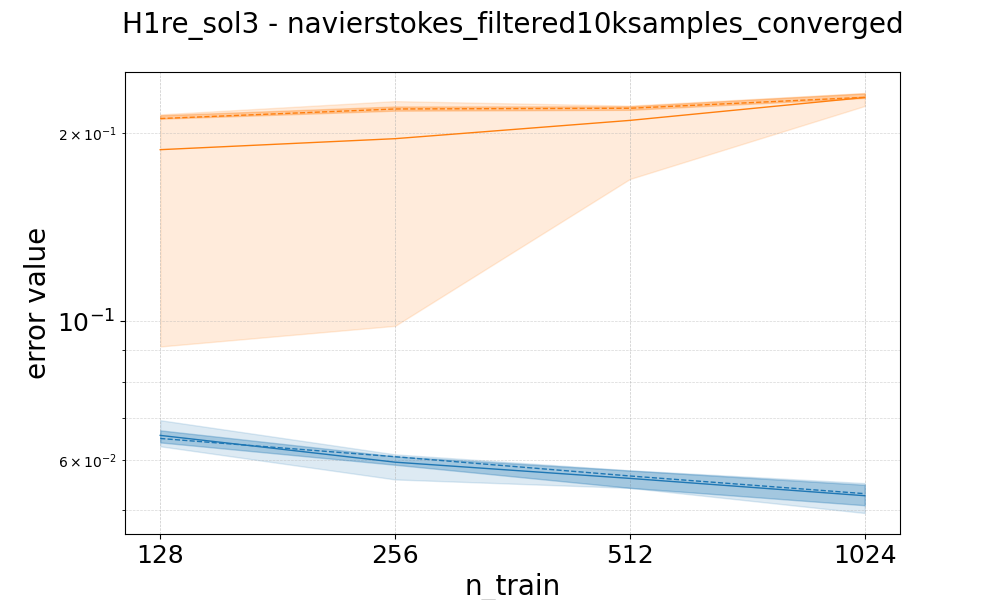}
    }
    \caption{
    Experiment 3 (Navier-Stokes) - Error statistics with respect to five random seeds for weight initialization, varying the training set size. In blue, the \pGINN performance, in orange the \pFCNN performance. Light colored areas represent the Min-Max range of values, the dark colored areas represent the values between fist and third quartiles. The continuous lines represent the average errors (same as Table \ref{tab:results_exp3_sol3}), the dotted lines represent the medians. Results for the pressure solution $p$.
    }
    \label{fig:results_exp3_sol3}
\end{figure}

\begin{figure}[htb!]
    \centering
    \subcaptionbox{Best case - $\solh$}{\includegraphics[trim=1.95cm 4.cm 7cm 5cm,clip,width=0.28\textwidth]{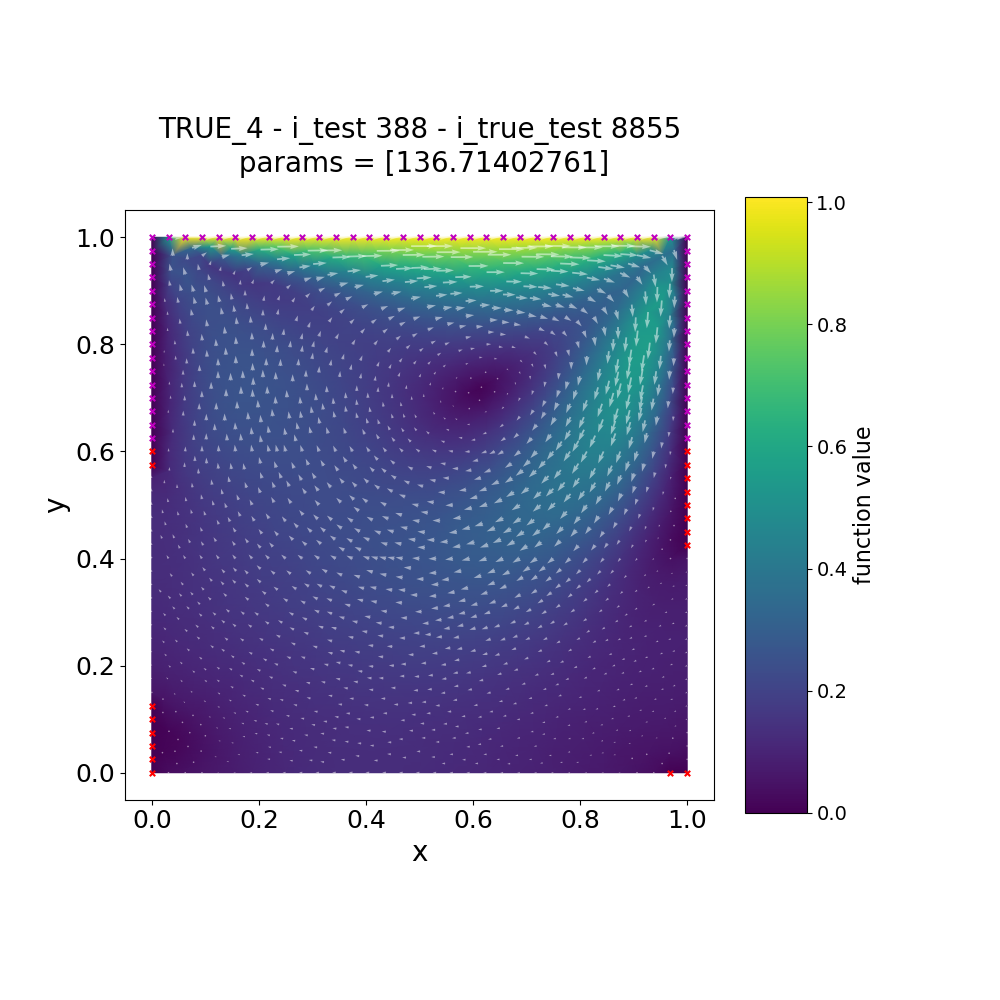}
    }
    \subcaptionbox{Best case - $\nnsolh$}{\includegraphics[trim=1.95cm 4.cm 3.75cm 4.75cm,clip,width=0.335\textwidth]{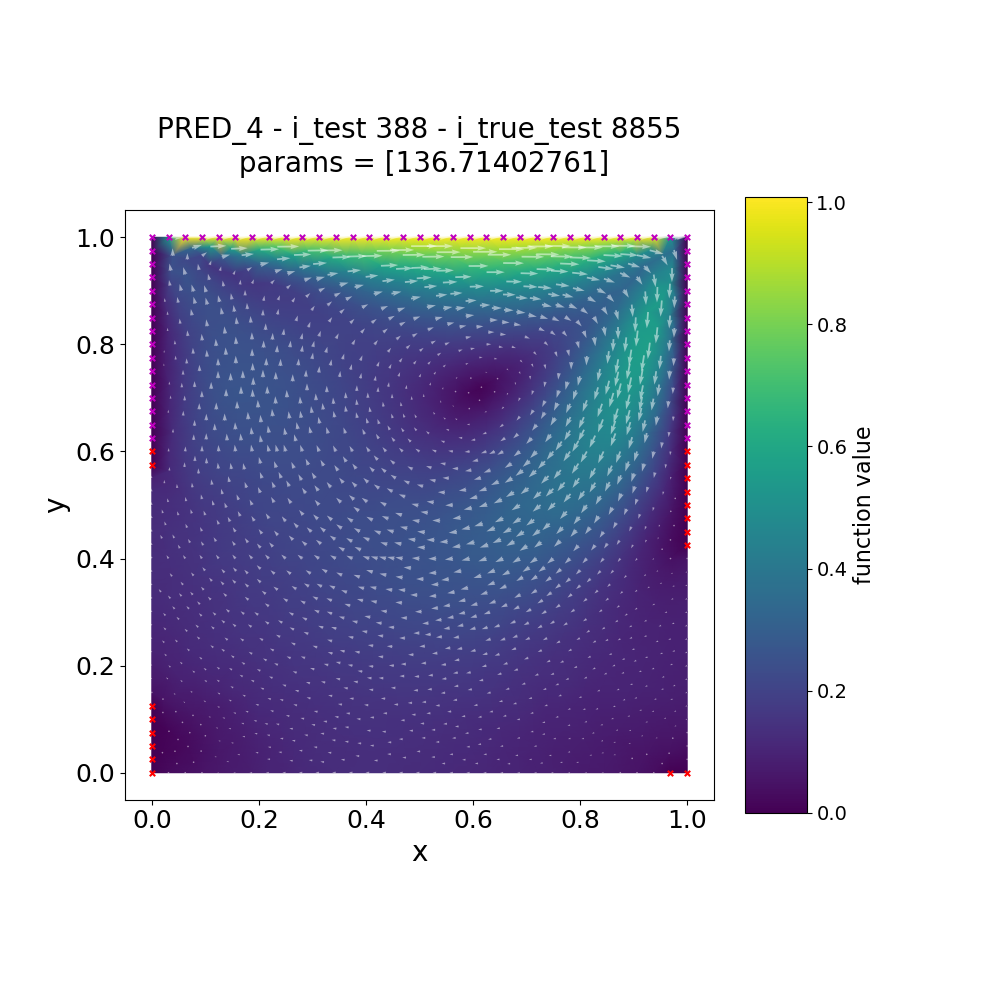}
    }
    \subcaptionbox{Best case - $|\solh - \nnsolh|$}{\includegraphics[trim=1.95cm 4.cm 3.25cm 4.75cm,clip,width=0.345\textwidth]{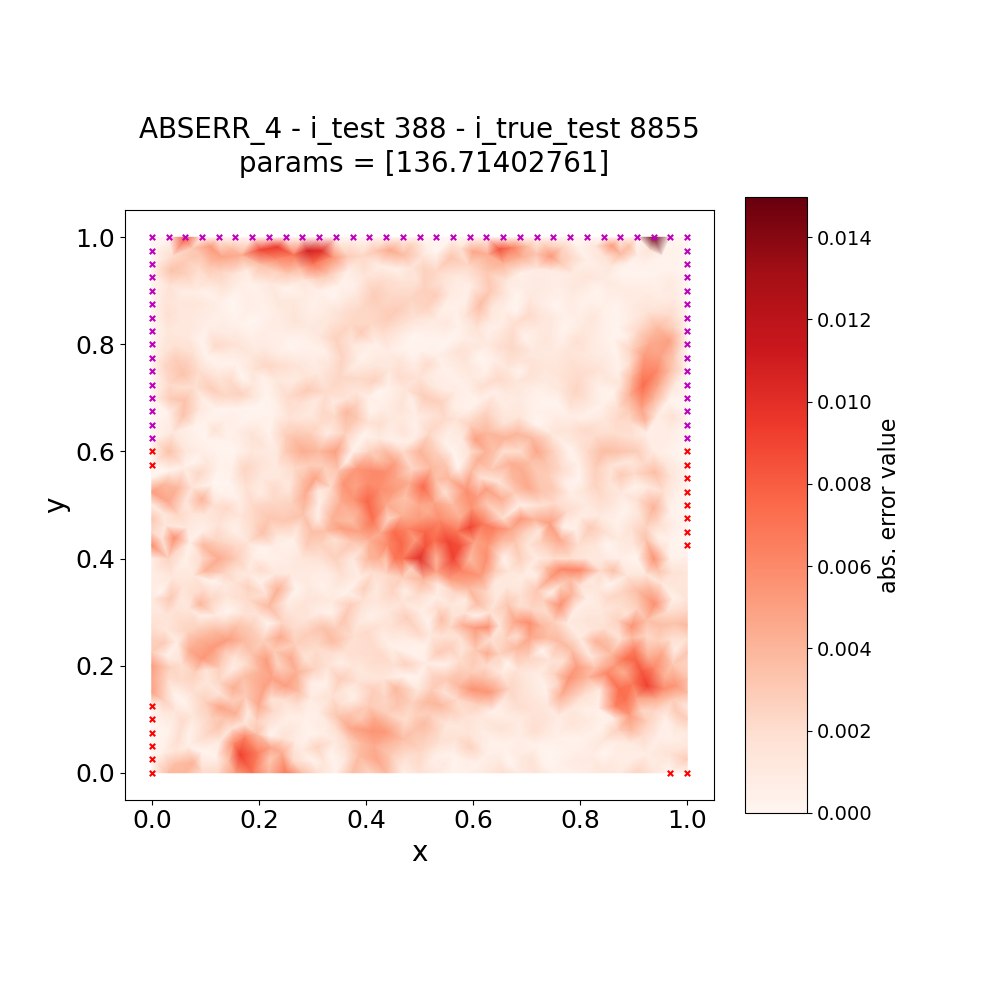}
    }
    \\
    \subcaptionbox{Median case - $\solh$}{\includegraphics[trim=1.95cm 4.cm 7cm 5cm,clip,width=0.28\textwidth]{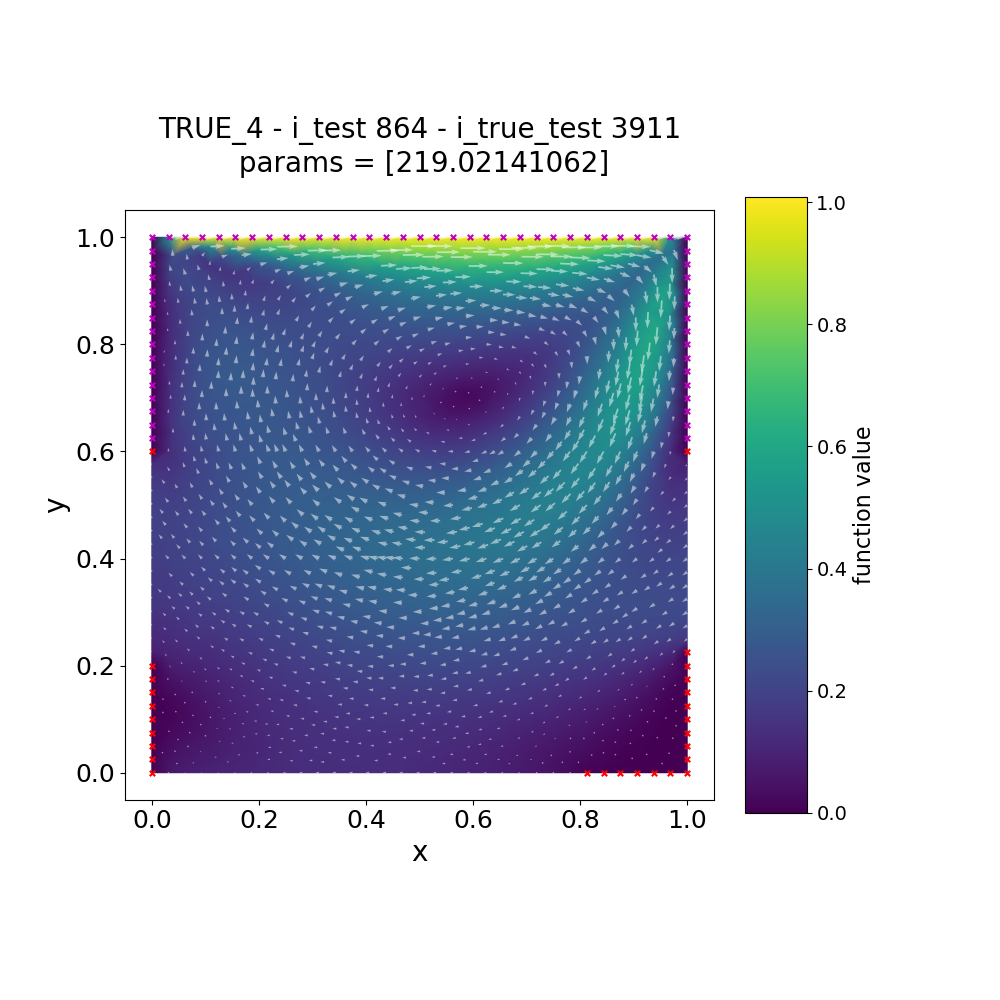}
    }
    \subcaptionbox{Median case - $\nnsolh$}{\includegraphics[trim=1.95cm 4.cm 3.75cm 4.75cm,clip,width=0.335\textwidth]{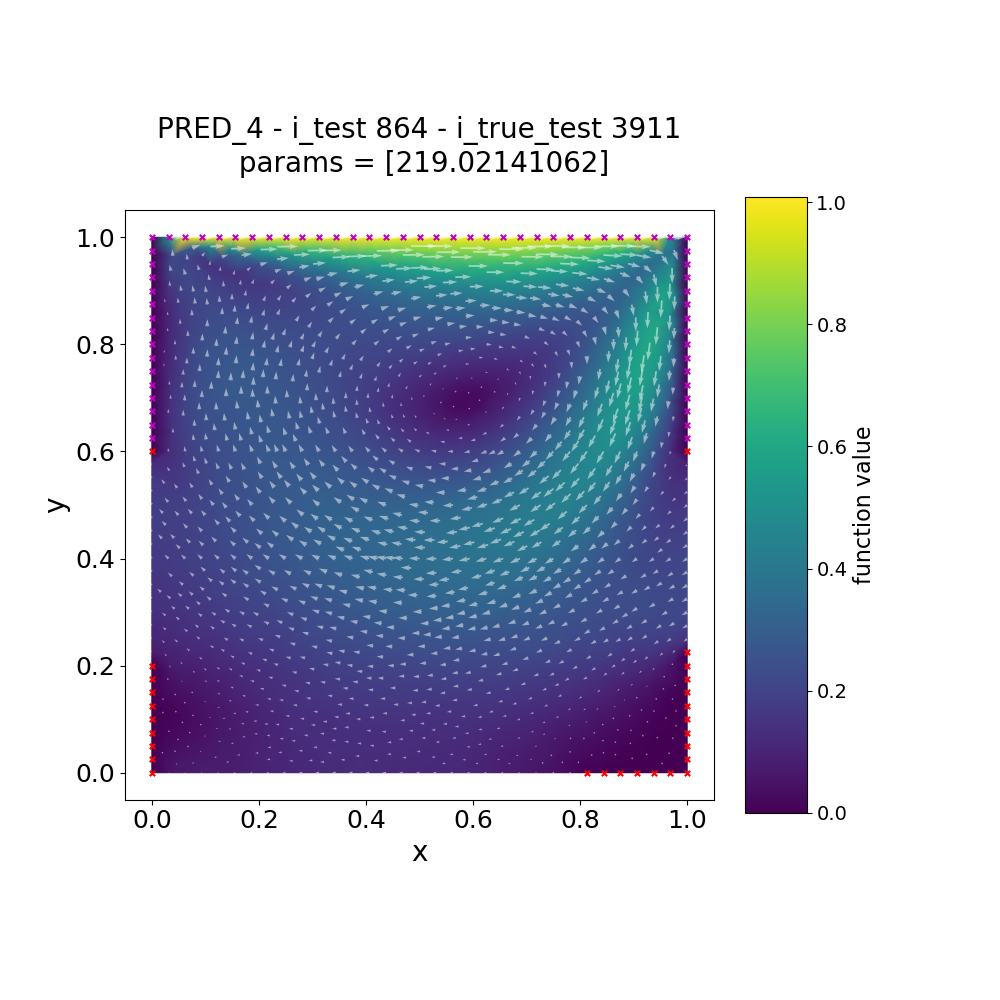}
    }
    \subcaptionbox{Median case - $|\solh - \nnsolh|$}{\includegraphics[trim=1.95cm 4.cm 3.25cm 4.75cm,clip,width=0.345\textwidth]{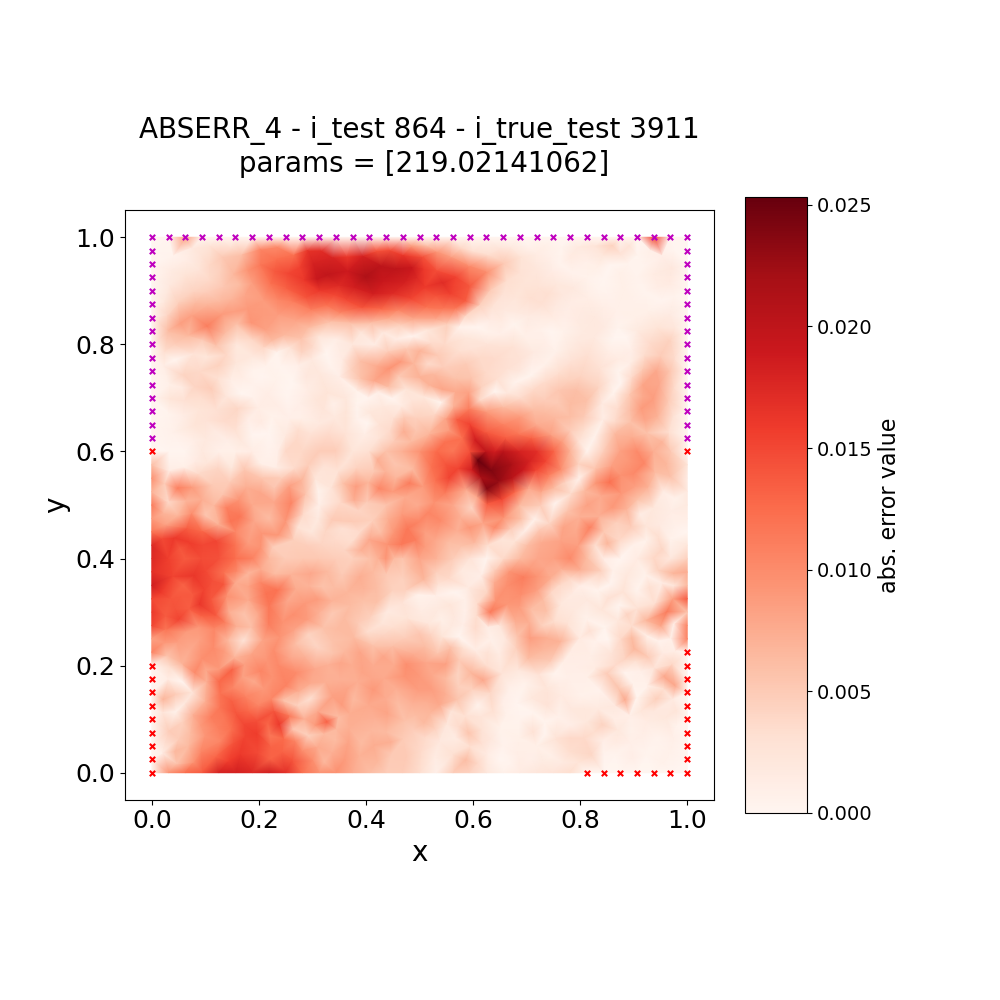}
    }
    \\
    \subcaptionbox{Worst case - $\solh$}{\includegraphics[trim=1.95cm 4.cm 7cm 5cm,clip,width=0.28\textwidth]{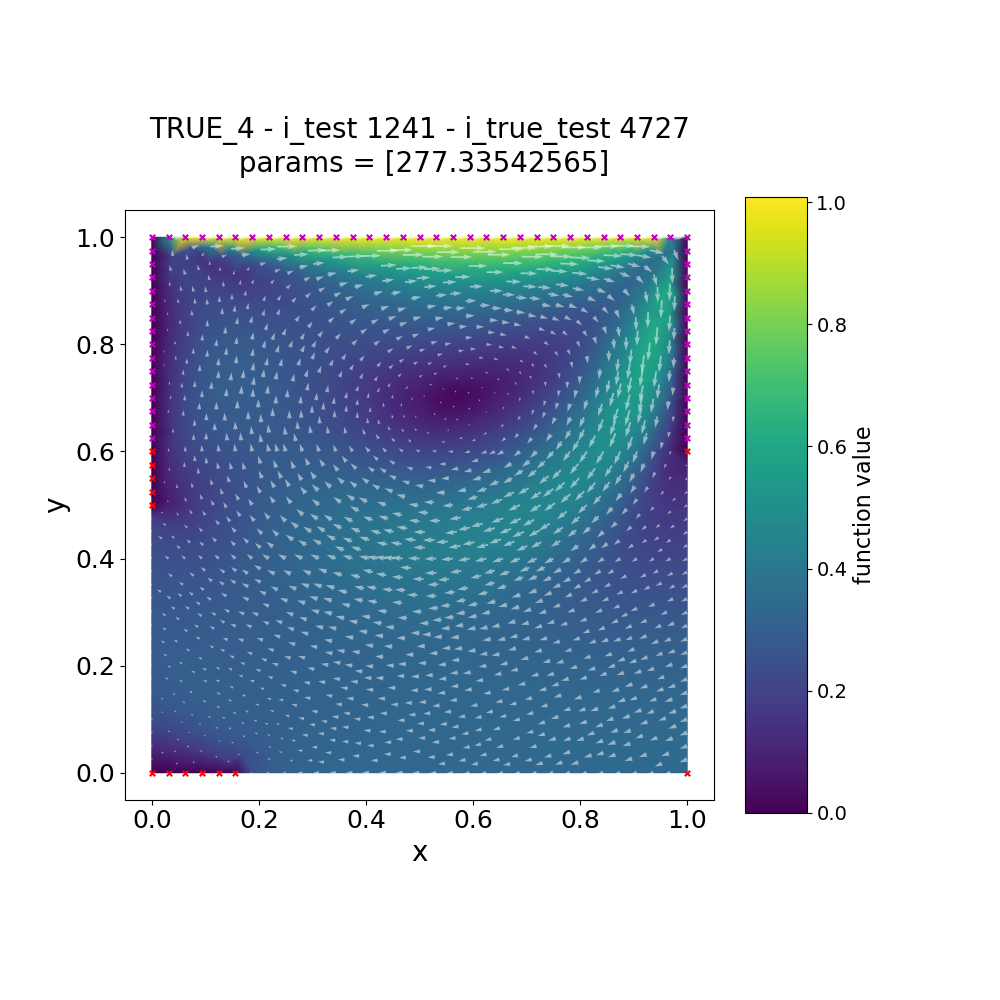}
    }
    \subcaptionbox{Worst case - $\nnsolh$}{\includegraphics[trim=1.95cm 4.cm 3.75cm 4.75cm,clip,width=0.335\textwidth]{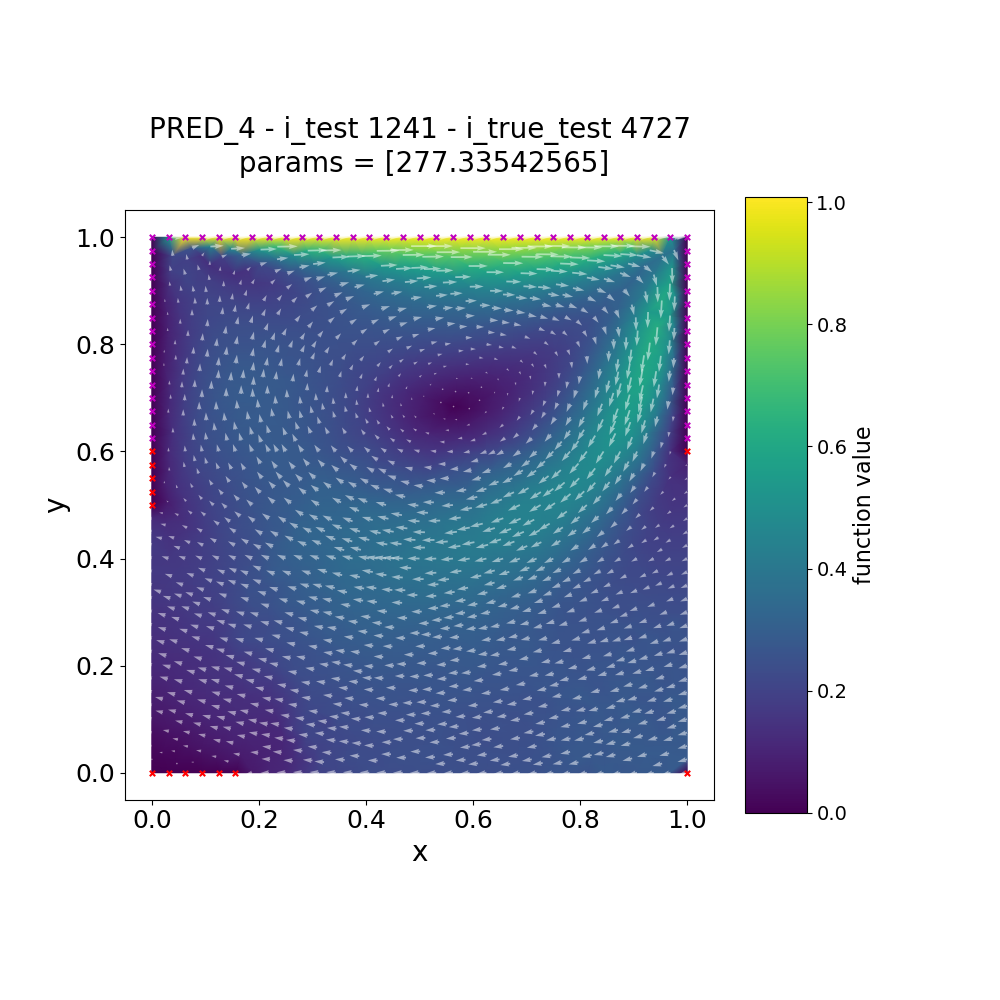}
    }
    \subcaptionbox{Worst case - $|\solh - \nnsolh|$}{\includegraphics[trim=1.95cm 4.cm 3.25cm 4.75cm,clip,width=0.345\textwidth]{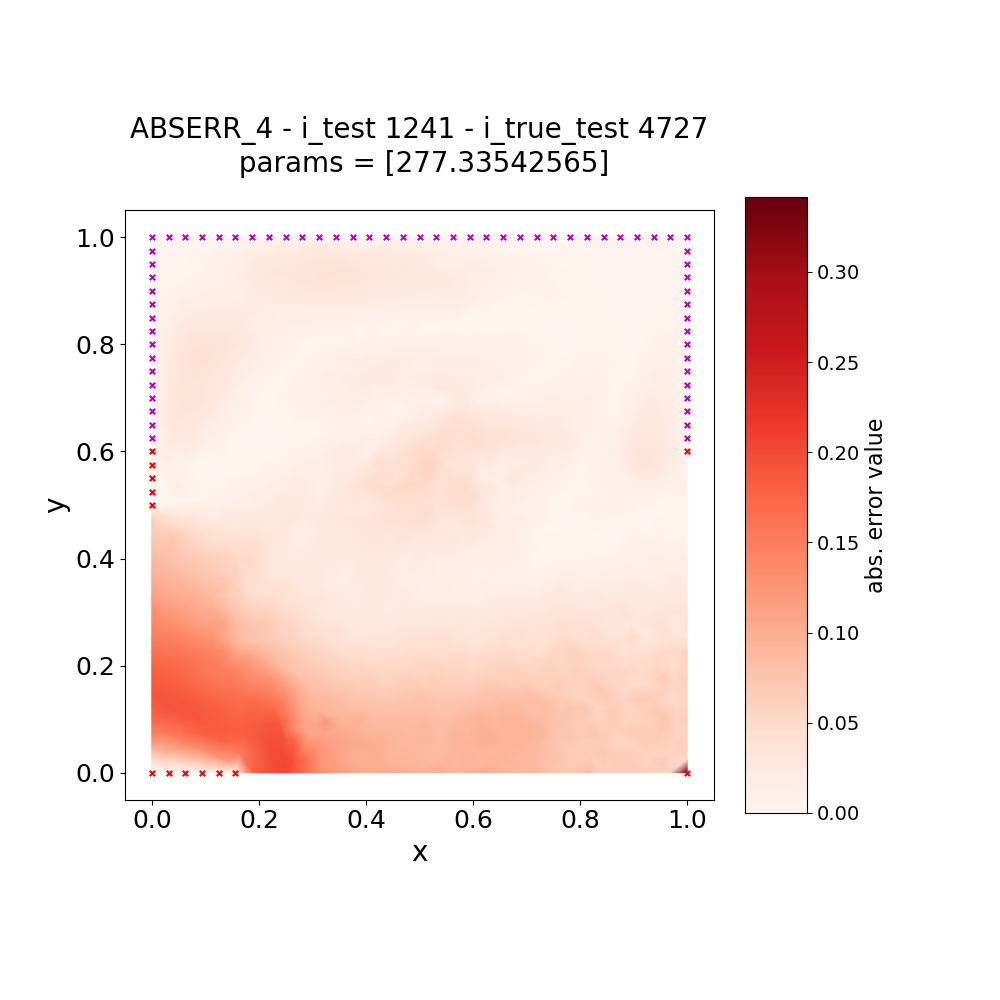}
    }
    \caption{
    Experiment 3 (Navier-Stokes) - Best, median, and worst prediction cases with respect to the test set $\Xi$, ranked with respect to $\MLRE$ averaged on the three solutions. Magenta dots denote fixed homogeneous Dirichlet BCs, while red dots denote the homogeneous Dirichlet BC defined by the boundary parametrization of the problem.
    }
    \label{fig:examples_exp3}
\end{figure}

The model performance is evaluated by considering the average errors illustrated in \eqref{eq:err_mean}, computed for the first and second components of the velocity and the pressure field. These errors, together with information about training epochs and training time, are reported in Tables \ref{tab:results_exp3_sol1}-\ref{tab:results_exp3_sol3} and illustrated in Figures \ref{fig:results_exp3_sol1}-\ref{fig:results_exp3_sol3}. In Figure \ref{fig:examples_exp3}, we report some prediction examples taken from the test set $\Xi$.

% ----------------------------
% SOLUTION 1
% ----------------------------
\begin{table}[htbp]
\centering
\resizebox{0.99\textwidth}{!}{
{
\begin{tabular}{l|l|l||r|r|r||r|r|r||r|r}
Model & n. weights & tr. size ($T$) 
& $\MlE$ & $\MLE$ & $\MHE$
& $\MlRE$ & $\MLRE$ & $\MHRE$
& batch tr.time (s) & tr. epochs \\
\hline
\hline
\multirow{4}{0pt}{\pFCNN} & \multirow{4}{0pt}{3.354e+07} & 128  
& 1.309e+00 & 2.698e-02 & 2.190e+00 
& 1.333e-01 & 1.422e-01 & 1.229e-01 
& 4.711e-01 & 541.4 \\
 & & 256  
& 1.297e+00 & 2.668e-02 & 2.180e+00 
& 1.317e-01 & 1.401e-01 & 1.221e-01 
& 4.688e-01 & 470.4 \\
 & & 512  
& 1.331e+00 & 2.727e-02 & 2.253e+00 
& 1.351e-01 & 1.431e-01 & 1.261e-01 
& 4.693e-01 & 578.2 \\
 & & 1024 
& 1.348e+00 & 2.758e-02 & 2.284e+00 
& 1.369e-01 & 1.449e-01 & 1.280e-01 
& 4.655e-01 & 449.6 \\
\hline
\multirow{4}{0pt}{\pGINN} & \multirow{4}{0pt}{2.655e+07} & 128  
& 7.518e-01 & 1.565e-02 & 1.248e+00 
& 7.639e-02 & 8.232e-02 & 6.987e-02 
& 2.311e+00 & 471.4 \\
 & & 256  
& 6.625e-01 & 1.368e-02 & 1.115e+00 
& 6.727e-02 & 7.188e-02 & 6.241e-02 
& 1.884e+00 & 518.4 \\
 & & 512  
& 5.629e-01 & 1.141e-02 & 9.688e-01 
& 5.705e-02 & 5.979e-02 & 5.418e-02 
& 1.564e+00 & 599.8 \\
 & & 1024 
& 4.926e-01 & 9.875e-03 & 8.579e-01 
& 4.999e-02 & 5.181e-02 & 4.803e-02 
& 1.504e+00 & 489.4 \\
\end{tabular}
}
}
\caption{
{Experiment 3 (Navier-Stokes) - Average performance of the models, with respect to five random seeds for weight initialization, varying the training set size. Results for the first component of solution $\v{u}$.}
}
\label{tab:results_exp3_sol1}
\end{table}

% ----------------------------
% SOLUTION 2
% ----------------------------
\begin{table}[htbp]
\centering
\resizebox{0.99\textwidth}{!}{
{
\begin{tabular}{l|l|l||r|r|r||r|r|r||r|r}
Model & n. weights & tr. size ($T$) 
& $\MlE$ & $\MLE$ & $\MHE$
& $\MlRE$ & $\MLRE$ & $\MHRE$
& batch tr.time (s) & tr. epochs \\
\hline
\hline
\multirow{4}{0pt}{\pFCNN} & \multirow{4}{0pt}{3.354e+07} & 128  
& 8.112e-01 & 1.587e-02 & 1.460e+00 
& 1.288e-01 & 1.350e-01 & 1.234e-01 
& 4.711e-01 & 541.4 \\
 & & 256  
& 8.101e-01 & 1.583e-02 & 1.461e+00 
& 1.284e-01 & 1.345e-01 & 1.232e-01 
& 4.688e-01 & 470.4 \\
 & & 512  
& 8.330e-01 & 1.625e-02 & 1.504e+00 
& 1.320e-01 & 1.380e-01 & 1.269e-01 
& 4.693e-01 & 578.2 \\
 & & 1024 
& 8.502e-01 & 1.658e-02 & 1.536e+00 
& 1.348e-01 & 1.409e-01 & 1.296e-01 
& 4.655e-01 & 449.6 \\
\hline
\multirow{4}{0pt}{\pGINN} & \multirow{4}{0pt}{2.655e+07} & 128  
& 3.949e-01 & 7.914e-03 & 6.878e-01 
& 6.242e-02 & 6.705e-02 & 5.788e-02 
& 2.311e+00 & 471.4 \\
 & & 256  
& 3.562e-01 & 7.180e-03 & 6.143e-01 
& 5.634e-02 & 6.088e-02 & 5.174e-02 
& 1.884e+00 & 518.4 \\
 & & 512  
& 3.093e-01 & 6.193e-03 & 5.383e-01 
& 4.889e-02 & 5.248e-02 & 4.528e-02 
& 1.564e+00 & 599.8 \\
 & & 1024 
& 2.804e-01 & 5.581e-03 & 4.910e-01 
& 4.430e-02 & 4.728e-02 & 4.128e-02 
& 1.504e+00 & 489.4 \\
\end{tabular}
}
}
\caption{
{Experiment 3 (Navier-Stokes) - Average performance of the models, with respect to five random seeds for weight initialization, varying the training set size. Results for the second component of solution $\v{u}$.}
}
\label{tab:results_exp3_sol2}
\end{table}

% ----------------------------
% SOLUTION 3
% ----------------------------
\begin{table}[htbp]
\centering
\resizebox{0.99\textwidth}{!}{
{
\begin{tabular}{l|l|l||r|r|r||r|r|r||r|r}
Model & n. weights & tr. size ($T$) 
& $\MlE$ & $\MLE$ & $\MHE$
& $\MlRE$ & $\MLRE$ & $\MHRE$
& batch tr.time (s) & tr. epochs \\
\hline
\hline
\multirow{4}{0pt}{\pFCNN} & \multirow{4}{0pt}{3.354e+07} & 128  
& 5.109e-01 & 9.276e-03 & 9.968e-01 
& 1.744e-01 & 1.659e-01 & 1.881e-01 
& 4.711e-01 & 541.4 \\
 & & 256  
& 5.239e-01 & 9.512e-03 & 1.023e+00 
& 1.806e-01 & 1.710e-01 & 1.958e-01 
& 4.688e-01 & 470.4 \\
 & & 512  
& 5.545e-01 & 1.003e-02 & 1.087e+00 
& 1.919e-01 & 1.806e-01 & 2.094e-01 
& 4.693e-01 & 578.2 \\
 & & 1024 
& 5.976e-01 & 1.076e-02 & 1.179e+00 
& 2.070e-01 & 1.936e-01 & 2.278e-01 
& 4.655e-01 & 449.6 \\
\hline
\multirow{4}{0pt}{\pGINN} & \multirow{4}{0pt}{2.655e+07} & 128  
& 1.899e-01 & 3.784e-03 & 3.303e-01 
& 6.741e-02 & 6.950e-02 & 6.575e-02 
& 2.311e+00 & 471.4 \\
 & & 256  
& 1.711e-01 & 3.409e-03 & 2.981e-01 
& 6.087e-02 & 6.262e-02 & 5.959e-02 
& 1.884e+00 & 518.4 \\
 & & 512  
& 1.590e-01 & 3.146e-03 & 2.803e-01 
& 5.659e-02 & 5.774e-02 & 5.613e-02 
& 1.564e+00 & 599.8 \\
 & & 1024 
& 1.489e-01 & 2.934e-03 & 2.633e-01 
& 5.294e-02 & 5.380e-02 & 5.267e-02 
& 1.504e+00 & 489.4 \\
\end{tabular}
}
}
\caption{
{Experiment 3 (Navier-Stokes) - Average performance of the models, with respect to five random seeds for weight initialization, varying the training set size. Results for the pressure solution $p$.}
}
\label{tab:results_exp3_sol3}
\end{table}

In this last, third experiment, the results we obtain are similar to the ones obtained for the first experiment. Specifically, we still have the evidence of the advantage in using \pGINNs instead of \pFCNNs, by looking at the behaviors and values of the average errors of the trained models. The \pGINN models have better predictive performance on the test set for both $\v{u}$ and the pressure $p$. As in the previous experiments, this performance is stable with respect to the random initialization of the weights. In contrast, \pFCNNs only rarely achieve slightly better results for some initializations (and still not better than those of \pGINNs).

Like in Experiment 1, \pFCNNs suffer from poor generalization properties for all the training set sizes, with an average value of the errors that is almost constant. On the contrary, as already noticed in the previous experiments, \pGINNs exhibit a clear error reduction as the training set size increases and low error values even with few hundred training samples.

Regarding computational costs, in this case \pGINNs have the same order of magnitude in the number of trainable parameters as \pFCNNs, although still slightly fewer (2.655e+07 weights instead of 3.357e+07). On the other hand, the training time of \pGINNs is longer (approximately one order of magnitude greater).
This similarity in computational costs between Experiment 1 and Experiment 3 is due to the similar number of mesh nodes $N_h$ (1141 in Experiment 1 and 1444 in Experiment 3).

For more details about the scalability of the models, see Section \ref{sec:gen_comments_exp} below.

\subsection{{General Comments for Experiment Results}}\label{sec:gen_comments_exp}

Overall, across the three experiments, a clear and consistent trend emerges in favor of \pGINNs over \pFCNNs. Concerning the aspects of predictive performance, \pGINNs systematically achieve lower average errors on the test set, for all the errors considered. Moreover, their performance is remarkably stable with respect to the random initialization of the weights, indicating a more robust and reliable training process. In contrast, \pFCNNs exhibit higher variability across different initializations and only exceptionally attain slightly improved results, which, however, never surpass those obtained with \pGINNs.

A second key aspect concerns generalization with respect to the training set size. In all experiments, \pFCNNs display poor scalability in terms of accuracy: the average error remains almost constant as the number of training samples increases, highlighting limited generalization capabilities. On the contrary, \pGINNs consistently show a clear error decay as the training set size grows, together with satisfactory accuracy levels even in low-data regimes (i.e., when only a few hundred training samples are available). This behavior suggests that the graph-instructed architecture is significantly more data-efficient and better suited to capturing the underlying structure of the problem.

\begin{figure}[htb]
    \centering
    \includegraphics[trim=1.5cm 0cm 2.5cm 1.8cm,clip,width=0.75\textwidth]{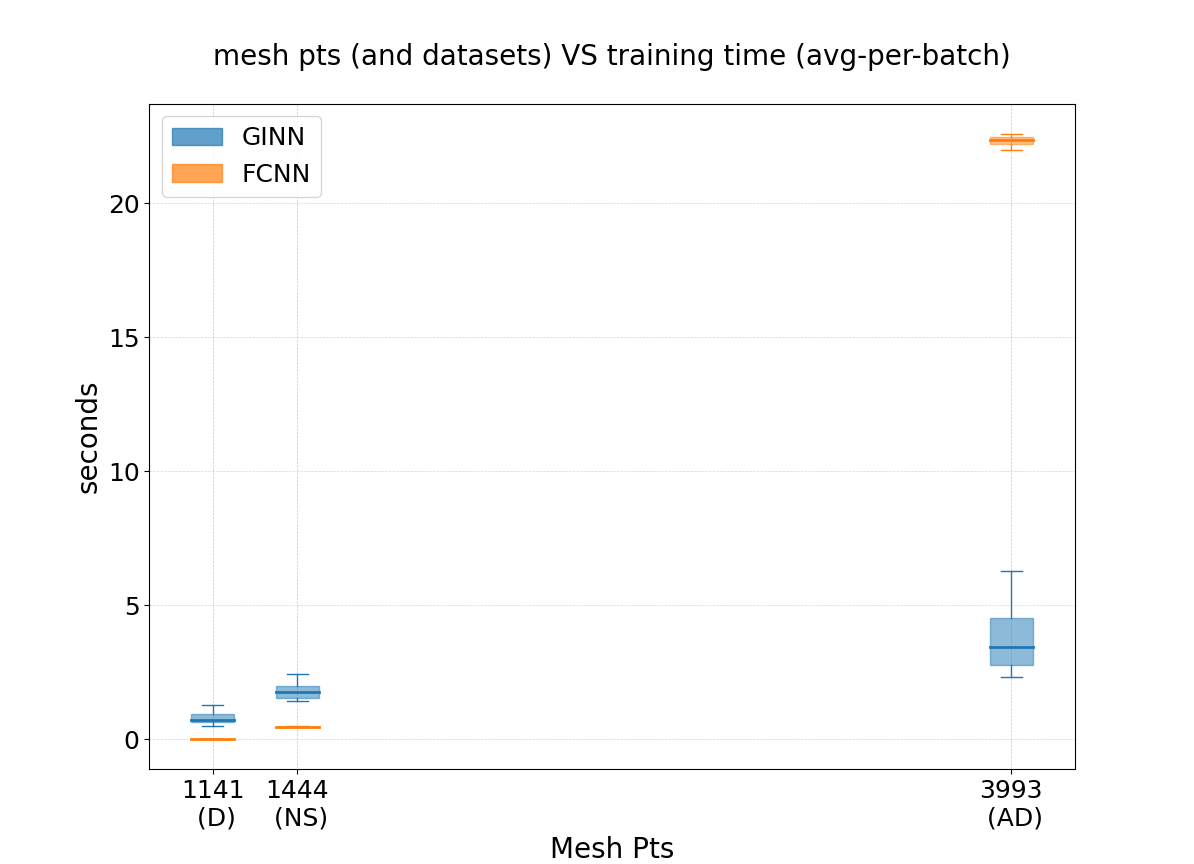}
    \caption{{Statistics of the average training time per mini-batch (in seconds), varying the experiments (Experiment 1, 2 and 3 are denoted by (D), (AD), and (NS), respectively); i.e., varying the number $N_h$ of mesh nodes. In blue the training times of \pGINNs, in orange the training times of \pFCNNs.}}
    \label{fig:trainingtimes}
\end{figure}

From a computational perspective, the comparison is more nuanced and strongly influenced by the number of mesh nodes $N_h$. In terms of trainable parameters, \pGINNs always require fewer weights than \pFCNNs, with differences ranging from moderate to one order of magnitude, depending on the experiment. However, training times seem not to not follow a uniform trend: for smaller meshes, \pGINNs may require longer training times despite having fewer parameters, whereas for larger meshes they become significantly more efficient. Actually, this behavior reflects the different scaling properties of FC layers and GI layers with respect to $N_h$, with the FC layers exhibiting a markedly worse scalability both in parameter count and computational cost. See Figure \ref{fig:trainingtimes} for better understanding the scalability of the training times.

Taken together, these results indicate that \pGINNs provide a more robust, scalable, and data-efficient modeling strategy, particularly in regimes characterized by increasing mesh resolution or limited training data availability.

%%%%%%%%%%%%%%%%%%%%%%%%%%%%%%%%%%%%%%%%%%%%%%%%%%%%%%%%%%%%%%%%%%%%%%%%%%%%%%%%%%%%%%%%%%%

\section{Conclusions}
\label{sec:conc}

In this work, we propose a novel methodological setting based on GINNs to tackle parametric varying PDEs, the \pGINNs, where the parameter acts on both the physics of the problem and the imposition of the BCs, surpassing the classical surrogate and reduced models proposed in the literature. We remark that the varying boundary problem has been limitedly investigated, both in the ROM and surrogate model communities. For the first time, an architecture related to the local and sparse information of the mesh discretization is employed, establishing a real-time approach that maps the parametric information to the PDE solution.

Our numerical experiments demonstrate that the \pGINNs are an accurate and efficient surrogate model for these types of problems, thanks to a thorough validation in several settings based on linear and nonlinear PDEs. Moreover, the GINN-based framework demonstrates more robust and accurate results compared to its fully connected counterpart, the \pFCNNs. This work represents a first methodological step and its results suggest that \pGINNs are a promising asset for many-query scenarios.

A natural step forward in our research is the extension of \pGINNs to more challenging problems where fast-response simulations are needed, such as design and control. Another interesting research direction is related to varying geometries and complex topologies. We believe that this work paves the way for a more sustainable and efficient integration of mathematical models directly in simulation science and its application.

%%%%%%%%%%%%%%%%%%%%%%%%%%%%%%%%%%%%%%%%%%%%%%%%%%%%%%%%%%%%%%%%%%%%%%%%%%%%%%%%%%%%%%%%%%%

\section*{Acknowledgments}
This study was carried out within the ``20227K44ME - Full and Reduced order modelling of coupled systems: focus on non-matching methods and automatic learning (FaReX)" project – funded by European Union – Next Generation EU  within the PRIN 2022 program (D.D. 104 - 02/02/2022 Ministero dell’Università e della Ricerca). This manuscript reflects only the authors’ views and opinions, and the Ministry cannot be considered responsible for them. Moreover, MS thanks the INdAM-GNCS Project ``Metodi numerici efficienti per problemi accoppiati in sistemi complessi" (CUP E53C24001950001).
F.D. acknowledges that this study was carried out within the FAIR-Future Artificial Intelligence Research and received funding from the European Union Next-GenerationEU (PIANO NAZIONALE DI RIPRESA E RESILIENZA (PNRR) – MISSIONE 4 COMPONENTE 2, INVESTIMENTO 1.3---D.D. 1555 11/10/2022, PE00000013). This manuscript reflects only the authors’ views and opinions; neither the European Union nor the European Commission can be considered responsible for them.

% \appendix

%%%%%%%%%%%%%%%%%%%%%%%%%%%%%%%%%%%%%%%%%%%%%%%%%%%%%%%%%%%%%%%%%%%%%%%%%%%%%%%%%%%%%%%%%%%

 \bibliographystyle{elsarticle-num} 
 \bibliography{references,references_fds}





% \end{thebibliography}
\end{document}